\documentclass[letterpaper,11pt,oneside,reqno]{article}
\usepackage{bm}
\usepackage{mathrsfs}
\usepackage{amsfonts,amsmath, amssymb,amsthm,amscd,stmaryrd}
\usepackage[height=9.6in,width=5.95in]{geometry}
\usepackage[mpexclude,DIV13]{typearea}
\usepackage{verbatim}
\usepackage{hyperref}
\usepackage{graphicx,xcolor}
\usepackage[latin1]{inputenc}
\usepackage{latexsym}
\usepackage{lscape}
  \ifx\LabelFigloaded\MYundefined\relax
  \else
   \fi

  \def\LabelFigloaded{\relax}

  \chardef\LabelFigCatAt\the\catcode`\@
  \catcode`\@=11

 \let\LabelFigwlog@ld\wlog
 \def\wlog#1{\relax}

 \ifx\\\MYundefined@
    \let\\\relax
 \fi

  \def\ms@g{\immediate\write16}

 \def\N@wif{\csname newif\endcsname }
 \def\Temp@ {\N@wif\ifIN@}
 \ifx\INN@\MYundefined@
    \else \let\Temp@\relax
 \fi
 \Temp@

  \def\IN@{\expandafter\INN@\expandafter}
  \long\def\INN@0#1@#2@{\long\def\NI@##1#1##2##3\ENDNI@
    {\ifx\m@rker##2\IN@false\else\IN@true\fi}%
     \expandafter\NI@#2@@#1\m@rker\ENDNI@}
  \def\m@rker{\m@@rker}
 
  \newtoks\Initialtoks@  \newtoks\Terminaltoks@
  \def\SPLIT@{\expandafter\SPLITT@\expandafter}
  \def\SPLITT@0#1@#2@{\def\TTILPS@##1#1##2@{%
     \Initialtoks@{##1}\Terminaltoks@{##2}}\expandafter\TTILPS@#2@}

 \def\Shifted@@#1#2#3{\setbox0=\hbox{#3}%
   \raise -\dp0\vbox {\kern-#2%
       \hbox {\kern#1\unhbox0\kern-#1}%
           \kern#2}}

 \newcount\gridcount
 \newbox\auxGridbox@ \newbox\hGridbox@ \newbox\vGridbox@
 \newbox\Labelbox@ \newbox\auxLabelbox@
 \newbox\Coordinatebox@
 \newtoks\Labeltoks@
 \newdimen\Wdd@ \newdimen\Htt@
 \newdimen\Wddd@ \newdimen\Httt@
 
 \def\Wr@{\immediate\write16}

 \newdimen\GL@wd
 \def\GridLineWidth#1{\GL@wd=#1}

 \def\gobble#1{}
 \def\EdgeErr@{\Wr@{}%
      \Wr@{\string\Edges\space argument
      1, 10, 100 or 1000 please\string!}%
      }

 \newcount\Edgect@

 \def\Sweepup#1\endSweepup{}

 \def\SetEdges@{%
    \edef\Zr@@s{\expandafter\gobble\number\Edgect@\empty}%
        \count255=0\Zr@@s\relax
        \ifnum\count255=\z@\else\EdgeErr@\show\tailtest\fi
        \count255=1\Zr@@s\relax
        \ifnum\count255=\Edgect@\relax\else\EdgeErr@\show\leadtest\fi
    \EdgGl@b\edef\Zr@s{\expandafter\gobble\Zr@@s\empty}
    \ifnum\Edgect@>\@ne\relax\EdgGl@b\let\L@Dc\empty
        \else\EdgGl@b\edef\L@Dc{\string.}\fi
    \ifnum\Edgect@>\@ne\relax
        \EdgGl@b\edef\Edgescale@##1{\divide##1 by \Edgect@}%
        \else\EdgGl@b\edef\Edgescale@##1{}\fi
    }

 \def\Edges#1{\Edgect@=#1\relax
     \let\EdgGl@b\global \SetEdges@}

 \Edges{1}

 \def\hhrule{\hrule height \GL@wd\vskip-.\GL@wd}

 \def\hRule@{%
   \advance\gridcount -2%
   \vfil\hhrule\vfil
   \llap{\smash{\raise -2.5pt
     \hbox{\L@Dc\number\gridcount\Zr@s\kern2pt}}}%
   \hhrule
   }

\def\vvrule{\vrule width \GL@wd \kern-\GL@wd}

 \def\vRule@{\advance\gridcount 2%
   \hfil\vvrule\hfil
   \setbox\auxGridbox@=\vbox to 0pt
      {\vskip \Htt@\vskip 2pt
        \hbox to 0pt{\hss\L@Dc\number\gridcount\Zr@s\hss}\vss}%
      \wd\auxGridbox@=0pt \box\auxGridbox@
   \vvrule
   }

 \def\PlaceGrid@@{\gridcount=10 
  \setbox\hGridbox@=\hbox{%
        \hbox{%
             \hskip-.4pt\vrule
             \vbox to \Htt@{%
               \offinterlineskip\parindent=\z@\relax
               \hbox to \Wdd@{\hfil}
               \hRule@\hRule@\hRule@\hRule@
               \vfil\hhrule\vfil}%
             \vrule\hskip-.4pt}
    }%
  \gridcount=0%
  \setbox\vGridbox@=\hbox{%
      \vbox{\offinterlineskip\parindent=0pt\hsize=0pt
         \vskip-.4pt\hrule%
         \hbox to \Wdd@{%
                 \vtop to \Htt@{\vfil}%
                 \vRule@\vRule@\vRule@\vRule@
                 \hfil\vvrule\hfil}%
         \hrule\vskip-.4pt}}%
  \wd\hGridbox@=0pt\ht\hGridbox@=0pt
  \wd\vGridbox@=0pt\ht\vGridbox@=0pt
  \hbox{\box\hGridbox@\box\vGridbox@}%
  }

 \def\LabelsGlobal{\def\LabGl@b{\global}}
 \def\LabelsLocal{\def\LabGl@b{}}
 \LabelsGlobal 

 \def\SetLabels#1\endSetLabels{%
   \LabGl@b\Labeltoks@={#1()\\}%
   }

 \LabGl@b\Labeltoks@={()\\}

 \def\ShowGrid{\LabGl@b\let\PlaceGrid@\PlaceGrid@@}
 \def\HideGrid{\LabGl@b\let\PlaceGrid@\relax}
 \def\Grids{\ShowGrid\LabGl@b\let\GridSwitch@\ShowGrid}
 \def\noGrids{\HideGrid\LabGl@b\let\GridSwitch@\HideGrid}

 \noGrids

 \def\bAdjust@@{%
     \setbox\auxLabelbox@=\hbox{\raise \dp\auxLabelbox@
            \box\auxLabelbox@}}
 \def\bAdjust@{\let\vAdjust@\bAdjust@@}

 \def\eAdjust@@{\dimen0=-.5\ht\auxLabelbox@
     \advance\dimen0 by .5\dp\auxLabelbox@
     \setbox\auxLabelbox@=
            \hbox{\raise\dimen0\box\auxLabelbox@}}
 \def\eAdjust@{\let\vAdjust@\eAdjust@@}

 \def\tAdjust@@{%
     \setbox\auxLabelbox@=\hbox{\raise-\ht\auxLabelbox@
            \box\auxLabelbox@}}
 \def\tAdjust@{\let\vAdjust@\tAdjust@@}

 \let\vAdjust@\relax

 \def\lAdjust@{\let\hAdjust@\rlap}
 \def\rAdjust@{\let\hAdjust@\llap}

 \let\hAdjust@\relax\let\vAdjust@\relax

 \def\FetchLabel@#1(#2)#3\\{%
     \IN@0#2@@\ifIN@
        \setbox0=\hbox{\ignorespaces#1#3\unskip}%
        \ifdim\wd0>0pt
           \ms@g{}%
           \ms@g{ !!! Bad label(s)? !!!}%
           \message{ #1(#2)#3}%
        \fi
        \def\LabelMole@##1\endFetchLabel@{%
            \IN@0()\\@##1@%
            \ifIN@\def\Temp@{\FetchLabel@##1\endFetchLabel@}%
            \else\def\Temp@{}%
            \fi
            \Temp@
           }%
     \else
       \ignorespaces#1\unskip
       \setbox\auxLabelbox@=%
         \hbox to 0pt{\hss\ignorespaces\hAdjust@
          {\ignorespaces#3\unskip}\hss}%
       \vAdjust@
       \let\hAdjust@\relax\let\vAdjust@\relax
       \AugmentLabelBox@@{#2}%
       \ht\Labelbox@=0pt\dp\Labelbox@=0pt
       \let\LabelMole@\FetchLabel@%
     \fi\LabelMole@}

 \newtoks\XYSep@ 
 \def\SetXYSeparator#1{%
     \IN@0#1@@\ifIN@\XYSep@{*}%
     \else
     \XYSep@{#1}%
     \fi
     }

 \SetXYSeparator*

 \def\AugmentLabelBox@@#1{%
     \IN@0\the\XYSep@ @#1@\ifIN@
       \SPLIT@0\the\XYSep@ @#1@%
       \setbox\Labelbox@=\hbox to 0pt{%
         \unhbox\Labelbox@
         \Shifted@@{\the\Initialtoks@\Wddd@}%
         {\the\Terminaltoks@\Httt@}%
         {\box\auxLabelbox@}}%
     \else
         \ms@g{}%
         \ms@g{ !!! Bad insertion point. !!!}%
         \message{ (#1\ this point was rejected.)}%
     \fi
    }

 \def\FetchOption@#1[#2]#3\endFetchOption@{%
    \def\temp{#1}
    \ifx\temp\empty
       \Edgect@=#2\relax
       \let\EdgGl@b\relax
       \SetEdges@
       \Cleaner@#3%
    \fi}

 \def\Cleaner@#1[@]{\Labeltoks@{#1}}
     
 \def\PlaceLabels@@{\mathsurround=0pt
     \def\Cr@{\\}%
     \let\L\lAdjust@\let\R\rAdjust@
     \let\B\bAdjust@\let\E\eAdjust@\let\T\tAdjust@
     \expandafter\FetchOption@\the\Labeltoks@[@]\endFetchOption@
     \Wddd@=\Wdd@ \Edgescale@\Wddd@ 
     \Httt@=\Htt@ \Edgescale@\Httt@
     \expandafter\FetchLabel@\the\Labeltoks@\endFetchLabel@
     \box\Labelbox@
     }%

 \let \PlaceLabels@\PlaceLabels@@

 \def\AffixLabels#1{\setbox\Coordinatebox@=\hbox{#1}%
      \Wdd@=\wd\Coordinatebox@ \Htt@=\ht\Coordinatebox@
      \advance\Htt@ \dp\Coordinatebox@
      \hbox{\copy\Coordinatebox@\kern-\Wdd@ 
           \Shifted@@{0pt}{-\dp\Coordinatebox@}%
           {\PlaceLabels@\PlaceGrid@}%
           \kern\Wdd@}%
      \GridSwitch@ 
      \LabGl@b\Labeltoks@{()\\}%
      }
 
   \let\wlog\LabelFigwlog@ld   
   \catcode`\@=\LabelFigCatAt  

                                By

              Raymond S\'eroul <A18645@FRCCSC21.BITNET>
                                and 
              Laurent Siebenmann <lcs@topo.math.u-psud.fr>
    
              VERSIONS: July 1991, Oct 1991, Jan 1992, July 1992

INTRODUCTION

      This labelling package is intended for TeX users who
rely on non-TeX sources for for their graphics inserts.  It
provides means for adding TeX labels to such inserts with a
minimum of fuss. 

       For most labels, TeX users have in the past found it
reasonably convenient to rely on non-TeX sources. Typical
occasions when an inescapable need for TeX labels seemed to
arise are

 (a) when the graphics program lacks certain exotic or complex
mathematical symbols

 (b) when the very highest typographical quality is wanted for the
labels

 (c) when labels included with the graphics fail to print, 
 and you cannot figure out why (cf. boxedeps.doc).  The labels
 provided by labelfig.tex are 100

       Since this package first appeared, many users, who in the
past scarcely dreamed of using TeX labels, have come to use
nothing but.  So it is now appropriate to add

Intoxication Warning:  TeX labels may be addictive and expensive. 

     If you have a fast preview you may disagree, and even find
that this package provides an agreeable paste-up environment; see
extra applications at end.

     Note to publishers: It is possible and convenient to ultimately
export the TeX labels produced by labelfig.tex to become an integral
part of the EPS file. This is often desired by a publisher who typically
uses an "upmarket" graphics or page layout program, with which the
staff is skilled in perfecting figures.  See Appendix I for
a recipe.

     The authors are grateful to Patrick Ion of Math Reviews for
helpful comments and encouragement.

BASIC INSTRUCTIONS

    After reading in the macro file using

preview or proof your figure with a coordinate grid printed on
top, by typing the following:

    \ShowGrid  
    \AffixLabels{<the graphics insertion>}

Here <the graphics insertion> is what you would type to insert
the graphics object alone without the grid.  This must provide
for the space around it. For example <the graphics insertion>
might well be \BoxedEPSF{MyFigure scaled 700} using the
boxedeps.tex macro package (from same source); this provides a
TeX box containing the encapsulated PostScript insert specified by
the file MyFigure. \AffixLabels{...} provides the grid (supposing
\ShowGrid is present) and later, once you have specified labels
using the grid, it will "tack on" the labels.

     The grid is a sort of (usually elongated) checkerboard of
ten rows and ten columns and its (internal) partitions are by
default numbered  .1, ... ,.9  both horizontally (X-coordinate
running left to right) and vertically (Y-coordinate running bottom
to top).  Thus the points enclosed by the grid correspond to the
points of the unit square in the cartesian "X-Y" plane, the lower
left corner corresponding to the origin (0,0).  By extrapolation,
the full page corresponds to a larger rectangle in the plane.

     These coordinates serve to position labels as follows.
Before the \AffixLabels{...} command type label specifications:

  \SetLabels
   (<X-coordinate>*<Y-coordinate>) <first label> \\
   .
   .
   .
   (<X-coordinate>*<Y-coordinate>)  <last label> \\
  \endSetLabels

Each row specifies one label and is terminated by \\.  In each
row, the position indicator comes first; it is written as a
standard cartesian point except that the X- and Y- coordinates
are separated by * rather than a comma because TeX allows a
comma as decimal point. There are no dimension units to specify
as the unit is the grid itself.

     By default, this cartesian point specifies where the middle
of the baseline of the label will be located.  However if you precede
the point by \L [or \R] the left [or right] edge of the baseline will
be located there. Similarly you may also precede the point by \T, \E,
or \B to vertically align the top equator or bottom of the label box
at the specified point.  This gives nine standard positions of
the label with respect to the insertion point --- corresponding to
the eight principle points of the compas and the center

                     \L\T     \T      \R\T

                     \L\E     \E      \R\E

                     \L\B     \B      \R\B

But this neglects the default "baseline" level of TeX,
giving potentially three more positions

                     \L    <no tag>   \R

For text, the baseline level is often the preferred. Its relation to
the others is variable. It will often coincide with the bottom level,
as happens for "X".  But it is often distinct, as for "g", in which
case you have in all 12 distinct positions rather than 9.

     It is convenient to think of this specification of label
position as attaching the label by a thumb-tack to the coordinate
grid. There are up to twelve positions of the thumb-tack on the
label, while the position of the thumb-tack on the coordinate grid is
arbitrary.  Normally, one choses the position of the thumb-tack on
the label to be the one that is the closest to the item being
labeled.  There are good reasons for this "rule of thumb":

   (a)  It facilitates correct positioning at first try.

   (b)  If the scale of the figure must be altered after labels
have been affixed, the labels have a good chance of remaining well
positioned.

   (c)  The visible grid need not extend beyond the "bounding box"
for the figure, because the best preferred position is always
(at least almost) within the bounding box .

The second reason is particularly important. Indeed it often
happens that scale has to be altered after labelling begins, in
order to either provide space for the labels, or to adjust
proportions between the labels and the figure.  (The size of labels
is unaffected by scaling.)

     Here is an artificial but self-contained test which uses
TeX rules to make a graphics object.

TEST

    Do not skip this!

 \def\FrameIt#1{\hbox{\vrule$\vcenter {\hrule\kern3pt%
             \hbox {\kern3pt #1\kern3pt}%
               \kern3pt\hrule}$\relax\vrule}}

 \def\Caption#1#2{\FrameIt{%
       \vtop {\hsize=#1\relax \parindent=0pt
         \leftskip=0pt \rightskip=0pt plus15pt
         \parfillskip=0pt
         \lineskip=1pt\baselineskip=0pt
         #2}}}

 \def\FirstQuadrant{\hbox to 100pt{\vrule\vbox to 100pt{%
        \hbox to 100pt{\hfil}\vfil\hrule}\hss}}

  \SetLabels
    \R(.5*.2) $\zeta\,\cdot$\\
    (.9*-.10) $\xi$\\
    \R(-.03*.9) $\eta$\\
    \T(.5*.9) \Caption{70pt}{%
          \it The norm of
          $g(\xi+i\eta)$ is indicated on
          contours of this invisible surface.}\\
  \endSetLabels

  \AffixLabels{\FirstQuadrant}

  \end

  Note that the coordinates to use for labels are indicated on the
edges of the grid (when visible) corresponding to the conventional
x- and y- axes of the Cartesian plane. By default the grid is
1-by-1. However, by the command \Edges{100}, you can change this
to 100-by-100 and many users find this alternative most
convenient. Place the command \Edges{...} in your style file (or
header) since its effect is is global. Other possible edge values
are 10 and 1000.

  If you use the command \Edges{...} at all, do so with care.  For
if you accidentally delete an \Edges{...} command your labels will
abruptly be badly misplaced and may logically but mysteriously
generate "dimension too big" errors under TeX and "off page" errors
under your driver.  

  You can dictate the edgescale for an individual figure by giving
the scale in brackets immediately after \SetLabels.  Thus, to
import into an article using say \Edge{100} a figure labelled using
another edgescale, say the original 1-by-1 default, you can use
\SetLabels[1]...\endSetLabels.

GETTING IT DOWN PAT

     Complicated labeling deserves the same respect as
complicated mathematics.  Do not expect it to come out perfect the
first time!  What is needed in either case is a mechanism to
repeatedly typeset troublesome pieces.

     One mechanism is always available.  One does complicated
labelling in a separate "test" file involving just the figure being
labelled;  a texpert will know how to \dump TeX's current state as
a temporary format that restarts rapidly at each retry.  Usually,
one then pastes the completed labelled figure back into the main
TeX file, but, of course, one can also \input it as an auxiliary
file.

     If you do not have a TeXpert at handy, here is a first
approximation to an efficient setup. By deletions reduce a copy
of your article to just a few lines before and after the figure.
Now label the figure, and finally, copy and paste the labelled
figure to the original article. Then copy the next figure to label
into this testbed and repeat. The TeXpert can improve the  speed
at which TeX starts up, by compiling a format specifically for
your article; just one caution: best NOT include in the format
ephemeral details of setup like \Set<mydriver>ArtSpecials (from
boxedeps.tex because this reads  figure dimensions which you may
change during your work session.

     An improved mechanism to repeatedly typeset troublesome
pieces is now available on the Macintosh; it is called LinoTeX;
see the same ftp sources.  It could be set up on many types
of computer.

     Before using labelfig.tex to attach labels to a graphics
object inserted using boxedeps.tex or BoxedArt.tex, make it a
firm rule to carefully adjust the bounding box using the trimming
commands of these packages, and also at least tentatively scale
and position the object. Beware of changing the grid inadvertently
after the labels have been positioned.  For example, correcting
the bounding box of a PostScript graphics object can foul up the
labels by changing the coordinate grid to which the labels are
attached. This is particularly true for the trimming  commands of
boxedeps.tex and BoxedArt.tex. However, as noted already, change
of scale is much less disruptive, and modest adjustments should be
well tolerated.

     Sometimes the labels protrude so far from the bounding box
of a figure that the figure has to be repositioned.  Best do this
by ad hoc spacing, say using \hglue and \vglue; altering the
bounding box would create a vicious circle.

     Remember that you are responsible for preventing labels
from overlapping. You are responsible for all label typography
including size and style. A label is really just about anything
that can be put in a TeX box. Note that spaces at the beginning
and end of labels will normally be suppressed; if you really want
them you must protect them with TeX braces.

     This package temporarily sets the \mathsurround parameter
of TeX to zero  while the labels are being affixed. This is done
because nonzero \mathsurround space would influence the position
of left and right aligned labels; then, when a texpert or printer
modifies mathsurround, diagram labeling might be disastrously
altered. There is a small price to pay involving labels that are
formatted as caption boxes including mathematics: you  may want or
need to specify an explicit mathsurround space within the caption
box; it will not influence anything outside.

     Those hostile to the use of * as separator between
the X and Y coordinates of label insertion points, are free to
impose another using \SetXYSeparator{<the new separator>}.  
Americans may prefer "," to "*" since they never use a 
comma as a decimal point; on the other hand, * may be more visible.

APPENDIX (I)  MERGING labelfig.tex LABELS INTO AN EPSF GRAPHICS OBJECT.

     As promised in the introduction, here is a recipe useful for
publishers. It works at least on Macintosh and at least for vectorized
graphics and Adobe type1 fonts.  (There is surely a similar recipe for
PCs under MSWindows.)

 (a)  Use boxedeps.tex utility to integrate the figure given by the eps
file, "x.eps" say, with a visible frame around it.  See
\ShowDisplacementBoxes command in boxedeps.tex.  To get precise results
automatically it is important to use the \Trim... commands of
boxedeps.tex making the "DisplacementBox" neatly fit the figure.

 (b)  Use the TeX printer driver and LaserWriter (versions >= 8.1.1) to
export to an EPSF the DVI page containing the integrated, labelled
figure. You now have an EPS file  "xx.eps"  that contains too much, and at
the wrong scale, and at wrong position.

 (c)  Convert the EPSF to an Adode Illustrator format EPSF using
the shareware utility called epsConvert by Sam Weiss
1993-- (currently $25).

 (d)  In Illustrator (or a compatible program), group the labels and the
"DisplacementBox"; copy them to the clipboard and paste them into "x.ps".
This step requires that all the label fonts be "visible to the Macintosh.

 (e)  Translate and scale the pasted group consisting of the labels plus
the "DisplacementBox" so as to make the "DisplacementBox" the bounding
box of (labelless) figure represented by "x.eps".  At this point the
labels will be correctly placed on the figure "x.eps".

 (f)  Ungroup and delete the "DisplacementBox".  The result is the
desired single EPS file, "x+.eps" say, It contains the original figure
plus its labels.  

     Using grouping and ungrouping appropriately in "x+.eps", a
publisher's staff can very efficiently improve label positions etc.

APPENDIX II)  SOME EXOTIC APPLICATIONS

     The grid of labelfig.tex is analogous to a light-table in
classical page makeup with wax or latex glue.  In principle, you
can use it to compose any page from its indivisible parts.  This
even has some of the artisanal charm of classical paste-up
provided you have a fast screen preview to make the process
"interactive".

     In practice labelfig.tex is a tool for nonstandard jobs.
Here are a few going beyond the labelling already discussed.

(I)  GRAPHICS INTEGRATION.

     This is accomplished by treating the imported graphics
objects as labels.  The underlying graphics object is then
typically an empty  \vbox to <dimension>{\vfill} in a TeX
\midinsert...\endinsert construction.  A label line
might be of the form

   (.1*.1) \\

The exact form of the special command varies from driver to
driver.  However, in the case of encapsulated PostScript graphics
(EPSF norm), by relying on boxedeps.tex, one can have the
following standard syntax (independant of driver  (see
boxedeps.doc for details.
  
  (.1*.1) \BoxedEPSF{MyFigure scaled <scale in mils>}\\

This may be slow since it requires TeX to read the PostScript
file to read bounding box using many complex macros.  So you
may want to try

  (.1*.1) \EPSFSpecial{MyFigure}{<scale in mils>}\\

which is fast and driver independant, but it squashes the
bounding box, normally to its lower left corner.

     Similarly for graphics of the Macintosh PICT norm ---
using BoxedArt.tex (same sources) in place of boxedeps.tex.

     This approach to integration is to be recommended when
one is assembling a composite graphics object.

 (II)  COMMUTATIVE DIAGRAM ENHANCEMENT

     Commutative diagrams or arrays of mathematical objects
connected by arrows of various sorts are common in mathematics.
The mathematical objects require the use of TeX.  Recently TeX
acquired a good collection of arrows of all slopes --- that of
LamSTeX --- plus pwerful macros to build the diagrams.

     However, even the LamSTeX collection is often
inadequate; it lacks for example double shafted arrows, dotted
arrows and curved arrows. Fortunately it is possible to produce
such arrows on an individual basis using sophisticated graphics
programs such as Illustrator and AldusFreehand (both serving
the EPSF norm) or using Metafont (with its public domain norm).
Since the creation of each new arrow is a work of love, you
probably want to limit the number of arrows by using LamSTeX
for most arrows. The 40K commutative diagram module of LamSTeX
has been adapted to work with AmSTeX and a copy may be posted
with LabelFig and related files. Unfortunately no one has yet
offered a version that works with Plain TeX or LaTeX.

       Suffice it here to say that when the exotic arrow has
been somehow imported into TeX, labelfig.tex treats it as a
label that one affixes to the commutative diagram.  Two other
steps will be treated in separate notes, namely the matter of
extracting the dimension specifications for the arrow and the
construction of the arrow --- for these steps are far from
unique and often depend intimately on your computer environment. 
Notes for the Macintosh-Textures-Illustrator combination are
found in the file ExoticArrows.doc.

 (III) NESTING 

Ingenuity pays off in exploiting labelfig.tex. One can
mix graphics and typography quite freely.  labelfig.tex is good
for freeform or overlapping arrangements, while boxedeps.tex (or
BoxedArt.tex) is best for regimented non-overlapping
arrangements --- and the two can be combined.

     The default behavior of labelfig.tex is not ideal 
for nesting objects, because to prevent trouble for beginners
the register for labels is globally cleared when \AffixLabels
concludes.  But there are switches available

      \LabelsGlobal      \LabelsLocal

which change this.  To understand this, extend the above test 
by something like:

 \LabelsLocal

 \SetLabels
    (.5*.5) AAA\\
 \endSetLabels

 {
 \SetLabels
    (.5*.5) ZZZ\\
 \endSetLabels
   \AffixLabels{\FirstQuadrant}
 }

   \AffixLabels{\FirstQuadrant}

     There are however potential pitfalls.  Neither
labelfig.tex nor boxedeps.tex has been tested under extreme
conditions. Problems may occur if their procedures are
indiscriminately nested. For boxedeps.tex (not labelfig.tex)
there is a precise cause for worry, namely many of its
variables are "global", which means that TeX braces will not
provide the protection one might expect.

COMMAND SUMMARY FOR labelfig.tex

  Here [...] means optional (one or zero)
       [...]* means any number of such constructs

  \SetLabels
    [[<P>](<X><Sep><Y>) <label> \\]*
  \endSetLabels
  \ShowGrid  
  \AffixLabels{<the figure>}

   --- <P> is tack position, one of eleven or empty
              order irrelevant

                   \L\T      \T      \R\T

                   \L\E      \E      \R\E

                     \L               \R

                   \L\B      \B      \R\B

   --- (<X><Sep><Y>) insertion point;
  <Sep> is separator, = * by default;
  \SetXYSeparator{<Sep>} changes it.
   <X> and <Y> are real numbers

  --- <label> a label to attach 

  --- <the figure> the figure to label 

  \GlobalLabels (default)     
  \LocalLabels  setting for nested constructs.

 \Grids makes ALL grids appear; \HideGrid then makes just next disappear.
 \noGrids returns to default.  The commands are always global.

 \GridLineWidth{<dimension>} adjusts width of grid lines. Default is very
small, to give "hairline" effect. If your grid lines are missing try
setting \GridLineWidth{1pt}.

 \Edges#1 globally changes the edge size of all grids to the numerical 
value #1, which must be 1, 10, 100, or 1000.  The default is 1.

VERSION HISTORY.
 --- Jan 1993: \Edges#1 and [??] option after \SetLabels
 --- July 1992: \Grids, \noGrids, \HideGrid;
       Gridlines become hairlines; \GridLineWidth{<dimension>}.
 --- Oct 1991, Jan 1992: \SetXYSeparator{<Sep>},  \LabelsGlobal,
       \LabelsLocal.
 --- July 1991: first release

Address for bugs and other feedback:

        Raymond S\'eroul
        IREM and Lab. de Typographie Informatise
        Univ. Rene Descartes
        Strasbourg

    Tel 33-88-41-63-45
    Email:  A18645@FRCCSC21.BITNET

        Laurent Siebenmann
        Mathematique, Bat. 425,
        Univ de Paris-Sud,
        91405-Orsay,
        France

    Tel 33-1-6941-7949; 
    Email: lcs@topo.math.u-psud.fr

\usepackage{setspace}

\usepackage{parskip}

\usepackage{titlefoot}

\newcommand{\e}[0]{\epsilon}

\newcommand{\PP}{\ensuremath{\mathbb{P}}}

\newcommand{\N}{\ensuremath{\mathbb{N}}}
\newcommand{\R}{\ensuremath{\mathbb{R}}}

\newcommand{\Z}{\ensuremath{\mathbb{Z}}}

\newcommand{\T}{\ensuremath{\mathbb{T}}}
\newcommand{\E}[0]{\mathbb{E}}

\newtheorem{theorem}{Theorem}[section]

\newtheorem{lemma}[theorem]{Lemma}
\newtheorem{proposition}[theorem]{Proposition}
\newtheorem{corollary}[theorem]{Corollary}

\newcommand{\ignore}[1]{}
\newcommand{\yujienew}[1]{\textcolor{blue}{#1}}
\newcommand{\yujie}[1]{\textcolor{brown}{\texttt{yujie:} #1}}

\newtheorem{definition}[theorem]{Definition}

\theoremstyle{definition}

\theoremstyle{definition}

\theoremstyle{definition}

\theoremstyle{definition}

\def\lora{\longrightarrow}

\newcommand{\intint}[1]{\llbracket 1,#1 \rrbracket}

\newcommand{\mc}{\mathcal}

\newcommand{\stake}{{\rm Stake}}

\newcommand{\val}{{\rm Val}}

\newcommand{\pay}{\mathsf{Pay}}
\newcommand{\game}{{\rm Game}}
\newcommand{\tugofwar}{{\rm TugOfWar}}
\newcommand{\stateofplay}{{\rm StateOfPlay}}

\newcommand{\totvar}{{\rm TotVar}}

\newcommand{\fin}{F}

\newcommand{\din}{d_{\rm in}}
\newcommand{\dout}{d_{\rm out}}
\newcommand{\thespan}{{\rm span}}
\newcommand{\tjdg}{T_{\rm JD}}
\newcommand{\vjdg}{V_{\rm JD}}
\newcommand{\ejdg}{E_{\rm JD}}
\newcommand{\hms}{h_{\rm MS}}
\newcommand{\hps}{h_{\rm PS}}
\newcommand{\msdepth}{K}
\newcommand{\Kdepth}{K}

\newcommand{\gameneps}{\mathrm{Game}_n(\e)}
\newcommand{\meanpayoff}{M}

\newcommand{\minaconf}{P_-^{\rm conf}}
\newcommand{\minaconfn}{P_-^{\rm conf} \rfloor_n}

\newcommand{\minaconfmacro}{\minaconf}

\newcommand{\maxineconf}{P^{\rm conf}_+}
\newcommand{\maxineconfn}{P^{\rm conf}_+ \rfloor_n}

\newcommand{\maxinehighfortuneinfinity}{\mathsf{HighFortune}_+(D,\infty)}
\newcommand{\maxineconfmacro}{\maxineconf}

\newcommand{\splusinfinity}{S_+}
\newcommand{\sminusinfinity}{S_-}

\def\sta{\mathsf{s}}
\def\mov{\mathsf{m}}

\newcommand{\omegmac}{w}

\newcommand{\nwithzero}{\N}
\newcommand{\nwithoutzero}{\N_+}

\newcommand{\Pminus}{P_-}
\newcommand{\Pplus}{P_+}

\newcommand{\jointstrategy}{strategy pair } 
\newcommand{\jointstrategies}{strategy pairs } 
 
\newcommand{\jointstrategyperiod}{strategy pair. }

\newcommand{\stakemacro}{s} 

\newcommand{\boundarymac}{B} 
\newcommand{\openmac}{O}
\newcommand{\nonrootmac}{B^*}

\begin{document}


\title{Stake-governed tug-of-war \\ and the biased infinity Laplacian}

\author{Yujie Fu \qquad Alan Hammond \qquad G\'abor Pete}

\AtEndDocument{
  \bigskip
  \small
  \par

\textsc{Yujie Fu} \\
  \textsc{Department of Mathematics, U.C. Berkeley} \\
   \textsc{775 Evans Hall, Berkeley, CA, 94720-3840, U.S.A.} \\
  \textit{Email:} \texttt{yujie\_fu@berkeley.edu} \\

\vspace{-3mm}

   \textsc{Alan Hammond} \\
  \textsc{Departments of Mathematics and Statistics, U.C. Berkeley} \\
   \textsc{899 Evans Hall, Berkeley, CA, 94720-3840, U.S.A.} \\
  \textit{Email:} \texttt{alanmh@berkeley.edu} \\

\vspace{-3mm}
  
  \textsc{G\'abor Pete} \\
  \textsc{HUN-REN Alfr\'ed R\'enyi Institute of Mathematics}\\
  \textsc{Re\'altanoda u. 13-15., Budapest 1053 Hungary}\\
  \textsc{and}\\  
  \textsc{Department of Stochastics, Institute of Mathematics, Budapest University of Technology and Economics}\\
  \textsc{M\H{u}egyetem rkp. 3., Budapest 1111 Hungary}\\
  \textit{Email:} \texttt{gabor.pete@renyi.hu}\\
}

\ignore{
\author[Y. Fu]{Yujie Fu}
\address{Y. Fu\\
  Department of Mathematics and Statistics\\
 U.C. Berkeley \\
??? Evans Hall \\
  Berkeley, CA, 94720-3840 \\
 U.S.A.}
\email{???}

\author[A. Hammond]{Alan Hammond}
\address{A. Hammond\\
  Department of Mathematics and Statistics\\
 U.C. Berkeley \\
  899 Evans Hall \\
  Berkeley, CA, 94720-3840 \\
 U.S.A.}
\email{alanmh@berkeley.edu}

\author[G. Pete]{G\'abor Pete}
\address{G. Pete\\
Alfr\'ed R\'enyi Institute of Mathematics\\
Re\'altanoda u. 13-15., Budapest 1053 Hungary\\
and\\
Institute of Mathematics, Budapest University of Technology and Economics\\
M\H{u}egyetem rkp. 3., Budapest 1111 Hungary
}
\email{gabor.pete@renyi.hu}
}

\maketitle

\begin{abstract} 
 In tug-of-war, two players compete by moving a counter along edges of a graph, each winning the right to move at a given turn according to the flip of a possibly biased coin. The game ends when the counter reaches the boundary, a fixed subset of the vertices, at which point one player pays the other an amount determined by the boundary vertex. Economists and mathematicians have independently studied tug-of-war for many years, focussing respectively on resource-allocation forms of the game, in which players iteratively spend precious budgets in an effort to influence the bias of the coins that determine the turn victors; and on PDE arising in fine mesh limits of the constant-bias game in a Euclidean setting. 
 
In this article, we offer a mathematical treatment of a class of tug-of-war games with allocated budgets: each player is initially given a fixed budget which she draws on throughout the game to offer a stake at the start of each turn, and her probability of winning the turn is the ratio of her stake and the sum of the two stakes. We consider the game played on a tree, with the boundary being the set of leaves, and the payment function being the indicator of a single distinguished leaf. We find the game value and the essentially unique Nash equilibrium of a leisurely version of the game, in which the move at any given turn is cancelled with constant probability after stakes have been placed. We show that the ratio of the players' remaining budgets is maintained at its initial value~$\lambda$; game value is a biased infinity harmonic function; and the proportion of remaining budget that players stake at a given turn is given in terms of the spatial gradient and the $\lambda$-derivative of game value. We also indicate examples in which the solution takes a different form in the non-leisurely game.

\end{abstract}


\unmarkedfntext{\emph{Key words and phrases.} Allocated-budget games,  bidding games on graphs, Colonel Blotto games,  multi-turn games,  random-turn games,  resource allocation,  stake-governed games,  strategic move evaluation, tug-of-war, Tullock contests.}%


\setcounter{tocdepth}{2}
\tableofcontents

\section{Introduction}

In 1987, Harris and Vickers~\cite{HarrisVickers87} considered a model of two firms competing in research effort over time, each seeking to be the first to secure a patent.  In a model they called tug of war,
a counter is located initially at the origin. It will move adjacently on~$\Z$ at each of a sequence of turns, until it reaches either $N$ or $-N$ (where $N \in \N_+$), with the patent then awarded to one or other player.
At each turn of the race, each player is asked to nominate a non-negative effort rate. A certain random mechanism then selects the winner of the turn, with a higher effort rate for a given player improving her chances of turn victory. The turn victor moves the counter one unit in the direction that brings the patent award closer for her. Each player must pay a cost associated to the rate and duration of her effort at the turn.
Harris and Vickers' article initiated a research endeavour among economists concerning resource-allocation tug-of-war games that continues to the present day. 

In 2009, Peres, Schramm, Sheffield and Wilson~\cite{PSSW09} studied tug-of-war games in which the victor at each turn is decided according to a fair coin flip. Although this decision rule is trivial in comparison to those considered by economists, the geometric setting is much richer:~\cite{PSSW09} found a relation between tug-of-war games in domains in Euclidean space and infinity-harmonic functions, prompting a wave of attention from researchers in probability and PDE.

The economists' and mathematicians' research efforts have been vigorous,
with~\cite{HarrisVickers87}
and~\cite{PSSW09}
both garnering several hundred citations, 
but also remarkably disjoint: according to google scholar, no article had cited both papers before the arXiv release of ours; in fact, the authors of~\cite{PSSW09} even appear to have alighted on the name `tug-of-war' independently (as Scott Sheffield told one of us). It is valuable to consider how these two long strands of research may profitably be woven together,
by finding and analysing games that unify the richness of random decision rules from the economics strand with the geometric complexity in the mathematicians'. 
Indeed, although the research leading to the first version of the present article was performed in ignorance of the economics strand, it is natural to view the contribution that it makes in this vein: for here, we specify a class of resource-allocation tug-of-war games played on trees, which have been considered by Klumpp~\cite{Klumpp} on finite integer intervals, in which players at the outset receive fixed budgets, to be spent throughout the game by the respective players as they seek to gain control on the local movement of the counter and to guide it to a favourable terminal location, with never-spent funds lost.
 We solve leisurely versions of these games (in which some moves are randomly cancelled), showing how the essentially unique Nash equilibrium in each game constitutes a compromise between the demand to spend big now (in the hope of controlling the counter at its next move) and to hold funds in reserve (so that the capacity to control gameplay is not too quickly exhausted).

A companion article~\cite{H2022} presents an analysis of a resource-allocation tug-of-war game where the players fund themselves (as they do in Harris and Vickers' original formulation). 
The two articles may be viewed as providing detailed and rigorous treatments of games that interpret and develop models at the heart of the economics tug-of-war literature. Viewed in terms of mathematical game theory, our contribution is to find a setup that may be rigorously solved (namely, introducing the leisurely version and focusing on trees) and then to identify rigorously the essentially unique Nash equilibrium, tasks that we have found quite non-trivial. In addition, our Nash equilibrium enriches the structure of infinity-harmonic functions on graphs, with likely interesting PDE scaling limits. We also hope that our work will bring the PDE point of view to the attention of economists.

In the introduction, we specify our model and state our results. In Section~\ref{s.econmath}, we give a short survey of the separate developments of tug-of-war, in economics and mathematics.  
A conceptual overview that motivates the form of our main results and clarifies the main elements needed to prove them is offered in Section~\ref{s.picture}, which ends with  
 a guide to the structure of the later part of the article.

\subsection{The constant-bias tug-of-war game}
 
Though not the original form of tug-of war, the game with constant bias is the simplest to define. Reflecting the discrete context in which we will work, we start by specifying it for finite graphs.

\begin{definition}\label{d.graph}
Let a finite graph $G = (V,E)$ be given, 
with the vertex set $V$ written as the disjoint union  of two non-empty sets: the field of {\em open} play, $\openmac,$ and the set of {\em boundary} vertices,~$\boundarymac$.
The payoff function is a given map $p:\boundarymac \lora [0,\infty)$. We will refer to the triple $(V,E,p)$ as a {\em boundary-payment graph}.
\end{definition}

We specify tug-of-war with a coin of fixed bias, and its counter evolution $X: \llbracket 0, F \rrbracket \lora V$. (We denote by  $\llbracket i,j \rrbracket$ the integer interval   $\{ k \in \Z: i \leq k \leq j\}$.)
\begin{definition}\label{d.tugofwar}
The game $\tugofwar(q,v)$ has parameters $q \in [0,1]$ and $v \in V$. Set $X_0 = v$, so that the counter starts at $v$. Let $i \in \nwithoutzero$ 
be a positive integer.
At the start of the $i$\textsuperscript{th} turn, the counter is at~$X_{i-1}$. A coin is flipped that lands heads with probability~$q$. If it lands heads, Maxine wins the right to move the counter to a vertex in $V$ adjacent to its present location~$X_{i-1}$; if not, Mina wins the same right. The resulting location is $X_i$. The finish time $F$ is the minimum value of $i \in \nwithoutzero$ such that $X_i$ lies in the boundary $\boundarymac$. The game ends when turn $F$ does. A final payment of $p(X_F)$ units is made by Mina to Maxine. (Our currency is the `unit'.)
\end{definition}
The notion of game value will be recalled in Definition~\ref{d.purevalue}; in essence, it is the mean final payment when the game is correctly played. It is quite intuitive that the value $V(q,v)$ of $\tugofwar(q,v)$ exists and solves the system of equations
\begin{equation}\label{e.v}
 V(q,v)  \, = \, (1-q) \min_{\omegmac \sim v} V(q,\omegmac) \, + \,  q \max_{\omegmac \sim v} V(q,\omegmac) \, ,
\end{equation}
subject to $V(q,w) = p(w)$ for $w \in \boundarymac$. Tug-of-war was considered in the symmetric case where $q=1/2$ by~\cite{PSSW09}. The biased game, with $q \neq 1/2$, was treated by~\cite{PPS10,PeresSunic}; we will later review the fast algorithm given in~\cite{PeresSunic} for computing biased infinity harmonic functions (which solve the displayed equations).

\subsection{The stake-governed game}\label{s.sgg}
Maxine wins each turn with probability $q \in (0,1)$ by fiat in $\tugofwar(q,v)$. In many applications, decision makers face choices of how to spend a precious and limited resource in order to gain strategic advantage in an evolving random situation. {\em Allocated-budget stake-governed tug-of-war}  offers a model for such choices. Maxine and Mina will each initially be given a certain amount of money. Our convention will be that Maxine holds $\lambda$ units and Mina one unit at this time. A counter is placed on the game-board $(V,E)$ at $v \in V$, as it was in constant-bias tug-of-war. Each turn opens with a request that the two players offer stakes. These are amounts  that each player selects, drawn from her remaining reserves. The amounts are withdrawn from the respective reserves and are not returned. Each player will win the right to move the counter to an adjacent location of her choosing with a probability that equals the ratio of the amount that she has just staked and the total amount just staked by the two players.
 The game ends as it does in the constant-bias game, with the arrival of the counter in $\boundarymac$ (at $w$, say). The remaining reserves of the players are swept from them and, as in the earlier game, Mina pays Maxine $p(w)$ units. 
Thus the initial funds allocated to the players are a non-renewable resource that each player must spend during the lifetime of the game in an effort to gain control on the resting place of the counter. The funds dictate the players' capacity to control play and have no other utility, because they cannot be saved beyond the end of the game. 

How should this version of tug-of-war be played? What is the value of this game? What is its relation to constant-bias tug-of-war on a finite graph? These are the principal questions that we seek to address. Players may stake random amounts, but an inkling of a solution is offered by restricting to deterministic stakes. If Maxine always stakes the $\lambda$-multiple of the non-random stake offered by Mina, then the win probability at every turn is $\tfrac{\lambda}{1+\lambda}$. The relative budget of the players holds steady at $\lambda$; the game reduces to a copy of $\tugofwar(q,v)$ with $q = \tfrac{\lambda}{1 + \lambda}$; and the game's value is the last displayed $V(q,v)$ (which corresponds to $h(\lambda,v)$ in the notation we will adopt, as set out in Definition~\ref{d.puregamevalue}). So the answers to the questions seem simple, in the sense that the new game appears to project onto the old one. But the picture is more complicated. If Mina stakes randomly, then Maxine cannot reliably maintain the budget ratio by the strategy of stake-proportion mimicry. And, if stakes are in fact non-random  under optimal play, so that the projection to tug-of-war occurs, there is a further natural question. The players are staking a shared and non-random proportion of their present reserves at any given turn. We call this proportion the {\em stake function}. What is it? 

We will see that the suggested, simple-minded, projection to constant-bias tug-of-war is false for several simple graphs. Our principal results show nonetheless that the picture sketched above is correct for a certain class of graphs and payment functions, and we identify the stake function and explain how to play the game optimally in these cases (by characterizing the Nash equilibria). 
And in fact these results will be shown to hold not for allocated-budget stake-governed tug-of-war but for a leisurely variant thereof, in which each move in the game is cancelled by an independent event of probability $1-\e \in [0,1)$
whose occurrence is revealed to the players after stakes for the turn in question have been submitted. Our formula~(\ref{e.stake}) for the stake function is a ratio of a spatial derivative of $h(\lambda,v)$
and a multiple of a $\lambda$-derivative. This form reflects the competition between the desire to gain territorial advantage, which pushes for a high stake value via the numerator, and the need to keep money, which places downward pressure on the stake via the denominator. It is equally a competition between the short and the long term: spend now and push the counter as you wish; but pay for this tomorrow, with a depleted budget for later gameplay.
An alternative stake formula~(\ref{e.altstake}) will be proved in which the latter interpretation is manifest. 
 A battle between space and money is waged in this two-person zero-sum game, which is fought over territory and mediated by the initial provision of limited finances. 

We turn next to specifying allocated-budget stake-governed tug-of-war more carefully. This will permit us to state our main conclusions later in the introduction.

\subsection{The rules of the game}\label{s.specifying}

A parameter $\e \in (0,1]$ called the {\em move probability} is given. The starting fortune $\lambda$ is a parameter taking a value in $[0,\infty)$.
(`Fortune' is a synonym of `budget' in our description.)

Maxine and Mina prepare for the start of play. Maxine receives $\lambda$ units and Mina, one unit. A counter is placed at a given vertex $v \in \openmac$ in the field of open play.
In a sequence of turns, the counter will move along adjacent vertices in $G$ until it reaches a boundary vertex, in $\boundarymac$. If $X_i$ denotes the counter location at the end of the $i$\textsuperscript{th} turn, 
then the process $X: \llbracket 0, \fin \rrbracket \lora V$, $X_0 = v$, encodes the counter locations throughout gameplay. Here, $\fin$, the finish time, is $\min \big\{ i \in \nwithoutzero: X_i \in \boundarymac \big\}$. 

Maxine's fortune at the start of the $i$\textsuperscript{th} turn is recorded as $\lambda_{i-1}$ units, where $\lambda_{i-1} \in [0,\infty)$; thus $\lambda_0 = \lambda$. Mina's fortune at  the start of this, and every, turn equals one unit. As we describe gameplay, note that the randomness constituting any tossed coin is independent of any other randomness.

Let $i \in \nwithoutzero$. The $i$\textsuperscript{th} turn is now about to begin. The data $\stateofplay = (\lambda_{i-1},X_{i-1}) \in (0,\infty) \times V$ encodes the present state of play of the game.
 The upcoming turn takes place in four steps.

{\em First step: bidding.}
Maxine stakes a value $a_{i-1}$ that lies in $[0,\lambda_{i-1}]$. Mina stakes a value $b_{i-1} \in [0,1]$. These amounts are deducted from the fortunes of the two players, who now hold the remaining amounts $\lambda_{i-1} - a_{i-1}$ and $1-b_{i-1}$ in reserve. Maxine's relative fortune is now $\tfrac{\lambda_{i-1} - a_{i-1}}{1 - b_{i-1}}$, and we denote this by $\lambda_i$.
The value of $\stateofplay$ is updated to be $( \lambda_i , X_{i-1})$.

{\em Second step: determining if a move will take place.} The croupier tosses a coin that lands heads with probability $\e$. If it indeed lands heads, he declares that a move will take place, and the third step begins. If the coin lands tails, no move happens. In this case, we set $X_i$ equal to the present location $X_{i-1} \in \openmac$. Follow from `Before the next move' below to see what happens next. 

{\em Third step: selecting who wins the turn.} 
A coin is tossed that lands heads with probability $\tfrac{a_{i-1}}{a_{i-1} + b_{i-1}}$.
If it lands heads, Maxine is the turn victor; otherwise, Mina is. Thus a player's stake relative to the total stake  at the turn
is her probability of winning.

{\em Fourth step: the counter moves.} 
The victorious player selects an element $X_i$ in $V$ that is a neighbour of $X_{i-1}$---so that $(X_{i-1},X_i)$ is an element of the edge-set $E$ (a notion we record with the notation $X_{i-1} \sim X_i$). 

If $X_i \in \boundarymac$, the game is declared over at the end of the $i$\textsuperscript{th} turn. In this case, we set $\fin = i$. Any remaining holdings of the two players are swept from them. Mina makes a payment of $\pay$ units to Maxine. The value of $\pay$ is the evaluation $p(X_\fin)$ of the payoff function at the terminal vertex.

{\em Before the next move.} If $X_i \in \openmac$, play continues to the   $(i + 1)$\textsuperscript{st} turn. A simple coordinate change---a currency revaluation---is now made. The change does not affect later gameplay but it makes for simpler notation.  Recall that Maxine has $\lambda_{i-1} - a_{i-1}$ units in reserve, and Mina has $1 - b_{i-1}$. We revalue so that Mina has one unit. Under this accounting device, the already specified $\lambda_i = \tfrac{\lambda_{i-1} - a_{i-1}}{1 - b_{i-1}}$ is equal to Maxine's fortune
 as  the $(i + 1)$\textsuperscript{st} turn arrives. In short, $\stateofplay$ is now set equal to $(\lambda_i,X_i)$.

{\em The zero-bid and reset rules.} The attentive reader may have noticed two problems in the specification of gameplay. First, if the players stake $a_{i-1} = b_{i-1} = 0$, who wins the right to make any move that takes place at the $i$\textsuperscript{th}  turn? The rule presented in the third step above involves a badly specified $0/0$ probability. Second, if the players both go for broke, so that $a_{i-1} = \lambda_{i-1}$ and $b_{i-1} = 1$, and the game continues via $X_i \in \openmac$, how can play continue given two bankrupt players? Some arbitrary rules are needed to permit play to continue in the two scenarios: a {\em zero-bid} rule to cope with the first problem, and a {\em reset} rule to address a pair of players with zero combined fortune. In the {\em status quo} reset rule applied as the $i$\textsuperscript{th} turn begins, we set $\lambda_{i-1} = \lambda_{i-2}$ if the amounts held in reserve, $\lambda_{i-2} -a_{i-2}$ and $1 - b_{i-2}$, were both zero at the $(i-1)$\textsuperscript{st} turn. (The rule makes sense and resolves a problem only when $i \geq 2$.) In the {\em status quo} zero-bid rule, we take the win probability  $\tfrac{a_{i-1}}{a_{i-1} + b_{i-1}}$  (in the third step above) at the $i$\textsuperscript{th} turn (for $i \geq 1$) equal to $\tfrac{\lambda_{i-1}}{1 + \lambda_{i-1}}$
in the case that $a_{i-1} = b_{i-1} =0$.  Note that it may be that both rules are invoked at the $i$\textsuperscript{th} turn, if $i$ is at least two: two bankrupt players will have their fortunes restored, but they may then both stake nothing. 
 We will work with the two status quo rules. Doing so slightly simplifies some proofs, but our results may be proved under a broad range of zero-bid and reset rules. 
 We will briefly return to this matter in some closing comments, in Section~\ref{s.resetrules}.

{\em The payment rule for unfinished games.} A remaining ambiguity is that gameplay may never finish due to the counter failing to reach~$\boundarymac$. 
We stipulate that, in this event,
the value of $\pay$ equals $\max \big\{  p(v): v \in \openmac \big\}$. This maximum equals one for the class of $p:\openmac \rightarrow [0,\infty)$ we will consider. This choice is more proscriptive than the mere need to render play well-defined dictates: see  Section~\ref{s.unfinishedgames}.
 

If the move probability $\e$ equals one, then a move takes place at every turn; we may simply omit the second step to obtain a description of this form of the game, which will call the {\em regular game}
and denote by $\game(1,\lambda,v)$. The first argument is the value of $\e$; the latter two  give the initial value of $\stateofplay$. When $\e \in (0,1)$, we call the game the {\em leisurely game} (with move probability~$\e$, of course), and denote it by $\game(\e,\lambda,v)$. We may thus refer to all the games by the condition that $\e \in (0,1]$. We occasionally refer to $\game(\e)$ in verbal summary or 
when the value of $\stateofplay$ is variable.

\subsection{Principal definitions}
We specify basic notions concerning strategy, game value, biased infinity harmonic functions, Nash equilibria and graphs.

\subsubsection{Strategy}
Before any given turn (in either form of the game), Maxine and Mina glance at the balance sheet and at the board. They see Maxine's fortune and the counter location encoded in the vector  $\stateofplay \in (0,\infty) \times \openmac$ and (let us say) they know how many turns have passed. How should they play?
\begin{definition}\label{d.strategy}
Let $\mc{P}[0,\infty)$ denote the space of probability measures on $[0,\infty)$. We view this set as a metric space given by total variation distance; and as a measurable space by endowing it with the Borel $\sigma$-algebra that arises from the system of metric open sets.
 A {\em mixed} strategy for Maxine is a pair of measurable maps, the stake $\sta_+: [0,\infty) \times \openmac \times \nwithoutzero \lora \mc{P}[0,\infty)$ and the move $\mov_+:  [0,\infty) \times \openmac  \times \nwithoutzero \lora V$,
that satisfy that the supremum of the support of $\sta_+(x,v,i)$ is at most $x$, and $\mov_+(x,v,i) \sim v$, for $x \geq 0$ and $v \in \openmac$ and $i \in \nwithoutzero$.
Suppose that Maxine adheres to the strategy $(\sta_+,\mov_+)$. To play at
 the  $i$\textsuperscript{th} turn, she randomly selects her stake  of $a_{i-1}$ units by setting $a_{i-1}$ equal to an independent  
 sample of the law  $\sta_+(\lambda_{i-1},X_{i-1},i)$. She nominates the move  $\mov_+(\lambda_i,X_{i-1},i)$. 

A mixed strategy for Mina is a pair $\sta_-: [0,\infty) \times \openmac \times \nwithoutzero \lora \mc{P}[0,\infty)$ and $\mov_- : [0,\infty) \times \openmac  \times \nwithoutzero \lora V$,
with the supremum of the support of $\sta_-(x,v,i)$ being at most one, and with $\mov_-(x,v,i) \sim v$, for $x \geq 0$ and $v \in \openmac$. In following this strategy, Mina stakes $b_{i-1}$, an independent sample of the law $\sta_-(\lambda_{i-1},X_{i-1},i)$, and nominates the move $\mov_-(\lambda_i,X_{i-1},i)$, at the $i$\textsuperscript{th} turn. 

Note that, by use of mixed strategies, the players are permitted to randomize their stakes at any given turn. The definition does not permit random choices of move nomination, however. The reason for this restriction is purely technical: it does not significantly reduce the set of sensible strategies (in particular, we will still find a Nash equilibrium), but it will make it easier to talk about the essential uniqueness of the equilibrium.

If every image probability measure  in a strategy is a Dirac delta, so that the concerned player never places a randomly chosen stake, then the strategy is called {\em pure}. 
 
Let $\mc{P}_+$ and $\mc{P}_-$ denote Maxine's and Mina's spaces of pure strategies, and let $\mc{S}_+$ and $\mc{S}_-$ denote the counterpart mixed strategy spaces.
Thus, $\mc{P}_+ \subset \mc{S}_+$ and  $\mc{P}_- \subset \mc{S}_-$. 
 Each element  $(S_-,S_+) \in \mc{S}_- \times \mc{S}_+$ specifies the evolution of gameplay begun from any given state of play $(\lambda,v) \in [0,\infty) \times \openmac$, where of course the two players follow their selected strategies. Let $\e \in (0,1]$. The {\em mean payoff} $M(\e,\lambda,v,S_-,S_+)$ of the mixed strategy pair $(S_-,S_+)$ in $\game(\e,\lambda,v)$ equals the mean value of the payoff $\pay$ made at the end of the game dictated by~$(S_-,S_+)$.  Often we will abbreviate $M(S_-,S_+) = M(\e,\lambda,v,S_-,S_+)$. 
\end{definition}

The basic perspective of Definition~\ref{d.strategy} is that players' decisions may be made in light only of the present state of play of the game and the turn index. Note how the stake component at the $i$\textsuperscript{th} turn depends on the input $(\lambda_{i-1},X_{i-1},i)$, while move nomination depends on $(\lambda_i,X_{i-1},i)$. The rules of the game stipulate that stakes are made at the start of the turn, with $\stateofplay$ equal to $(\lambda_{i-1},X_{i-1})$. The players hand over their stakes and the turn victor is decided, with $\stateofplay$ updated to $(\lambda_i,X_{i-1})$, before the move victor makes her move; thus she does so in view of the relative fortune $\lambda_i$ at this time.

\subsubsection{Game value and Nash equilibria}

\begin{definition}\label{d.purevalue}
Consider $\game(\e)$ for $\e \in (0,1]$. Note that the inequality
 \begin{equation*}
 \sup_{S_+ \in \mathcal{S}_+}  \inf_{S_- \in \mathcal{S}_-} M(\e,\lambda,v,S_-,S_+ ) \geq 
  \inf_{S_- \in \mathcal{S}_-}  \sup_{S_+ \in \mathcal{S}_+} M(\e,\lambda,v,S_-,S_+)
 \end{equation*}
always holds. If there is equality here, then the single quantity is called the {\em value} $\val(\e,\lambda,v) \in [0,\infty)$ of the game begun with $\stateofplay = (\lambda,v)$. If the inequality is strict, the value does not exist.
%
\end{definition}

\begin{definition}\label{d.nash}
A pair $(S_-,S_+) \in \mc{S}_- \times \mc{S}_+$ is a Nash equilibrium if, for all  choices of the initial data $(\lambda,v) \in [0,\infty) \times \openmac$,
the bounds 
$$
M(\e,\lambda,v,S'_-,S_+) \geq M(\e,S_-,S_+) \, \, \, \textrm{and} \, \, \,   M(\e,,\lambda,v,S_-,S'_+) \leq M(\e,\lambda,v,S_-,S_+)
$$ 
hold for every
$S'_- \in \mc{S}_-$ and $S'_+ \in \mc{S}_+$. A Nash equilibrium is called {\em pure} if it lies in $\mc{P}_- \times \mc{P_+}$ and {\em mixed} if it does not.
\end{definition}

It is easy to see that, if a Nash equilibrium $(S_-,S_+)$ exists, then the mean payoff $M(\e,\lambda,v,S_-,S_+)$ is independent of the choice of this equilibrium, and is given by the value~of~the~game.

\subsubsection{Biased infinity harmonic functions}

\begin{definition}\label{d.puregamevalue}
For $\lambda \in [0,\infty)$, let $h(\lambda,\cdot): V \lora [0,\infty)$ denote the $\lambda$-biased\footnote{Peres--\v{S}uni\'c would call this function $\lambda^{-1}$-biased. Our parametrization follows an albeit-perhaps-unfair tendency to adopt the point of view of the maximizing player, who associates a higher probability of winning to a higher value of the bias.} infinity harmonic function on $V$ with boundary data $p: \boundarymac \lora [0,\infty)$.
This function is the solution to the system of equations
\begin{equation}\label{e.h}
 h(\lambda,v) =
\begin{cases}
\, \, p(v) \, \, \textrm{for} \, \, v \in \boundarymac  \, , \\
\, \,  \tfrac{\lambda}{\lambda +1} \max_{u \sim v} h(\lambda,u) +  \tfrac{1}{\lambda +1} \min_{u \sim v} h(\lambda,u)  \, \, \textrm{for} \, \, v \in \openmac \,   ,
\end{cases}
\end{equation}
subject to $h(v) = p(v)$ for $v \in \boundarymac$. 
Recall that, by $u \sim v$, we mean that $u \in V$ with $(u,v) \in E$.
\end{definition}

The solution exists and is unique: see~\cite[Section~1.5]{PeresSunic}. We will recall an efficient algorithm of Peres and \v{S}uni\'c \cite{PeresSunic} to compute this value in Section~\ref{s.tools}. The structure of this solution is particularly simple for root-reward trees, which we next define. As we will see, the key feature of these trees is that the solution can be iteratively computed using a decomposition of the tree into paths in a way that does not depend on the value of $\lambda$.

%

\subsubsection{Root-reward trees}\label{s.rrt}
The focus will lie on a class of boundary-payment graphs. These are trees with a unique leaf supporting the payment function.

\begin{definition}\label{d.uniquepaymenttree}
Recall that a boundary-payment graph comprises a triple $(V,E,p)$ where the finite graph $(V,E)$ has vertex set~$V$ written as a disjoint union $\openmac \cup \boundarymac$ and where the payment function~$p$ maps $\boundarymac$ to $[0,\infty)$. Such a triple is called a {\em root-reward tree} when the following conditions are met:
\begin{itemize}
\item the graph $(V,E)$ is a rooted tree whose root~$r$ is a leaf in the tree;
\item the boundary set $\boundarymac$ is equal to the set of leaves of the tree; and
\item the payment function equals ${\bf 1}_r$.
\end{itemize}
The non-root boundary $B^*$ equals $B \setminus \{ r \}$.

In this way, a terminal payment, of one unit, 
is made if the game ends by the counter reaching $r$, rather than~$B^*$. The root~$r$ may also be 
called the {\em reward} vertex, if a name that reflects Maxine's wishes is preferred.
 
\end{definition}

\subsection{Main results}\label{s.mainresults}
We will study the leisurely game on root-reward trees. Here is our first principal conclusion.

\begin{theorem}\label{t.leisurely}
 Let $T = (V,E,{\bf 1}_r)$ be a root-reward tree. There exists $\e_0 \in (0,1)$ such that, when $\e \in (0,\e_0)$ and $\lambda >0$, the following hold for each $v \in \openmac$. 
 \begin{enumerate}
\item 
The value of $\game(\e,\lambda,v)$ exists and equals $h(\lambda,v)$.
\item
The game $\game(\e,\lambda,v)$ has a Nash equilibrium.
Any Nash equilibrium is a pair of pure strategies 
that lead the players  to  stake non-random amounts that maintain the fortune ratio: that is, $\lambda_i = \lambda$ for $i \in \intint{\fin -1}$ in the gameplay that arises 
at the Nash equilibrium.
\end{enumerate}
The value of $\e_0$ may be taken to be $\vert V \vert^{-2 \vert \boundarymac \vert -1}/4$.
\end{theorem}


Theorem~\ref{t.leisurely}(2) indicates that, under jointly optimal play, a shared and non-random proportion of each player's present reserves will be staked. The question is obvious: what is this proportion? 
Our second principal result, Theorem~\ref{t.nashform}, offers an answer. To state it, we first need a technical result.  

\begin{proposition}\label{p.differentiable}
On root-reward trees, for $v \in \openmac$, the function $\lambda \mapsto h(\lambda,v)$ is differentiable on $\lambda \in (0,\infty)$.
\end{proposition}

\begin{definition}\label{d.stake}
Set $\Delta: [0,\infty) \times \openmac  \lora [0,1]$ according to
$$
\Delta(\lambda,v) \, =  \, \max_{\omegmac \sim v} h(\lambda,\omegmac) \, - \, \min_{\omegmac \sim v} h(\lambda,\omegmac) \, .
$$
Then, Proposition~\ref{p.differentiable} permits us to define the {\em stake} function ${\rm Stake}: (0,1] \times  [0,\infty) \times \openmac  \lora [0,\infty)$ as
\begin{equation}\label{e.stake}
{\rm Stake}\big(\e,\lambda,v \big) = \frac{\e \, \Delta(\lambda,v)}{(\lambda + 1)^2\tfrac{\partial}{\partial \lambda} h(\lambda,v)} \, .
\end{equation}
\end{definition}

 

\begin{definition}\label{d.inout}
A directed edge $(u,v)$ in~$T$ is called {\em in} if $d(r,v) = d(r,u) - 1$ and {\em out} if $d(r,v) = d(r,u) + 1$. A path in $T$ is called out if each of its edges is out. 
For $v \in V$, let the in-distance $\din(v)$ equal $d(r,v)$; and the out-distance $\dout(v)$, the minimum length of an out path from $v$ to  the non-root boundary $\nonrootmac$.

The parent $v_+$ of $v \in V \setminus \{ r \}$ is the unique neighbour~$u$ of $v$ such that $\din(u) = \din(v) -1$. For $v \in V \setminus B^*$, the set
 $$
 \mc{V}_-(v) \, = \, \Big\{ \, u \in V: u \sim v \, , \,   u \neq v_+ \, , \,   \dout (u) =  \dout(v) -1 \,  \Big\} \, .
 $$
 distinguishes some of the children of~$v$.
\end{definition}

\begin{definition}\label{d.conforming}
Consider a \jointstrategy $(S_-,S_+) \in \mc{S}_- \times \mc{S}_+$ in $\game(\e,\lambda,v)$ and the gameplay~$X:\llbracket 0, F \rrbracket \lora V$ governed by this \jointstrategy in this game.
Let $i \in \nwithoutzero$. 
\begin{enumerate}
\item
Maxine {\em conforms} at the $i$\textsuperscript{th} turn (under this gameplay) if at this turn  she stakes the quantity~$\lambda_{i-1} \stake \big( \e,\lambda_{i-1},X_{i-1} \big)$ 
and nominates the move $\big( X_{i-1} \big)_+$. The strategy $S_+$ is said to be conforming against $S_-$ if Maxine almost surely conforms at every turn.
\item  Mina conforms  at the $i$\textsuperscript{th} turn (under the same gameplay) if she stakes   $\stake \big( \e,\lambda_{i-1},X_{i-1} \big)$ 
and nominates a move to an element of $\mc{V}_- \big( X_{i-1} \big)$ at this turn. The strategy $S_-$ is said to be conforming against $S_+$ if Mina almost surely conforms at every turn.
\end{enumerate}  
The \jointstrategy $(S_-,S_+)$ is said to be conforming if each component is conforming against the other.
\end{definition}
Here then is the promised result that elucidates the form of Nash equilibria from Theorem~\ref{t.leisurely}(2).
\begin{theorem}\label{t.nashform}
Let $(V,E,{\bf 1}_r)$ be a root-reward tree.
With $\e_0 \in (0,1)$ specified in Theorem~\ref{t.leisurely}, the condition $\e \in (0,\e_0)$ implies that each Nash equilibrium in $\game(\e,\lambda,v)$ is
a conforming \jointstrategyperiod
\end{theorem}
In the stake formula~(\ref{e.stake}) a battle between space and money is evident in the relative sizes of numerator and denominator. This is also a battle between the short and long terms, as we will see shortly, in Proposition~\ref{p.totvar}. A little notation is needed, so that the stake formula denominator may be identified as the mean time remaining in a suitably calibrated clock.
When the counter lies at a vertex $v \in \openmac$ for which $\big\vert \mc{V}_-(v) \big\vert \geq 2$, Mina may choose among the members of $\mc{V}_-(v)$ when she nominates her move with the counter at~$v$ while retaining a conforming strategy. We now specify a collection of gameplay processes that arise from different strategies for Mina compatible with such choices. 
Proposition~\ref{p.totvar} shows that in fact it makes no difference, at least for the purpose of measuring the mean calibrated time remaining in the game,  which of these choices Mina makes; this result thus provides an alternative formula for the stake function~$\stake(\e,\lambda,v)$.
\begin{definition}\label{d.thetafirst}
Let $(V,E,{\bf 1}_r)$ be a root-reward tree. Let $\Theta$ denote the set of mappings $\theta: (0,\infty) \times \openmac \times \nwithzero  \lora V$ such that $\theta(\lambda,v,i) \in \mc{V}_-(v)$ for each $v \in \openmac$. (We write $\nwithzero = \nwithoutzero \cup \{ 0\}$.)
\end{definition}

\begin{definition}\label{d.theta}
 For $\e \in (0,1]$, $\theta \in \Theta$ and $\lambda \in (0,\infty)$, 
let $X_\theta(\e,\lambda,\cdot): \llbracket 0,F_\theta \rrbracket \lora V$ denote {\em $(1-\e)$-lazy $\lambda$-biased} walk on $V$ (with index $\theta$).
 This is the Markov process such that  $X_\theta(\e,\lambda,0) =v \in \openmac$  and, whenever $X_\theta(\e,\lambda,k) \in \openmac$ for $k \in \nwithzero$,
 $$
 X_\theta(\e,\lambda,k+1) =
\begin{cases}
\,  X_\theta(\e,\lambda,k)  \, \, \textrm{with probability $1-\e$}  \, , \\
\,  \big( X_\theta(\e,\lambda,k)\big)_+ \, \, \textrm{with probability $\e \tfrac{\lambda}{1 + \lambda}$}  \, , \\
\, \theta\big(\lambda,X_\theta(\e,\lambda,k),k\big) \, \, \textrm{with probability $\e \tfrac{1}{1 + \lambda}$}  \,   .
\end{cases}
$$
The process is stopped on arrival at $\boundarymac$, so that $F_\theta = \min \big\{ j \geq 0: X_\theta(\e,\lambda,j) \in \boundarymac \big\}$.

The first two arguments will usually be $\e$ and $\lambda$, and we will often use the shorthand $X_\theta(\cdot) = X_\theta(\e,\lambda,\cdot)$.

The next result expresses the conclusion of Theorem~\ref{t.nashform} in terms of counter evolution.
\end{definition}
\begin{corollary}\label{c.gameplay}
 Under the hypotheses of Theorem~\ref{t.nashform},  the set of gameplay processes governed by Nash equilibria in $\game(\e,\lambda,v)$ is given by $\big\{ X_\theta: \theta \in \Theta \big\}$. 
\end{corollary}
Set $\totvar(\e,\lambda,\theta,v) = \sum_{i=0}^{F_\theta -1} \Delta \big(\lambda,X_\theta(i)\big)$ where $X_\theta(0) = v$.

\begin{proposition}\label{p.totvar}
 Let $T = (V,E,{\bf 1}_r)$ be a root-reward tree.
For $\lambda \in [0,\infty)$,  the value of the mean $\E \, \totvar(\e,\lambda,\theta,v)$ is independent of $\theta \in \Theta$, and, in an abuse of notation, we denote it by  $\E \, \totvar(\e,\lambda,v)$.
Then $\E \, \totvar(\e,\lambda,v) = \e^{-1} (\lambda + 1)^2\tfrac{\partial}{\partial \lambda} h(\lambda,v)$ for $\e \in (0,1]$. As such, we obtain an alternative formula for the stake function~(\ref{e.stake}): for $\e \in (0,1]$,
 \begin{equation}\label{e.altstake}
{\rm Stake}\big(\e,\lambda,v \big) \, = \, \frac{\Delta(\lambda,v)}{\E \, \totvar(\e,\lambda,v)} 
\, = \,  \frac{\e \, \Delta(\lambda,v)}{\E \, \totvar(1,\lambda,v)} \, .
\end{equation}
\end{proposition}

A further main conclusion is an explicit formula for $\stake(1,\lambda,v)$, $(\lambda,v) \in (0,\infty) \times \openmac$, on a root-reward tree~$(V,E,{\bf 1}_r)$. The formula is expressed in terms of the Peres-\v{S}uni\'c 
 decomposition of the tree, which is presented in Section~\ref{s.tools}, and it appears as Theorem~\ref{t.stakeformula} at the end of that section.  
We end the introduction by stating two consequences of Theorem~\ref{t.stakeformula} that do not require new notation. The first is an explicit formula for the special case where $V$ is an integer interval. The second states instances of the general formula where simple but interesting forms appear: 
these are the cases when $\lambda$ is low; when it is one; and when it is  high.

Let $n \in \nwithoutzero$. The {\em line graph} $L_n$ has vertex set $\llbracket 0, n \rrbracket$; its edges connect consecutive vertices. 
By choosing its root to be $n$, and taking $p = {\bf 1}_n$, we obtain a root-reward graph~$(L_n,\sim,{\bf 1}_n)$ (where the edge-set is abusively denoted by~$\sim$). 

\begin{proposition}\label{p.stake}
For $n \in \nwithoutzero$, consider the triple $(L_n,\sim,{\bf 1}_n)$.  
Let $i \in \intint{n-1}$.
\begin{enumerate}
\item For $\lambda \in (0,\infty) \setminus \{ 1 \}$, 
$$
 \stake(\e,\lambda,i) \, = \, \e \cdot \bigg( \tfrac{\lambda +1}{\lambda - 1} \cdot i \cdot \Big( 1 - \tfrac{n(\lambda^i - 1)}{i(\lambda^n - 1)} \Big) \bigg)^{-1} \, .
$$
\item  For $\lambda \in (0,\infty) \setminus \{ 1 \}$, , 
 $\stake(\e,\lambda,i) = 
 \stake(\e,\lambda^{-1},n-i)$.
\item We further have that  $\stake(\e,1,i) =  \e \cdot \tfrac{1}{i(n-i)}$.
\item 
The preceding results are equally valid when the payment function ${\bf 1}_n$ is replaced by any choice of $p:\{0,n\} \lora [0,\infty)$ for which $p(0) < p(n)$. 
\end{enumerate}
\end{proposition}



\begin{proposition}\label{p.threelambda}
Let $(V,E,{\bf 1}_r)$ be a root-reward tree, and  let $v \in V$. Writing $\stake(\lambda,v)$ in place of $\stake(1,\lambda,v)$, we have that
$$
 \lim_{\lambda \searrow 0} \, \stake (\lambda,v) = \din(v)^{-1}  
\, \, \, , \, \, \, \stake(1,v) = \dout(v)^{-1} \din(v)^{-1} 
$$
and
$$
 \lim_{\lambda \nearrow \infty} \, \lambda^{\dout(v) \, - \, \dout^{\rm min}(v)} \cdot \stake (\lambda,v) \, = \,  \Big( \vert J \vert \cdot \dout^{\rm min}(v) \Big)^{-1} \, , 
 $$
where $\dout^{\rm min}(v)$ equals the minimum of $\dout(w)$ as $w$ ranges over vertices on the geodesic path $[r,v]$, and~$J$ is the set of vertices~$w$ on $[r,v]$ such that $\dout(w) = \dout^{\rm min}(v)$. 
\end{proposition}
When $v \in \intint{n-1}$ lies on the line graph $L_n$, $\lim_{\lambda \nearrow \infty} \stake(\lambda,v)$ equals $\dout(v)^{-1}$, and there is a certain symmetry between the players in the preceding result. But in general $\stake(\lambda,v)$ may tend rapidly to zero in the limit of high $\lambda$, even as this function remains bounded away from zero as $\lambda$ tends to zero.

\subsubsection{Acknowledgments}
The authors thank Judit Z\'ador for help in preparing Figures~\ref{f.lthree},~\ref{f.tgraph} and~\ref{f.poisson}, and Anoushka Chitnis for helpful comments. 
They further thank a referee for useful input.
The first author is supported by  NSF grant DMS-2153359; the second,
by
NSF grants DMS-1855550, DMS-2153359 and DMS-2450726 and a 2022 Simons Fellowship in Mathematics; and the third by the ERC Consolidator Grant 772466-NOISE and Synergy Grant 810115-DYNASNET, and the Hungarian National Research, Development and Innovation Office grant K143468.

\section{Survey of the two strands of tug-of-war research}\label{s.econmath}

Harris and Vickers~\cite{HarrisVickers87} 
initiated a wave of attention to resource-allocation tug-of-war among economists;
about twenty years later, Peres, Schramm, Sheffield and Wilson~\cite{PSSW07,PSSW09} did something similar for constant-bias tug-of-war among mathematicians. 
We now present a brief overview of the two strands of research. 
The reader may consult the introduction of~\cite{H2022} for a slightly more expansive account.

\subsection{Tug of war in economics}\label{s.econ}

 In 1980, Lee and Wilde~\cite{LeeWilde} considered a one-step model for how research effort influences outcome. Mina and Maxine (as we call them)
 choose respective effort rates $x,y \in (0,\infty)$. 
 At an exponentially distributed random time~$M_+$ of mean $x^{-1}$, Maxine makes a discovery.
 Mina does so at an independent exponential time $M_-$ of mean~$y^{-1}$. The player who discovers first, at time $M = \min \{ M_-,M_+\}$, receives a reward. During their research efforts, Mina and Maxine incur costs at respective rates $c_-(x)$ and $c_+(y)$, so that  their costs are $M c_-(x)$ and $M c_+(y)$. The resulting mean costs with joint effort rates $(x,y)$ are $\tfrac{c_-(x)}{x+y}$ and  $\tfrac{c_+(y)}{x+y}$ because the mean of $M$ is $(x+y)^{-1}$. 
 The game is player funded, in the sense that, while a player operates under no constraint on budget, the cost she incurs is deducted from her reward. 
 
In introducing tug-of-war, Harris and Vickers~\cite{HarrisVickers87} employed  
 Lee and Wilde's model to decide each given turn of the game. Players specify effort rates at which turn, and the Lee-Wilde rule determines the outcome according to given
 cost functions $c_-$ and $c_+$ that are independent of the turn index. With Mina playing left and Maxine right according to who wins at each turn, the counter evolves on the finite integer-interval~$\llbracket 0,N \rrbracket$, with $N \in \N_+$, until Mina or Maxine wins 
 when the counter reaches $0$ or $N$. Again the game is player funded, with the sum of a player's running costs deducted from her terminal receipt. A Nash equilibrium is called symmetric if Mina's effort rate at $i$ is equal to Maxine's at $N-i$. And the game is called symmetric if $c_- = c_+$ and the players' terminal receipts are equal after interchange of endpoints. 
Harris and Vickers proved the existence of a symmetric equilibrium for a fairly broad range of cost functions, and computed this equilibrium in some special cases for the cost. Under the computed equilibrium, players invest greatly when the counter is in the middle of the gameboard, with rapid but asymmetric decays in effort as it moves away. (This conclusion is verified and developed in the classification of Nash equilibria in~\cite{H2022} for a player-funded tug-of-war game; the very staccato structure of an intense battle amid near peace stands starkly in contrast to the allocated-budget games that we analyse here, where budgets are spent down progressively over the lifetime of the game, and stakes tend to be higher near the boundary.)

Several authors~\cite{KonradKovenock05,AgastyaMcAfee,Konrad2012} have studied player-funded tug-of-war on a finite interval with the all-pay auction rule, under which the player who offers the greater stake at a given turn is the turn victor. In~\cite{Hafner17,Hafner16} team contest versions of tug-of-war are analysed: these are player-funded games in which each player is in fact a team of individuals, each of whom has only a limited liability in regard to the overall running costs of the team.  

There is a second geometric context often found in economists' treatment: the first-to-$N$,  or majoritarian objective game, which ends with victory for the first player who wins a given number of the stage contests. Player-funded versions of this game have been studied since Harris and Vickers~\cite{HarrisVickers87}, with sudden bursts of early expenditure offered as an explanation~\cite{KlumppPolborn}  of the `New Hampshire' effect in U.S. presidential primary contests: here we see roughly the same `battlefield' feature as we discussed a moment ago; the qualitative difference between player-funded and allocated-budget games is far more fundamental than between finite interval tug-of-war and the first-to-$N$ version. All-pay auction rules and intermediate-time payoffs~\cite{KonradKovenock09}  or discounting~\cite{Gelder} have been studied for player-funded games, as has a team contest version~\cite{FuLuPan}. 
Ewerhart and Teichgr\"aber~\cite{ET19} consider player-funded tug-of-war framework more generally than integer interval or majoritarian objective contests. Invoking assumptions that include a certain exchangeability in the order of moves  for gameplay during open play, they show that gameplay occurs in two phases: a first, in which no state may be revisited; and then a copy of tug-of-war.

The game we study here has allocated budgets, and such games have a long history. 
 In 1921, Emile Borel introduced games~\cite{Borel1921} that he illustrated with an example in which players $A$ and $B$ select an ordered triple of non-negative reals that sum to one. 
 The player who nominates at least two values that exceed the opponent's counterparts is the winner. 
In general, an officer and his opponent each distribute a given resource across a certain number of battlegrounds, with a given battle won by the player who devotes more resource there, and the 
war won according to the battle count. The reader is encouraged to think up sensible strategies for the case of three battlegrounds, and they may notice that it is far from straightforward to come up with Nash equilibria, let alone to find all of them. The first equilibria were found in \cite{BorelVille}, with a complete characterization (but still not a full list) obtained only in \cite{Roberson} for the original game with continuously divisible resources, and in \cite{Hart} in a discretized version. The officer has come to be called Colonel Blotto, and the  name `Blotto' is attached to games where players have fixed resources, 
which may be committed simultaneously as in Borel's game, or sequentially, as in the allocated-budget tug-of-war games which are our object of study. 

Two articles~\cite{KlumpKonradSolomon,Klumpp} treat allocated-budget tug-of-war and thus have some close points of contact with the present article. To describe these works, some notation regarding single-stage contests is useful. In a Tullock contest~\cite{Tullock}, player $A$ stakes $x \in [0,\infty)$ and player $B$, $y \in [0,\infty)$, and $A$ wins the contest with probability $\tfrac{x^\gamma}{x^\gamma + y^\gamma}$, where the Tullock exponent is $\gamma \in (0,\infty)$. The contest we consider, with $\gamma =1$, is called a lottery; 
while fair-coin tug-of-war or all-pay auction versions arise in the limits of low or high $\gamma$.
 In Section~\ref{s.sgg}, we briefly indicated 
  a premise for how to play the allocated-budget tug-of-war games that we study in this article, which are Blotto games with the lottery contest function $\tfrac{x}{x+y}$: 
  as we will see in developing the premise in Section~\ref{s.picture}, it proposes that each player
  stakes the same proportion, namely~(\ref{e.stake}) with $\e = 1$, of her remaining fortune. But as we will also later discuss, a graph as simple as $\llbracket 0,3 \rrbracket$ is a counterexample to this premise. It is by making the game leisurely that our study overcomes such counterexamples. A natural alternative would be to seek a class of graphs for which the premise is correct. The majoritarian objective games---tug-of-war played on rectangle graphs with north or east pointing edges---are a good starting point. As we will see in Section~\ref{s.asimplecase}, in a set of $n$ sequential contests where one player needs to win every round,
 the optimal stake proportion is $1/n$ and  validates~(\ref{e.stake}).  Klumpp, Konrad and Solomon~\cite{KlumpKonradSolomon} generalise this example to show that the game on any directed rectangle graph, and for a broad range of turn contest functions, has an equilibrium with an even split of resources across turns. 
 In~\cite{Klumpp}, Klumpp examines  allocated-budget tug-of-war on a finite integer interval. In note~$1$ on page $27$, he notes an example, similarly as we will in Subsection~\ref{s.ellthree}, in which a formula equivalent to~(\ref{e.stake}) with $\e = 1$
 is invalidated when the weaker player is one step from winning. While we turn to a leisurely version of the game to regularise this and a much larger class of examples, Klumpp lowers the Tullock exponent from $\gamma = 1$: from a numerical investigation, he indicates that the concerned formula will offer an equilibrium precisely when $\gamma$ is at most one-half.

\subsection{Mathematical contributions}

\subsubsection{Tug-of-war and the infinity Laplacian}
The article~\cite{PSSW09} played fair-coin tug-of-war in a Euclidean setting. On winning a turn, Mina or Maxine moves the counter---now a point in a domain in~$\R^d$---to a point of her choosing that is at most $\e$ from its present location, until one moves the counter into the domain's boundary, when Mina pays Maxine the evaluation of a real-valued function defined on that boundary. Peres, Schramm, Sheffield and Wilson proved that the value of this game in the low~$\e$ limit equals the infinity harmonic extension to the domain of the boundary data.  Namely, the limiting game value~$h$ solves the equation~$\Delta_\infty h = 0$ on the domain, where the operator $\Delta_\infty$ is the infinity Laplacian. 
This is a degenerate second-order differential operator; formally at least, a solution~$h$ of the stated equation has vanishing second derivative in the direction of its gradient. 
The infinity Laplacian is a subtle and beautiful object, whose study was given impetus by this game theoretic perspective. 
There is a notion of viscosity solution~\cite{CIL92} for $\Delta_\infty h = 0$ that is a form of weak solution governed by the maximum principle and that is characterized as the absolutely Lipschitz minimizing function that interpolates given boundary data~\cite{Jensen93}.  Solutions enjoy uniqueness~\cite{Jensen93} but limited regularity~\cite{Savin05,EvansSavin}. The game theoretic point-of-view that~\cite{PSSW09} offers led to simpler and intuitive proofs of uniqueness~\cite{ArmstrongSmart2012}; this point-of-view may be expressed by a dynamic programming principle~\cite{MPR2012}.
When the counter suffers a random order~$\e$ movement after each turn in $\e$ tug-of-war,  the $p$-Laplacian~\cite{Lindqvist}, 
which for $p \in (2,\infty)$ interpolates the usual Laplacian operator and the infinity version, 
describes game value in the low $\e$ limit~\cite{PeresSheffield}: see~\cite{Lewicka}
for a survey centred around this perspective.  
The abundant connections between tug-of-war and PDE are reviewed in the book~\cite{BlancRossi}. 

Likely, the simplest game-theoretic role of discrete infinity harmonic functions---solutions of~(\ref{e.h}) with~$\lambda$ equal to one---is  as game value for tug-of-war on graphs. 
This was not however the first connection to be found.
Richman games were defined in~\cite{LLPU} and further analysed in~\cite{LLPSU}. A first-price Richman game is a version of allocated-budget tug-of-war in which an auction occurs at each turn (so that the player who stakes more wins), and with the turn victor's payment being made to his opponent, who makes no payment at the turn.  The threshold ratio is a value~$t$ such that, if Maxine's budget constitutes a proportion of at least $t$ of the combined budget of the two players at the game's outset, then Maxine has a strategy that will permit her to win the game almost surely, while Mina has such a strategy if this proportion is less than~$t$. In~\cite{LLPU}, it is shown that the threshold ratio exists as a function of initial counter position and is equal to the discrete infinity harmonic extension of $f:B \lora \R$. Several variants---with all-pay rules for the auction, or an infinite duration for the game, or a poormen variant, where the turn victor's payment goes to the bank, not to the opponent---have been analysed~\cite{AJZ,AHC,AHI,AIT20} by theoretical computer scientists, who also discuss approximate algorithms and computation complexity: see~\cite{AH22} for a survey of these directions.

\subsubsection{Random-turn selection games} 
Hex is an alternating move two-player game that was invented by Piet Hein in 1942 and rediscovered by John Nash in 1948. 
A finite domain in the hexagonal lattice is delimited by four boundary segments, consecutively labelled red--blue--red--blue. Red and blue players alternately place like-coloured hexagons on as-yet-unplayed faces in the domain, in an effort to secure a crossing in the given colour between the opposing boundary segments of that colour. It is a classical fact that, in a suitably symmetric domain, the player who moves first has a winning strategy, but this strategy is unknown except on smaller gameboards; indeed, the game is played on boards of given height and width such as eleven  in international competitions. In 2007,  Peres, Schramm, Sheffield and Wilson~\cite{PSSW07} introduced a variant of Hex, in which the right to move at any given turn is awarded to the red or blue player according to the flip of a fair coin. A simple and striking analysis reveals the existence and explicit form of the optimal strategies for 
this {\em random turn} variant of Hex. Namely, the victorious player at a given turn 
should choose the hexagon that is most likely to be pivotal for forging the desired path when the unplayed gameboard is completed by an independent critical percolation (where each unplayed face is coloured red or blue according to a fair coin flip). Thus, in jointly optimal play, the players always choose the same hexagon, coloured according to the winner of the coin flip, and hence the gameplay can be considered as a carefully optimized random order in which an unknown fair random colouring of the board is revealed. This form of the explicit strategy holds true for any random turn {\em selection game} (i.e., where the winner is given by a Boolean function of the two-colourings of a base set), and therefore, jointly optimal play is an {\em adaptive algorithm} (also called a {\em randomized decision tree} \cite{OSSS}) to determine the output of a Boolean function on independent fair random input. These random-turn-game algorithms often appear to have interesting properties \cite{PSSW07}, such as low revealment, which is important, for instance, in proving sharp thresholds and noise sensitivity for Boolean functions (see \cite{OSSS,DCRT,SchSt,GarbanSteif}). However, these algorithmic connections are far from being well-understood in general. 

\subsubsection{Broader connections and conclusion}
To conclude this brief tug-of-war survey, we first mention two broad 
antecedents to stake-governed games. Isaacs' differential games~\cite{Isaacs65} include models where a pursuer seeks to capture an evader in Euclidean space, each instantaneously selecting control variables such as velocity, just as Mina and Maxine make ongoing choices of stakes. Shapley's stochastic games~\cite{Shapley1953} involve  transition probabilities governed step-by-step by decisions of two players, and the stake game framework stands under their umbrella in a broad sense.\footnote{Both differential and stochastic games have been the subject of more recent attention. Several chapters of the handbook~\cite{HandbookDynamicGameTheory} concern dynamic game theory, and~\cite{StochasticGamesRelatedTopics} collects articles that treat stochastic games.} 
Stake games concern how  geometry influences player decisions and position strength. This theme has been explored for iterative network-bargaining problems in~\cite{KleinbergTardos,Azar,Kanoria}.

These are fairly loose connections, however. Regarding tug-of-war itself, we have reviewed the two separate research strands. 
 Mathematicians have focused on constant-bias models,
exploring a much broader geometric setting than the  finite integer intervals and (as they may be called) the directed rectangle graphs seen in much of the economics literature, 
 and the implications that games in this setting have in probability and PDE. Tug-of-war originated in economics, from the outset with bidding rules far beyond the trivial constant-bias, and has received attention from economists for almost four decades, with the theme of resource allocation at the heart of  investigation throughout. 
This article offers a detailed mathematical treatment of a class of resource-based tug-of-war games with allocated budgets; its companion~\cite{H2022}
does so for player-funded games. 
We hope that  tug-of-war in its original resource-allocation guise, with its beautiful mathematical structure
and its modelling of economic behaviour, will attract the broader attention of mathematicians.

\section{The big picture: 
the global saddle hope and the Poisson game}\label{s.picture}

Here we offer some heuristics that motivate our main conclusions and their underlying hypotheses. Arguments sometimes lack rigour, and the rest of the paper does not rest logically on this section.

Our main assertions, Theorems~\ref{t.leisurely} and~\ref{t.nashform}, may be said to solve allocated-budget stake-governed tug-of-war for the applicable graphs and parameters, since these results determine game value and all practically important properties of Nash equilibria.
The results are proved only for leisurely games on root-reward trees, however. In order to offer a simple means of seeing how the   $\lambda$-biased infinity harmonic function may arise as game value in stake-governed tug-of-war\footnote{Here and later, we omit the term `allocated-budget' from the name `stake-governed tug-of-war'.} in larger generality, we state (and shortly prove) the next result.  

For $\lambda \in (0,\infty)$, the $\lambda$-idle zone of a boundary-payment graph $(V,E,p)$ is the set of $v \in \openmac$ such that $h(\lambda,w) = h(\lambda,v)$ for every neighbour $w \in V$ of $v$. 
We will see in Lemma~\ref{l.a}(3) that the $\lambda$-idle zone of any root-reward tree is empty for $\lambda \in (0,\infty)$. 

\begin{proposition}\label{p.purevalue}
Let $(V,E,p)$ be a boundary-payment graph, and
let $\e \in (0,1]$. For $(\lambda,v) \in [0,\infty) \times \openmac$, 
suppose that the $\lambda$-idle zone of $(V,E,p)$ is empty and that
there exists a pure Nash equilibrium in $\game(\e,\lambda,v)$. Then $\val(\e,\lambda,v)$ exists and equals $h(\lambda,v)$.  
\end{proposition}

This is in essence our only result proved in a generality beyond root-reward trees. (The assumption of idle zone emptiness is a minor convenience that will facilitate the proof.) 
The proposition is however conditional and is not really one of our principal conclusions.  We prove it next, in Section~\ref{s.strategymimicry}, by a simple argument of strategy mimicry. We present the result and its derivation because they offer a hypothesis for how to play stake-governed tug-of-war; later in Section~\ref{s.picture}, we will interrogate the likely validity  of the hypothesis, thus motivating the formulation of Theorems~\ref{t.leisurely} and~\ref{t.nashform}. The hypothesis in question, which may be contemplated for any boundary-payment graph, is that 
the Nash equilibrium existence assumption of Proposition~\ref{p.purevalue} indeed holds, and that
 the joint stake
 at the hypothesised pure Nash equilibrium in   $\game(\e,\lambda,v)$ takes the form
   $(\lambda,1) \cdot S$. This hypothesis is natural, given the proposition, because a constant fortune ratio of $\lambda$ is compatible with the $\lambda$-biased form~$h(\lambda,v)$ for game value; nevertheless, the proposition does not claim this. Proposition~\ref{p.purevalue} offers no clue either as to the form of the non-random stake function~$S$, while Theorem~\ref{t.nashform} and Proposition~\ref{p.totvar}, in the special cases where they apply, did identify this function as $\stake(\e,\lambda,v)$, as given in~(\ref{e.stake}) and~(\ref{e.altstake}). 
    In Section~\ref{s.stakeperturb}, we present a heuristic argument that purports to identify $S$ for the regular game as being equal to (\ref{e.altstake}) with~$\e=1$. We will then draw attention to some difficulties that arise in trying to make this argument rigorous; and, in Section~\ref{s.localglobal}, we will challenge the purported conclusion for the regular game by analysing aspects of this game on three simple graphs. From these examples, we will learn that it is unrealistic to seek to prove Theorem~\ref{t.nashform} and Proposition~\ref{p.totvar} for the regular game without significant new restrictions. A particular problem that we will find is that one or the other player may gain by deviating from the proposed equilibrium by {\em going for broke}, staking everything in an effort to win the next turn and end the game at that time. The leisurely game undermines the efficacy of go-for-broke staking. In Section~\ref{s.poisson}, we analyse a Poisson game, which may be viewed as a game $\game(0^+)$ of infinite leisure. In the Poisson game, time is continuous, and moves occur when Poisson clocks ring. We do not study this game rigorously, because significant questions arise in formulating and analysing it. Rather, we study it formally, where it acts as an $\e = 0^+$ idealization of the leisurely game. Indeed, we will see how formal calculations involving second derivatives of value for a constrained version of the Poisson game indicate that Theorem~\ref{t.nashform} and the stake formula~(\ref{e.stake}) may be expected to hold for this game. By the end of Section~\ref{s.picture}, we hope to have indicated clearly our motivations for studying the leisurely game and
the rough form of the principal challenges that lie ahead as we prepare to prove
 Theorem~\ref{t.nashform}. The section thus concludes with an overview of  the structure of the remainder of the paper.

\subsection{Strategy mimicry: proving Proposition~\ref{p.purevalue}}\label{s.strategymimicry}

\begin{lemma}\label{l.assumenashexists}
Consider $\game(\e,\lambda,v)$ for $\e \in (0,1]$, $\lambda \in (0,\infty)$ and $v \in \openmac$.
Suppose that the game has a pure Nash equilibrium. Suppose further that, for all $P_- \in \mc{P}_-$, there exists $P_+ \in \mc{P}_+$ such that $M(P_-,P_+) \geq h(\lambda,v)$; and that, for 
 all $P_+ \in \mc{P}_+$, there exists $P_- \in \mc{P}_-$ such that $M(P_-,P_+) \leq h(\lambda,v)$. Then $\val(\e,\lambda,v)$ exists and equals $h(\lambda,v)$.
\end{lemma}
{\bf Proof}. Denote the hypothesised pure Nash equilibrium by $(P_-^0,P_+^0) \in \mc{P}_- \times \mc{P}_+$. Write $m = M(P_-^0,P_+^0)$. It follows directly from Definitions~\ref{d.purevalue} and~\ref{d.nash}  that $m$ equals  $\val(\e,\lambda,v)$. 

By assumption, there exists $P_+ \in \mc{P}_+$ such that $M(P_-^0,P_+) \geq h(\lambda,v)$. But $m \geq M(P_-^0,P_+)$ since $(P_-^0,P_+^0)$ is a Nash equilibrium. Hence $m \geq h(\lambda,v)$. Similarly, there exists $P_- \in \mc{P}_-$ such that $M(P_-,P_+^0) \leq h(\lambda,v)$. Since $m \leq M(P_-,P_+^0)$, we find that $m \leq h(\lambda,v)$. Hence $m = h(\lambda,v)$, and the lemma is proved. \qed


{\bf Proof of Proposition~\ref{p.purevalue}.}
In view of Lemma~\ref{l.assumenashexists}, two assertions will suffice. First, for 
any  pure strategy $P_-$ for Mina in $\game(\e)$, there exists a pure strategy $P_+$ for Maxine in this game such that $\meanpayoff \big(P_-,P_+ \big) \geq h(\lambda,v)$. Second,  for
any given pure strategy $P_+$ for Maxine in $\game(\e)$, there exists a pure strategy $P_-$ for Mina in this game such that $\meanpayoff\big(P_-,P_+\big) \leq h(\lambda,v)$. 

The two assertions have very similar proofs and we establish only the second of them. 
Suppose then given a pure strategy $P_+$ for Maxine in $\game(\e)$. Let $\theta: \openmac \lora V$
be such that $\theta(v)$ is an $h(\lambda,\cdot)$-minimizing neighbour of $v$ for each $v \in \openmac$.
Let $P_-$ denote the strategy for Mina in which she stakes the $\lambda^{-1}$-multiple of the stake offered by Maxine at any given turn and proposes the move $\theta(v)$ when $\stateofplay = (\lambda,v)$. If $X$ denotes the resulting counter evolution, then $h(\lambda,X)$ is a supermartingale, because 
$$
 \E \, \big[ h(X_{i+1}) \big\vert X_i \big] \, \leq \, (1-\e) h(\lambda,X_i) + \e \tfrac{1}{1 + \lambda} \min_{u \sim X_i} h(\lambda,u) +  \e \tfrac{\lambda}{1 + \lambda} \max_{u \sim X_i} h(\lambda,u) \, = \, h(\lambda,X_i) \, .
$$
Since $h$ takes values in $[0,1]$, the supermartingale converges to a limiting value. Provided that the game ends in finite time almost surely, this limiting value is almost surely equal to $\pay$.
To verify that the game ends almost surely, note that,
when Mina wins the right to move, the value of $h(\lambda,\cdot)$ evaluated at the counter location strictly decreases, because the $\lambda$-idle zone of $(V,E,p)$ is empty.
 Let $d$ be the maximum length of a path in $\openmac$. We see then that, if Mina wins the right to move on $d$ consecutive occasions at which a move takes place, then the game will not endure beyond the last of these moves, because vertices cannot be revisited. Since the probability of Mina enjoying such a string of successes from a given move is positive (it is $(\lambda + 1)^{-d}$), we confirm that the game will indeed end almost surely. Thus the  payment that Mina makes is indeed almost surely $\lim_n h(\lambda,X_n)$, whose mean is by Fatou's lemma at most $\liminf \E \, h(\lambda,X_n)$, which is at most $\E \, h(\lambda,X_0) = h(\lambda,v)$ because $h(\lambda,X_n)$ is a  
  supermartingale. 
This confirms the second of the two assertions 
whose proof we promised in deriving Proposition~\ref{p.purevalue}. \qed 

\ignore{The next result is in essence a consequence of a simplification of the preceding derivation. We omit the proof but discuss differences between the derivations. 
 Mixed strategies in the stake game do not necessarily maintain  the fortune ratio, so Proposition~\ref{p.purevalue} requires the hypothesis that a pure Nash equilibrium exists so that its proof may merely address pure strategies in view of Lemma~\ref{l.assumenashexists}. For constant-bias tug-of-war, this difficulty with mixed strategies is absent, and indeed the classical fact~\cite{VonNeumann} that any two-person zero-sum game with finite strategy spaces has value may play the role of Lemma~\ref{l.assumenashexists}; thus,  the value of classical biased tug-of-war  (determined over mixed strategies) is unconditionally identified in the next result.
 \begin{corollary}\label{c.value}
Let $(V,E,p)$ be a boundary-payment graph, and let $\lambda \in (0,\infty)$ and $v \in \openmac$. Suppose that the $\lambda$-idle zone of $(V,E,p)$ is empty. 
The value of the game $\tugofwar\big(\tfrac{\lambda}{1 + \lambda},v\big)$ exists and equals~$h(\lambda,v)$.
\end{corollary}
The idle-zone emptiness assumption merely permits the adaptation of the proof of Proposition~\ref{p.purevalue}: Corollary~\ref{c.value} may be  readily obtained without this assumption from the study of biased tug-of-war in~\cite{PPS10} with suitable adaptations to handle the context of graphs.}

\subsection{The stake formula argued via perturbation}\label{s.stakeperturb}

Consider the regular game, with $\e = 1$. When the idle zone is empty, Proposition~\ref{p.purevalue} implies that the {\em assumption} that
 \begin{itemize}
\item  [A1]:
 there exists a pure Nash equilibrium in $\game(1,\lambda,v)$ for all $(\lambda,v) \in (0,\infty) \times V_O$
 \end{itemize}
 leads to the conclusion that
 \begin{itemize}
 \item [A2]:
 the game value $\val(1,\lambda,v)$ exists and equals $h(\lambda,v)$.
 \end{itemize}
 The proposition also strongly suggests, though it does not trivially imply, that
 \begin{itemize}
 \item [A3]:
 at any given turn, the players stake a non-random and shared proportion of their remaining fortunes and offer move nominations that maximize or minimize~$h(\lambda,\cdot)$ among neighbours of the counter location. 
 \end{itemize}
 We will not attempt to pin down whether [A3] follows from the basic premise [A1]. Our focus is on constructing and examining a picture of how the regular game would be played if we admit [A1,2,3], which we will collectively denote by [A], as well as certain further expressions of regularity or good behaviour which we will mention as the need for them arises.

 

Admit then [A], and denote the shared stake proportion by $S(\lambda,v)$ for $(\lambda,v) \in (0,\infty) \times \openmac$. 
The relative fortune of the players will remain constant, at its initial value $\lambda$, and the counter will move according to  $\lambda$-biased infinity harmonic walk. 

We will heuristically identify the form of the stake $S(\lambda,\cdot): \openmac \lora [0,1]$ as a function of the vertices in open play, finding  $S(\lambda,v)$ to equal the right-hand side of the alternative stake formula~(\ref{e.altstake}), (with $\e=1$, since we consider the regular game). This right-hand side is a ratio. In the numerator is $\Delta(\lambda,v)$, namely the difference in mean payment according to whether Maxine wins or loses the first turn. In the denominator is $\E \, \totvar(1,\lambda,v)$, which is the mean value of the sum of such differences as the counter evolves during gameplay. The derivation will develop a theme mooted in Section~\ref{s.mainresults}: the ratio may be interpreted as a short-term gain divided by a long-term cost.


Set $S=S(\lambda,v)$. Maxine and Mina will play according to the Nash equilibrium in~[A1]. 
 With $\stateofplay= (\lambda,v)$, Mina will stake $S$ and Maxine, $\lambda S$. They will nominate moves to minimize or maximize the value of $h(\lambda,\cdot)$ among neighbours of $v$. To find the form of $S$, we consider the possibility that, for the first turn only, Maxine perturbs her $S$-value by a small positive amount $\eta$, so that she instead stakes $\lambda(S+\eta)$. After the first turn, she adheres to the Nash equilibrium strategy.  How will she be affected by this change? The net effect on mean payment will be a difference $G-L$ of two positive terms, a gain~$G$ and a loss $L$. The gain $G$ is short term: Maxine may benefit by winning the first turn due to increased expenditure. The loss $L$ is the price that she pays due to depleted resources as the second term begins.

We want to find the forms of $G$ and $L$, of course. To compute $G$, note first that the increase in Maxine's first turn win probability due to her alteration is
$$
 \frac{\lambda(S+\eta)}{\lambda(S + \eta) + S} - 
 \frac{\lambda}{\lambda + 1}  \,\,  \, \, \textrm{which equals} \, \, \, \, \eta S^{-1} \tfrac{\lambda}{(1+\lambda)^2} \, ,
$$
where we will treat $\eta$ as an infinitesimal, so that $\eta^2 = 0$. If Maxine converts a loss into a win in this way, then her expected mean payment increases by 
$$
 \max_{u \sim v} h(\lambda,u) \, - \, \min_{u \sim v} h(\lambda,u)  
$$
which equals $\Delta(\lambda,v)$. The gain term $G$ is the product of the two preceding displays: 
\begin{equation}\label{e.gain.heur}
G =  \eta S^{-1} \tfrac{\lambda}{(1+\lambda)^2} \Delta(\lambda,v) \, .
\end{equation}

What of the loss term~$L$?  Maxine enters the second turn poorer than she would have been. Indeed, if 
$\lambda_{\rm alt}$ denotes the relative fortune at the second turn (and later) when Maxine alters her stake at the first turn, we have  $\lambda_{\rm alt} = \lambda - \lambda (1-S)^{-1} \eta$.


The counter evolves from the second turn as a $\lambda_{\rm alt}$-biased walk in view of~[A3]. The difference in the probability that Mina wins one of the later turns in the altered gameplay (when $\eta > 0$) but not in the original one (when $\eta = 0$)
equals  the positive quantity
\begin{equation}\label{e.altlambda}
\lambda/(1+\lambda) - \lambda_{\rm alt}/(1+\lambda_{\rm alt}) \, = \,  (1-S)^{-1} \lambda(1+\lambda)^{-2} \eta \, ,
\end{equation}
where a brief calculation shows the displayed equality. What price will Maxine pay for any lost opportunities to move at these later turns? An argument in the style of the derivation of Russo's formula for percolation will tell us the answer. Since $\eta$ is infinitesimal, we may neglect the possibility that Maxine loses out twice.
If she loses out when the counter is at $w \in \openmac$, the mean payment at the end of the game will drop by $\Delta(\lambda_{\rm alt},w)$. This equals  $\Delta(\lambda,w)$ for our purpose because such a loss is incurred only with probability~(\ref{e.altlambda}), and $\eta^2 = 0$, provided at least that $\Delta(\cdot,w)$ is differentiable at~$\lambda$.
By averaging over the counter trajectory $X: \llbracket 0 , F \rrbracket \lora V$, we see that the loss in mean payment after the first turn equals
\begin{equation}\label{e.loss.heur}
L \, = \, \eta \cdot (1-S)^{-1} \lambda(1+\lambda)^{-2}   \cdot \E \sum_{i=2}^F \Delta(\lambda,X_{i-1}) \, ,
\end{equation}
where note that the turn index~$i$ runs from two; also recall that $F$ is the time at which $X$ reaches~$\boundarymac$. 
(In fact, neglecting the possibility that Maxine loses out twice requires some assumption on the tail decay of $F$, such as that $\E F^2 < \infty$.)

In the pure Nash equilibrium, Maxine will indeed stake $\lambda S$ at the first turn, which means that $G$ equals $L$: if $G > L$, Maxine could choose a small $\eta > 0$; if $G < L$, a small $\eta < 0$; and in either case, she would gain. Thus,~(\ref{e.gain.heur}) equals~(\ref{e.loss.heur}), so that 
$$ 
  S^{-1} \Delta(\lambda,v)
 =  (1-S)^{-1}  \E \, \sum_{i=2}^F \Delta(\lambda,X_{i-1}) \, .
 $$
 Rearranging, we find that
 \begin{equation}\label{e.s}
 S \, =  \, \frac{\Delta(\lambda,v)}{\E \, \sum_{i=1}^F \Delta(\lambda,X_{i-1})} \, ,
 \end{equation}
 where the summation (or turn) index $i$ now begins at one because $X_0 = v$.
Since this denominator equals $\E \, \totvar(1,\lambda,v)$ as it is specified before Proposition~\ref{p.totvar}, we have thus completed a heuristic argument for the formula for the stake function $S = S(\lambda,v)$ given in~(\ref{e.altstake}) for $\e=1$. 
 

 We will note two criticisms of this heuristic or of the conclusions that we may be tempted to draw from it. The next definition is useful for expressing the first criticism and we will also use it later.
 \begin{definition}\label{d.saddle}
 Let $I$ and $J$ be two intervals in $\R$ and suppose given $F: I \times J \lora \R$. A point $(x,y) \in I \times J$ is a global saddle point in $I \times J$ if $F(x^*,y) \geq F(x,y) \geq F(x,y^*)$ for all $(x^*,y^*) \in I \times J$.
 A point is a local saddle if there exist some intervals $I^* \subseteq I$ and $J^* \subseteq J$ whose interiors respectively contain $x$ and $y$, and the previous bounds hold for $(x^*,y^*) \in I^* \times J^*$. The word `minimax' is a synonym of `saddle' for these definitions.
 \end{definition}

 \subsubsection{First criticism: the predicted saddle point may merely be local}\label{s.criticismone}
 
 The above argument is perturbative, and it claims only that $(\lambda S, S)$ is a local saddle point of the map that associates to each $(a,b) \in [0,\lambda] \times [0,1]$
 the value  of the game resulting from joint stakes of $(a,b)$ at the first turn. (We will introduce notation for this value very soon.) 
 The prospect that the saddle point is global needs further examination. Now, at present, we are operating under [A1], and the posited  pure equilibrium  guarantees the existence of a global saddle point of the one-step stake game at any given turn, 
 except in situations of go-for-broke type where stakes lie at the boundary.  
This structure does not allay the criticism because the global saddle may not be unique; 
and because the stake-pair induced by an equilibrium
may fail to coincide with  
the local saddle identified by the perturbative argument: notably, when this equilibrium prompts stakes that lie at the boundary,  the vanishing-derivative conditions that determine the latter saddle need not be satisfied. 




 \subsubsection{Second criticism: the saddle point formula may be badly defined}\label{s.criticismtwo}

 The posited pure Nash equilibrium may not be unique, and
one or other player may have a choice of move nomination at some $\stateofplay = (\lambda,v)$. One choice may typically lead $\game(\e)$ to end more quickly than another. This means that the denominator in~(\ref{e.altstake}) is not canonically determined by the data~$(\lambda,v)$, because different equilibrium gameplay processes~$X$ that specify $\totvar(\e,\lambda,v)$ may lead to different values of $\E  \, \totvar(\e,\lambda,v)$.
 When we spoke of averaging over the trajectory $X: \llbracket 0 , F \rrbracket \lora V$ to derive~(\ref{e.loss.heur}), we implicitly used a uniqueness assumption about this trajectory: that the mean value $\E  \, \totvar(\e,\lambda,v)$ is well defined as a function of the state $(\lambda,v)$, and unaffected by equilibrium selection. Even the identification of a local saddle point may fail if this assumption is false.

 
 \subsection{The global saddle hope consistency check}\label{s.localglobal}
 
 
Suppose that not only are the assumptions [A] from Section~\ref{s.stakeperturb} in force, but also that the identified local saddle point---namely, $(\lambda S,S)$ with $S$ given by~(\ref{e.s})---is global (allaying the first criticism) and well defined (likewise the second). When this happens, we  say that the {\em global saddle hope} is realized. We will interrogate this hope by examining the regular game on three simple graphs---and will see that the hope is often disappointed.

 Let $(V,E,p)$ be a boundary-payment graph and let  $(\lambda,v) \in (0,\infty) \times \openmac$.
 For $a\in [0,\lambda)$ and $b \in [0,1)$, consider the constrained version $\game(1,\lambda,v,a,b)$ of $\game(1,\lambda,v)$ in which Maxine must stake $a$, and Mina $b$, at the first turn. 
 (The appearance of five parameters in this order will be characteristic when first-turn-constrained games are considered.)
 Denote the value of the constrained game by $\val(1,\lambda,v,a,b)$. For the three graphs we will investigate, we will carry out what we may call a {\em GloshConch}: a {\em Gl}obal {\em s}addle {\em h}ope {\em Con}sistency {\em ch}eck.
To perform the GloshConch for a given graph, we suppose [A]: regular game value is the biased harmonic function, for every $(\lambda,v)$.
We will then be able to compute  constrained-game values $\val(1,\lambda,v,a,b)$
because, after the first turn, these games reduce to copies of the unconstrained game  $\game(1)$. We will be in a position to determine whether  the map $[0,\lambda) \times [0,1) \lora [0,\infty):(a,b) \mapsto \val(1,\lambda,v,a,b)$ has a saddle point at $(\lambda S,S)$, where $S = \stake(1,\lambda,v)$ 
 is specified in~(\ref{e.s}).
The global saddle hope promises that, whatever the value of $(\lambda,v) \in (0,\infty) \times \openmac$,
this point is a saddle, and indeed that the saddle is global. If this promise is delivered for all such parameters, we say that the GloshConch is passed; otherwise, it fails.  Analysing a graph and finding that the GloshConch fails may entail identifying interesting structure for the maps of constrained value. We cannot however make detailed inferences  from the GloshConch's failure: merely that [A] does not hold for all $(\lambda,v)$. Further work would be needed to elucidate the structure of Nash equilibria and find the value of games for which the consistency check fails.

 
 
 \subsubsection{Normal form contour plots}
 A simple graphical representation will allow us to perform the GloshConch in three simple graphs. 
 For a given boundary-payment graph, $\lambda \in (0,\infty)$ and $v \in \openmac$, the function  
 $$
 [0,\lambda]\times [0,1] \lora [0,1]: (a,b) \mapsto \val(1,\lambda,v,a,b)
 $$ 
 may be depicted as a contour plot in the $[a,b]$-rectangle. Mark points~$(a_0,b_0)$ in  $[0,\lambda]\times [0,1]$ as red if they are maximizers of $[0,\lambda] \lora (0,\infty): a \mapsto \val(1,\lambda,v,a,b_0)$. Mark points~$(a_0,b_0)$ in  $[0,\lambda]\times [0,1]$ as blue if they are minimizers of $[0,1] \lora (0,\infty): b \mapsto \val(1,\lambda,v,a_0,b)$. We refer to the resulting sketches as normal form contour plots: see Figures~\ref{f.lthree} and~\ref{f.tgraph} for several examples that we will shortly discuss. The resulting red and blue curves (which, as the sketches show, sometimes have discontinuities) are Maxine and Mina's respective best responses to first-turn stakes offered by Mina and Maxine at given vertical and horizontal coordinates. A global saddle point for the plotted function occurs when the red and blue curves meet, provided that the function is continuous at the intersection point. (We will see shortly that the curves may however meet at a point where this continuity is lacking.)

 \subsubsection{Example I: The line graph $(L_2,\sim,{\bf 1}_2)$}

The first of the three graphs has only one vertex in open play. From~$1$, the game will end in one turn. If Maxine wins, the counter moves to~$2$, and $\pay=1$. If she loses, the counter moves to~$0$ with no payment made. Note that $\val(1,\lambda,1,a,b)$ equals~$\tfrac{a}{a+b}$, with a global minimax achieved at~$(\lambda,1)$. This corresponds to the trivial conclusion that there is a unique Nash equilibrium where both players go for broke: what else could they do in a game with one turn?
 
 \subsubsection{Example II: The line graph $(L_3,\sim,{\bf 1}_3)$}\label{s.ellthree}
 
\begin{figure}[htbp]
\includegraphics[width=1\textwidth]{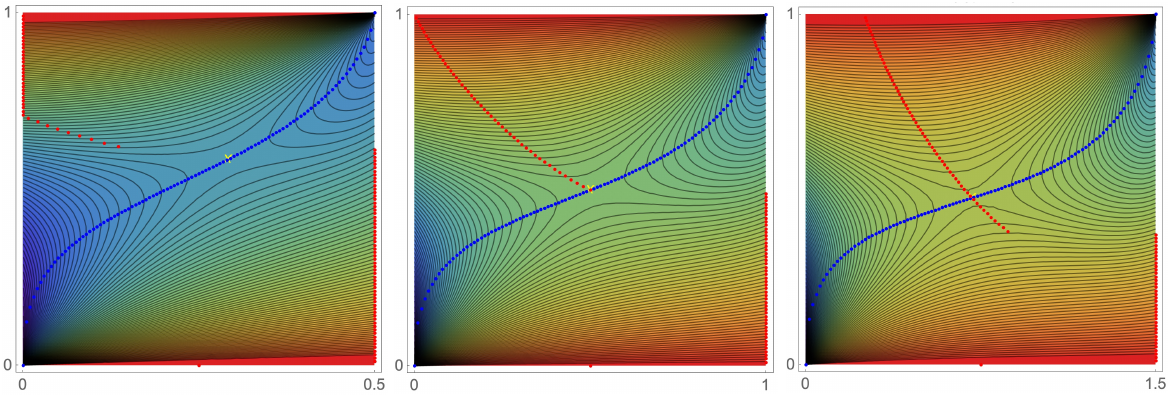}
\caption{Normal form contour plots  of $[0,\lambda] \times [0,1] \lora [0,1]: (a,b) \mapsto \val(1,\lambda,v,a,b)$ 
for vertex $v = 2$ in the line graph $(L_3,\sim,{\bf 1}_3)$, drawn under the assumptions [A] of the GloshConch. The $\lambda$-values are $1/2$, $1$ and $3$ from left to right.
}\label{f.lthree}
\end{figure}

 For this graph, the vertices $1$ and $2$ are in open play. Mina must play left and Maxine right, so move nomination is trivial. But stake decisions are less evident than they were in  the preceding example. The $\lambda$-biased infinity harmonic values are  
 $$
 h(\lambda,1) = \lambda^2 (\lambda^2 + \lambda + 1)^{-1} \, \, \, \, \textrm{and}   \, \, \, \, h(\lambda,2) = \lambda (\lambda +1) (\lambda^2 + \lambda + 1)^{-1} \, .
 $$
 For $a\in [0,\lambda)$ and $b \in [0,1)$,  $\val(1,\lambda,1,a,b) = \tfrac{a}{a+b} h (\lambda_1,2)$ and  $\val(1,\lambda,2,a,b) = \tfrac{a}{a+b} + \tfrac{b}{a+b}  h(\lambda_1,1)$, where $\lambda_1 = \tfrac{\lambda - a}{1 -b}$ is the relative fortune after the first turn in the constrained game. The predicted Nash equilibrium point is $(\lambda S,S)$ with $S$ given by~(\ref{e.altstake}) and simple computations for biased random walk on a line segment, as in Proposition~\ref{p.stake}, equal to $(1+\lambda + \lambda^2)(2\lambda +1)^{-1}(\lambda + 1)^{-1}$. (It is understood throughout this discussion of regular-game examples that $\e=1$ in~(\ref{e.altstake}.) See Figure~\ref{f.lthree} for three normal form contour plots for the constrained game at vertex~$2$ on this graph. While the yellow cross $(\lambda S, S)$ marks a local saddle point for each $\lambda \in (0,\infty)$, it is only when the condition $\lambda \geq 1$ (such as in the middle and right plots) is met that this saddle point is global. Indeed, when $\lambda < 1$ (as in the left plot) and Mina stakes $S$ with $\stateofplay = (\lambda,2)$, Maxine prefers $\lambda$,  the go-for-broke choice, to $\lambda S$ as her response. In an example as simple as $L_3$, the GloshConch fails---the global saddle hope is disappointed---for many values of $\lambda$. The first criticism, offered in Subsection~\ref{s.criticismone}, is valid.

 \subsubsection{Example III: The $T$ graph}\label{s.tgraph}

In the third example, move nomination as well as stake decisions are non-trivial for at least one player. Take copies of the line graphs $L_2$ and $L_3$ and identify vertex~$1$ in $L_2$ with vertex~$3$ in~$L_3$. The result is the $T$ graph. It has two vertices in open play, which we call the north vertex $N$ and the south vertex $S$, with $N$ being the just identified vertex. We label the three leaves in the $T$ graph $0$, $1$ and $2$, where $0$ is adjacent to $S$ and $1$ and $2$ are adjacent to $N$. We will consider the payment function given by the identity map on the leaves of $T$ with the just indicated notational convention: we name the leaves by the value of $\pay$ that they offer. 

Consider $\game(1,\lambda,N)$ on this boundary-payment graph. Mina has a choice between $1$ and $S$ when she nominates a move. This is a choice between a short and a long game, since the former nomination will end the game should she win the resulting move and the latter will keep the counter in open play, at $S$. 

To compute the form of $h(\lambda,N)$ and $h(\lambda,S)$ subject to the boundary condition that $h(\lambda,i) = i$ for $i \in \llbracket 0,2\rrbracket$, note that the $T$ graph may either be viewed as a copy $[1,N,2]$ of $L_2$ to which the path $[0,S,N]$ has been adjoined, or as a copy $[0,S,N,2]$ of $L_3$ to which $[1,N]$ has been adjoined. 
With the former view, we may specify $h_1(\lambda,N)$ as the value $\tfrac{2\lambda + 1}{\lambda + 1}$ of $\lambda$-biased infinity harmonic value at\footnote{The bold font disambiguates this vertex reference from that reserved for $1$ and~$2$.}~${\bm 1}$ on $L_2$ with boundary values $1$ and $2$, and then take  $h_1(\lambda,S) = \tfrac{\lambda}{\lambda + 1} h_1(\lambda,N)$. With the latter view, we set $h_2(\lambda,N) =  2\lambda (\lambda +1) (\lambda^2 + \lambda + 1)^{-1}$
and $h_2(\lambda,S) = 2\lambda^2 (\lambda^2 + \lambda + 1)^{-1}$ to be the $\lambda$-biased values at ${\bm 2}$ and ${\bm 1}$ on $L_3$ with boundary data $0$ and $2$.
On the $T$ graph, 
for $v \in \{N,S\}$,
$$
h(\lambda,v) =
\begin{cases}
\,  h_1(\lambda,v)  \, \, \, \textrm{when $h_2(\lambda,S) \geq 1$}  \, , \\
\,  h_2(\lambda,v) \, \, \, \textrm{when $h_2(\lambda,S) <1$}  \,   .
\end{cases}
$$
There is thus a critical value $\lambda_c$, which equals the golden ratio 
$\tfrac{\sqrt{5}+1}{2}$,
such that $h_1$-values are used when $\lambda \geq \lambda_c$, and $h_2$-values are used in the opposing case.

 For $a\in [0,\lambda)$ and $b \in [0,1)$, we will compute  $\val(1,\lambda,N,a,b)$ by use of~[A2].
  To do so, we need to understand whether Mina will choose to nominate $1$ or $S$ at the first turn of $\game(1,\lambda,N)$.
 Writing $\lambda_1 = \tfrac{\lambda - a}{1-b}$, she will nominate $1$ if $\lambda_1 > \lambda_c$, because $h(\lambda_1,1) = 1 \leq h(\lambda_1,S)$, and she will nominate $S$ if $\lambda_1 < \lambda_c$, because then  $h(\lambda_1,1) = 1 > h(\lambda_1,S)$.  (If $\lambda_1 = \lambda_c$, she may nominate either of the two moves. Note also that she is able to make the indicated choices, because the updated relative fortune $\lambda_1$ has been encoded in $\stateofplay$ by the time she makes her move, as we emphasised after Definition~\ref{d.strategy}.) Writing $\omega = \tfrac{a}{a+b}$, we see then that
$$
\val(1,\lambda,N,a,b) =
\begin{cases}
\, 2 \omega + 1-\omega \, \, \, \textrm{when $\lambda_1 > \lambda_c$}  \, , \\
\, 2 \omega + (1-\omega) \val(\lambda_1,S) \, \, \, \textrm{when $\lambda_1 < \lambda_c$}  \,   .
\end{cases}
$$
By~[A2],  $\val(1,\lambda_1,S) = h(\lambda_1,S)$. 
Thus, when $\lambda_1 < \lambda_c$, $\val(1,\lambda_1,S)$ equals $h_2(\lambda_1,S)$.
We find then that
$$
\val(1,\lambda,N,a,b) =
\begin{cases}
\,  \omega + 1  \, \, \, \textrm{when $\lambda_1 > \lambda_c$}  \, , \\
\, 2 \omega +  2(1-\omega) \frac{\lambda_1^2}{1 + \lambda_1 + \lambda_1^2} \, \, \, \textrm{when $\lambda_1 < \lambda_c$}  \,   ,
\end{cases}
$$
with the right-hand formulas coinciding to specify $\val(1,\lambda,N,a,b)$ in the critical case $\lambda_1 =\lambda_c$, since $\tfrac{2 \lambda_c^2}{1 + \lambda_c + \lambda_c^2}  = 1$.

\begin{figure}[htbp]
\includegraphics[width=1\textwidth]{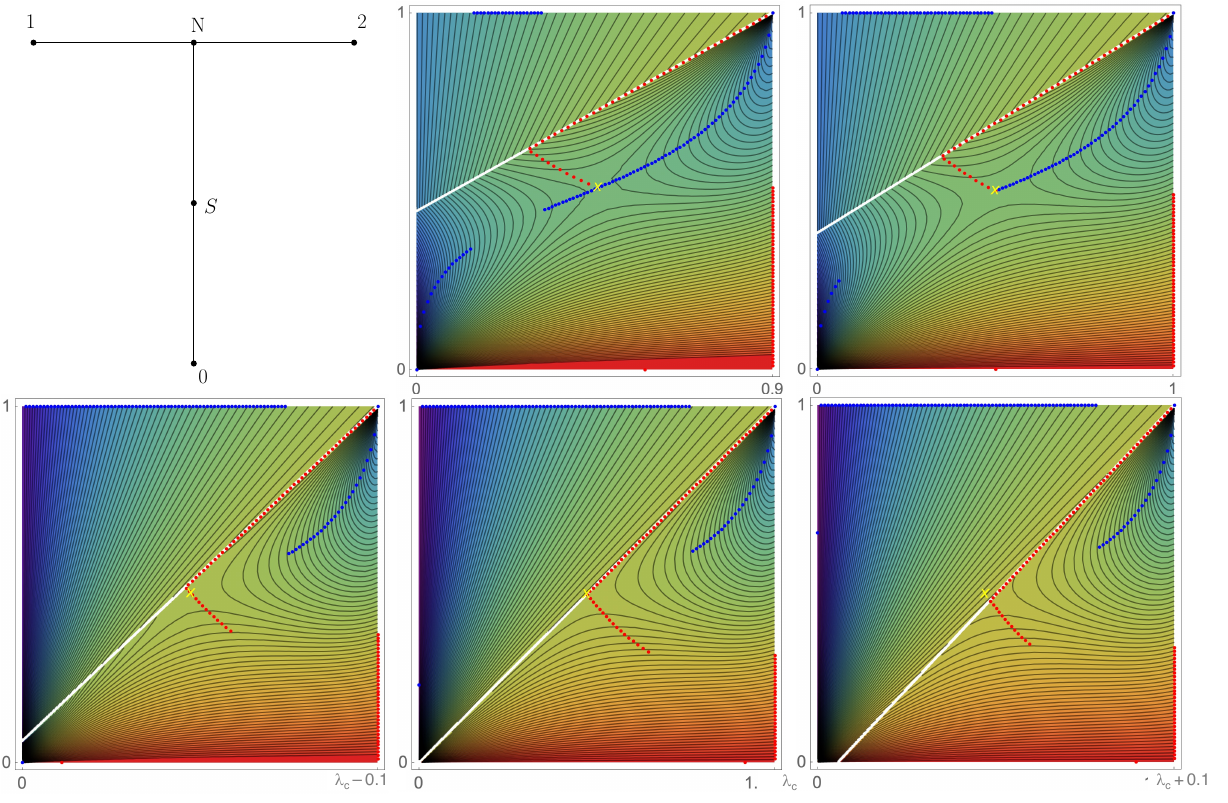}
\caption{The $T$ graph and five normal form contour plots  of $[0,\lambda] \times [0,1] \lora [0,2]: (a,b) \mapsto \val(1,\lambda,N,a,b)$, where $N$ is the north vertex  and the GloshConch's assumptons are in force.  The values of $\lambda$ in the upper middle and right are $0.9$ and~$1$; on the lower left, middle and right, they are $\lambda_c - 1/10$,
 $\lambda_c$ and $\lambda_c + 1/10$.
}\label{f.tgraph}
\end{figure}

The saddle point $(\lambda S,S)$ in $(a,b) \mapsto \val(1,\lambda,N,a,b)$ predicted in~(\ref{e.altstake})  is not uniquely defined, because the denominator of the right-hand side differs according to whether Mina nominates $1$ or $S$ as her move. 
(This means that the second criticism, in Subsection~\ref{s.criticismtwo}, is valid. In fact, the criticism is valid in a strong sense, because, as we will see, Mina is motivated to nominate a move that leads to counter evolutions that are not $\lambda$-biased infinity harmonic walks in the sense of Definition~\ref{d.theta}. This broadens the ambiguity in the definition of the denominator in~(\ref{e.altstake}).) If we assume that, at $v=N$, Mina either always nominates 1, or always nominates $S$, then there are two saddle point predictions: ($\lambda s,s)$ with $s=1$, or $s = (1 + \lambda + \lambda^2)(2\lambda + 1)^{-1}(\lambda + 1)^{-1}$, corresponding to Mina's nominations of $1$ or $S$, in the latter case working with the game on a line graph of length 3. The $s=1$ point is the go-for-broke location, where each player bids it all. Five normal form contour plots for the $T$ graph appear in Figure~\ref{f.tgraph}, with the latter of the two predicted saddle points appearing as a yellow cross in each sketch.  


The 
GloshConch
 fails in all five sketches, in the sense that neither of the predicted saddles is global in any case. (The first criticism voiced above is thus valid.) We discuss the lower-left plot, with $\lambda = \lambda_c - 0.1$, though several properties are shared with the other cases. The saddle $(\lambda s,s)$ is the yellow cross in the middle of the sketch. This saddle point is local but not global. Mina prefers the go-for-broke stake of one to the saddle-specified stake of~$s$   when Maxine stakes $\lambda s$. Should Mina stake~$1$, Maxine also prefers the go-for-broke~$\lambda$. But if Maxine stakes this, Mina has a response that Maxine lacks: she may stake the maximum~$1$ but hold an infinitesimal amount ${\rm d}x$ in reserve, and then nominate the move~$S$. If she wins the turn, she will automatically win the next turn also (since ${\rm d}x/0 = \infty$), and the game will end with the counter at vertex~$0$. Indeed, there is a discontinuity in  $[0,\lambda]\times [0,1] \lora [0,2]: (a,b) \mapsto \val(1,\lambda,N,a,b)$ at the northeast corner~$(\lambda,1)$ which Mina may exploit by playing the long game in this way. Thus the predicted saddle at $(\lambda,1)$ is not even a local saddle. Maxine reacts to a $1 - {\rm d}x$ stake by Mina by staking an amount qualitatively of the form $1 -  ({\rm d}x)^{1/2}$, and, if the players alternate in best responses, the cycle of increasing withholdings from maximum stakes would continue until the withheld amounts are of unit order; eventually, Mina once again goes for broke. In this example, then, Mina's luxury of move choice disrupts the predicted saddle at $(\lambda,1)$. There is no pure Nash equilibrium. (A parallel role for vertex $N$ in the $T$ graph is mentioned in Kleinberg and Tardos' work~\cite{KleinbergTardos} on network-bargaining problems, where it is noted that a player at $N$ who negotiates with neighbours over splitting rewards on intervening edges is in a strong position. Experimental designs with negotiating participants have been set up and performed~\cite{CookEmerson,CEGY} to test the strength of negotiators in such positions.)
 
\subsection{The Poisson game}\label{s.poisson}

We have seen that go-for-broke is a principal mechanism that disrupts the global saddle hope in the regular game. 
 The leisurely game is a variant that is designed largely in order to frustrate the efficacy of the go-for-broke strategy. Who would bet his life savings on the next turn when a shortly impending coin flip may dictate that no move will even take place? 
 We now introduce an idealized low-$\e$ limit of the leisurely game and carry out formal calculations that caricature the upcoming proof of Theorem~\ref{t.nashform}, in which the global saddle hope will be demonstrated for the leisurely game. In particular, we will see how assumption [A3] of Section~\ref{s.stakeperturb} and the stake formula~\eqref{e.stake} arise naturally.


The Poisson game on a boundary-payment graph $G = (V,E)$ takes place in continuous time. Let~$\mc{P}$ denote a unit intensity Poisson process on $[0,\infty)$. 
With $\stateofplay_0 = (\lambda,v) \in (0,\infty) \times \openmac$, Maxine's starting fortune is $\lambda_0 = \lambda$ and the starting counter location is $X_0 = v$.  
Maxine and Mina stake at time $t \in (0,\infty)$ at respective rates $a(t)$ and $b(t)$, where each of these processes is adapted with respect to the history of gameplay before time~$t$.
With instantaneous currency revaluation at time $t + {\rm d}t$ to ensure that Mina holds one unit at this time, we see then that Maxine's time $t+{\rm d}t$ fortune equals
\begin{equation}\label{e.lambdainf}
   \lambda(t + {\rm d}t) = \frac{\lambda(t) - a(t)dt}{1- b(t){\rm d}t} = \lambda(t) - \big( a(t) - \lambda(t) b(t) \big) {\rm d}t \, ,
\end{equation}
where we work formally with an infinitesimal calculus for which $({\rm d}t)^2$ equals zero.

A move happens when a Poisson clock rings: when $t$ reaches an element of $\mc{P}$. Indeed, this is the game $\game(0^+)$ of infinite leisure; except it has been speeded up so that moves occur at a unit-order rate. If a move happens at a given time $t$, Maxine wins with probability $\tfrac{a(t)}{a(t) + b(t)}$; otherwise, Mina does. Thus it is the present stake rates that dictate the move outcome. It is intuitively plausible that the strategy mimicry argument for Proposition~\ref{p.purevalue} in Section~\ref{s.strategymimicry}
would show that the game value\footnote{It may seem that the strategy mimicry argument in question would adapt to the Poisson case to show that the game value  $\val(\lambda,v)$ equals $h(\lambda,v)$ only if players are restricted to the use of pure strategies.  We elide the distinction between value with mixed or merely pure strategies in this heuristic discussion. In fact, even the question of how to define mixed strategy in the Poisson game may have several reasonable answers. To attempt nonetheless a rough summary:  the instantaneous monitoring of an opponent's strategy that the continuous time evolution in the Poisson game permits would make it impossible for a player to take advantage of the element of surprise that is a signature advantage of the use of mixed strategies. Tentatively, then, we may believe that the quantity $\val(\lambda,v)$ in the Poisson game does not change according to whether mixed strategies are permitted or prohibited.} $\val(\lambda,v) = \val (0^+,\lambda,v)$  is given by the biased infinity harmonic value in Definition~\ref{d.puregamevalue}:
 $$
 \val(\lambda,v) = h(\lambda,v) \,.
 $$
 We do not seek to make rigorous this assertion---there are challenges in formulating strategy spaces and resulting outcomes in continuous-time games---but rather focus on the prospect  that the more analytic nature of the Poisson game may help realize the global saddle hope.
 
Let $a,b \in [0,\infty)$. Consider a constrained game $\game(0^+,\lambda,v,a,b)$, where Maxine and Mina are obliged during $[0,{\rm d}t]$ to submit respective stake rates $a$ and $b$. They regain freedom at time ${\rm d}t$. Write $\val(\lambda,v,a,b)$ for the 
value of  $\game(0^+,\lambda,v,a,b)$. (We again omit a first argument $0^+$ in the notation for value.) Then we claim heuristically that 
\begin{equation}\label{e.phiclaim}
\val (\lambda,v,a,b) = \val(\lambda,v) +  \Phi(\lambda,v,a,b) {\rm d}t \, ,
\end{equation}
where  $\Phi:(0,\infty) \times \openmac \times (0,\infty)^2  \lora \R$ is given by 
\begin{equation}\label{e.phi}
 \Phi(\lambda,v,a,b) = - \val(\lambda,v) - (a-b\lambda) \val'(\lambda,v) + \tfrac{a}{a+b} \val(\lambda,v_+) + \tfrac{b}{a+b} \val(\lambda,v_-) 
\end{equation}
with $\val'(\lambda,v)  = \tfrac{\partial}{\partial \lambda} \val(\lambda,v)$, and where $v_+ = v_+(\lambda)$ and $v_- = v_-(\lambda)$ are neighbours of $v$ that maximize, or minimize, 
$\val(\lambda,\cdot) = h(\lambda,\cdot)$. (These neighbours may not be unique, but we only consider such quantities as $h(\lambda,v_+)$. They may depend on $\lambda$, and indeed this may pose problems, but we prefer to suppress this dependence in the notation, and elide this difficulty, in this heuristic discussion. A related assumption that we make is that the derivative  $\val'(\lambda,v)$ exists.)
To explain why~(\ref{e.phiclaim}) holds, consider the partition of the state space offered by the events of absence of a Poisson point in $[0,{\rm d}t]$; of the presence of such a point accompanied by a win for Maxine; and of such a presence alongside a win for Mina. By adding the values of the resulting subgames, we see that 
$$
\val (\lambda,v,a,b) 
  =  (1-{\rm d}t) {\rm Val}\big(\lambda({\rm d}t),v\big) + {\rm d}t \tfrac{a}{a+b} {\rm Val} \big( \lambda({\rm d}t), v_+ \big) +  {\rm d}t \tfrac{b}{a+b} {\rm Val}\big( \lambda({\rm d}t), v_- \big) \, .
$$
Using $(\rm d t)^2 = 0$, we may make three right-hand replacements $\lambda({\rm d}t) \to \lambda(0) = \lambda$. 
Also using $\lambda({\rm d} t) = \lambda - (a-\lambda b){\rm d}t$ (which is due to~(\ref{e.lambdainf})), we obtain~(\ref{e.phiclaim}).

The function $[0,\infty)^2 \lora \R: (a,b) \mapsto \Phi(\lambda,v,a,b)$ describes the infinitesimal reward in mean value for Maxine when the players jointly commit to $(a,b)$ rates for an instant from time zero. With shorthand $\Phi(a,b)=\Phi(\lambda,v,a,b)$, we thus expect that $(a_0,b_0)$, the stake pair to which the players infinitesimally commit under jointly optimal play (that is, in a pure Nash equilibrium), must be a saddle point, so that the first derivatives, given by
$$
\tfrac{\partial}{\partial a}
\Phi \big( a,b \big)  =  - \val'(\lambda,v) + \tfrac{b}{(a+b)^2}   \Big(  {\rm Val} \big( \lambda, v_+ \big)  -  {\rm Val} \big( \lambda, v_- \big)    \Big)   \, ,
$$
and
$$
\tfrac{\partial}{\partial b}
\Phi \big( a,b \big)  =   \lambda \val'(\lambda,v) - \tfrac{a}{(a+b)^2}   \Big(  {\rm Val} \big( \lambda, v_+ \big)  -  {\rm Val} \big( \lambda, v_- \big)    \Big)   \, ,
$$
must equal zero. Therefore,
\begin{equation}\label{e.threenonzero}
 a_0  \val'(\lambda,v)  =  \tfrac{a_0b_0}{(a_0+b_0)^2}   \Big(  {\rm Val} \big( \lambda, v_+ \big)  -  {\rm Val} \big( \lambda, v_- \big)    \Big) = b_0 \lambda  \val'(\lambda,v) \, .
\end{equation}
The quantity  ${\rm Val} \big( \lambda, v_+ \big)  -  {\rm Val} \big( \lambda, v_- \big) = h(\lambda,v_+) - h(\lambda,v_-)$ is typically positive (on root reward trees for example, this is a consequence of the upcoming Lemma~\ref{l.a}(4)), meaning that there is something to play for at this turn. Provided then $a_0b_0 \neq 0$---we omit an attempt to justify this detail---we have $a_0 = b_0 \lambda$, confirming the heuristic assumption [A3]: in jointly optimal play, the stake rates will satisfy 
 \begin{equation}\label{e.alambdab}
  a(t) = \lambda b(t) \, .
\end{equation}
This also implies, in~\eqref{e.phi}, that $\Phi(a_0,b_0) = 0$, which is also natural from the point of view that game value, being an expected later payment, should be a martingale under jointly optimal play, with the drift in~\eqref{e.phiclaim} vanishing.

As for the 
second derivatives, they equal
$$
\tfrac{\partial^2}{\partial a^2}
\Phi \big( a,b \big)  =  - \tfrac{2b}{(a+b)^3}   \Big(  {\rm Val} \big( \lambda, v_+ \big)  -  {\rm Val} \big( \lambda, v_- \big)    \Big)   
$$
and
$$
\tfrac{\partial^2}{\partial b^2}
\Phi \big( a,b \big)  =   \tfrac{2a}{(a+b)^3}     \Big(  {\rm Val} \big( \lambda, v_+ \big)  -  {\rm Val} \big( \lambda, v_- \big)    \Big)  \, .
$$

 We see then that
$\tfrac{\partial^2}{\partial a^2} 
\Phi \big( a,b \big)  < 0$ and 
$\tfrac{\partial^2}{\partial b^2} 
\Phi \big( a,b \big)  > 0$  for all $(a,b) \in (0,\infty)^2$. Thus, the saddle point $(a,b) \mapsto \Phi(a,b)$ at $(a_0,b_0)$, with $a_0 = \lambda b_0$, is a global minimax. We want to confirm the value of the stake rate~$b_0$ predicted by the right-hand side of~(\ref{e.stake}) after the removal of the factor of~$\e$ (a removal needed because we are concerned with stake {\em rates}). 
 By $a_0 = b_0 \lambda$ and the supposed non-vanishing of the factors in~(\ref{e.threenonzero}),  we find that
 $$
  \val'(\lambda,v) \, = \, \frac{1}{b_0(\lambda +1)^2}  \,  \Big(  {\rm Val} \big( \lambda, v_+ \big)  -  {\rm Val} \big( \lambda, v_- \big)    \Big) \, ,
 $$
whence
\begin{equation}\label{e.bzero}
b_0 = \frac{{\rm Val} \big( \lambda, v_+ \big)  -  {\rm Val} \big( \lambda, v_- \big)}{(\lambda +1)^2  \val'(\lambda,v)} \,,
\end{equation}
heuristically confirming the stake formula~(\ref{e.stake}).
 
 \subsection{The Poisson game by examples}\label{s.poissonexample}
 
 We reexamine in the Poisson case the three examples treated for the regular game in Section~\ref{s.localglobal}.
 
 \subsubsection{Example I: $(L_2,\sim,{\bf 1}_2)$}\label{s.exltwo}
 
 For the Poisson game on $(L_2,\sim, {\bf 1}_2)$ with $\stateofplay = (\lambda,1)$, we have $\val(\lambda,1) = \tfrac{\lambda}{\lambda +1}$ and 
 that  the global minimax $(a_0,b_0)$ equals $(\lambda,1)$, since~(\ref{e.bzero}) holds alongside $a_0 = \lambda b_0$.
 Since there is only one move on this graph, this choice of stake rates may be the closest the Poisson game comes to go-for-broke: the uncertainty of when the next Poisson clock will ring limits the players' tendency to stake at high rates. 

 \subsubsection{Example II:  $(L_3,\sim, {\bf 1}_3)$}
 
  We have that $\tfrac{\lambda}{1 + \lambda} \cdot  \val(\lambda,1) = \val(\lambda,2) =
 \tfrac{2 \lambda^2}{1 + \lambda + \lambda^2}$. 
 Thus, by~(\ref{e.phi}),
\begin{equation}\label{e.lthree}
 \Phi(\lambda,2,a,b) \, = \,  \frac{\lambda(\lambda + 1)}{1 + \lambda + \lambda^2} 
 \bigg( \frac{a}{a+b} - \frac{\lambda}{1 + \lambda} - \frac{a - b\lambda}{(1 + \lambda)^2} \Big( 2 + \tfrac{\lambda(1-\lambda)}{1 + \lambda + \lambda^2} \Big)
 \bigg) \, .
\end{equation} 
 The saddle point for~(\ref{e.lthree}) is global. It 
  equals $(\lambda b_0,b_0)$, where
$$
 b_0^{-1} \, = \, 2 +   \frac{\lambda(1-\lambda)}{1 + \lambda + \lambda^2} 
 \, .
$$

 \subsubsection{Example III: the $T$ graph}
 
 

Consider the graph~$T$ with the same boundary data as before. If $\lambda \neq \lambda_c$, then the updated value $\lambda({\rm d}t)$ in~(\ref{e.lambdainf}) lies on the same side of $\lambda_c$ as does $\lambda$. This permits the $(a,b)$-constrained 
value to be analysed by computing the relative values of $\lambda$ and $\lambda_c$ (provided that these are unequal),
rather than by the relative values of $\lambda_1$ and $\lambda_c$ that were seen in the regular case.


When $\lambda < \lambda_c$, Mina at $N$ nominates a move to $S$, and gameplay takes place on the $L_3$ copy $[0,S,N,2]$. 
 Thus, by~(\ref{e.lthree}),
$$
 \Phi(\lambda,N,a,b) \, = \,  \frac{\lambda(\lambda + 1)}{1 + \lambda + \lambda^2} 
 \bigg( \frac{a}{a+b} - \frac{\lambda}{1 + \lambda} - \tfrac{a - b\lambda}{(1 + \lambda)^2} \Big( 2 + \tfrac{\lambda(1-\lambda)}{1 + \lambda + \lambda^2} \Big)
 \bigg) \, .
 $$
When $\lambda > \lambda_c$, Mina at $N$ nominates a move to~$1$, seeking to end the game at the present turn, and the $L_2$ copy $[1,N,2]$ dictates the outcome. The formula for 
$\Phi(\lambda,N,a,b)$ now reads 
\begin{equation}\label{e.ltwo}
  \frac{a}{a+b} \, - \, \frac{\lambda}{1+\lambda} \, + \, \frac{a-b\lambda}{(1+\lambda)^2} \, .
\end{equation}
When $\lambda =\lambda_c$, the formula for $\Phi(\lambda,N,a,b)$  is the minimum of the two preceding forms. This is because it is Mina's prerogative at $N$ to nominate $1$ or $S$, and she does so to minimize subgame value.

\subsubsection{Contour plots and features seen in the Poisson examples}\label{s.sketchesandlessons}
We have noted that the $L_2$ saddle point $(a_0,b_0)$ equals $(\lambda,1)$. 
  
\begin{figure}[htbp]
\includegraphics[width=1\textwidth]{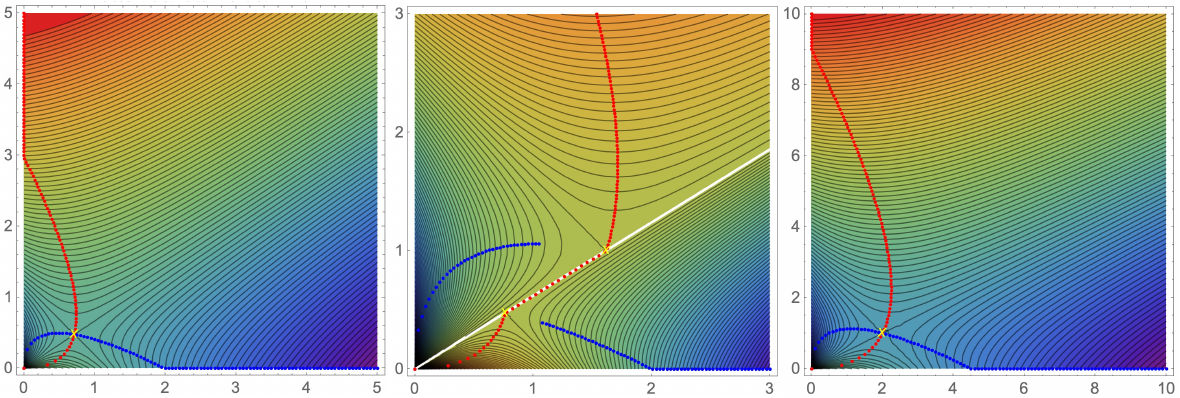}
\caption{Plotting the Poisson game: normal form contour plots of $[0,\infty)^2 \lora \R: (a,b) \mapsto \Phi(\lambda,N,a,b)$ for the $T$ graph. From left to right, the $\lambda$-values are $3/2$, $\lambda_c = 1.618 \cdots$, and $2$.
The left plot equally depicts  $\Phi(3/2,2,a,b)$ for the line graph~$L_3$, and the right one depicts  $\Phi(2,1,a,b)$ for~$L_2$.  The global minimax point $(\lambda b_0,b_0)$, with $b_0$ specified in~(\ref{e.bzero}), is marked with a yellow cross in the left and right plots.  The middle plot depicts the minimum of $\Phi(\lambda_c,2,a,b)$ for $L_3$ and   $\Phi(\lambda_c,1,a,b)$ for~$L_2$, and the two yellow crosses mark the predicted saddle point $(\lambda b_0,b_0)$ for these two functions.  The minimum reflects Mina's capacity to decide between nominating vertex $S$ or vertex~$1$.
}\label{f.poisson}
\end{figure}

Three normal form contour plots of $(0,\infty)^2 \lora \R: (a,b) \mapsto \Phi(\lambda,N,a,b)$ appear in Figure~\ref{f.poisson}. We indicate two basic features.

{\em Don't bet the house.} Note how play has a compact character: a player will cut her losses, staking at rate zero, if her opponent spends big.

{\em Don't sit the next one out.} The red and blue best response curves meet the origin tangentially to the coordinate axes in all three sketches. If the opponent offers no stake, the best response in each of the three cases is to offer a positive but infinitesimal stake.

 The complexity of the regular game plots in Figures~\ref{f.lthree} and~\ref{f.tgraph} has largely vanished in the off-critical plots, on the left and right in Figure~\ref{f.poisson}. The trace that remains is in the middle, critical, plot. Maxine's red response play runs the ridge between the lower-left $L_3$ yellow cross (the long game) to the upper-right $L_2$ yellow cross (the short game). We do not rigorously formulate the Poisson game in this article, but its study is well motivated, and we moot prospects for a rigorous inquiry in the final Section~\ref{s.directions}. To summarise our formal computations with and examples of this game: the global saddle hope appears to be realized under the assumptions (which we effectively  supposed) that 
 the set of  maximizers or minimizers of $h(\lambda,\cdot)$ among the neighbours of any given vertex is independent of $\lambda \in (0,\infty)$; and that the partial derivative of $h(\lambda,\cdot)$ in $\lambda$ exists. The prospects of the global saddle hope remain untested if one or other of these assumptions fails. Line graphs meet the assumptions but the $T$ graph does not, and further work would be needed to resolve these questions even in a graph as simple as the latter.

\subsection{A simple case where the stake formula is valid in the regular game}\label{s.asimplecase}
We have seen that the global saddle hope fails for the regular game for several simple graphs. Here we present an example where the stake formula~(\ref{e.stake}) is in fact valid for this game.

Let $n \in \nwithoutzero$. The half-ladder $\big(H_n,\sim,{\bf 1}_0\big)$ is the root-reward tree formed from the line graph $L_n$ and a collection of points $i^*$, $1 \leq i^* \leq n$, by attaching edges $i \sim i^*$, $i \in \intint{n}$. The reward vertex equals~$0$ and the field of open play is $\intint{n}$. Consider $\game(1,\lambda,n)$. The game has $n$ turns. Maxine must win at each turn so that the counter evolution $n \to n-1 \to \cdots \to 1 \to 0$ ends at $0$. A victory for Mina at the $i$\textsuperscript{th} turn will lead her to end the game with the move $n-i+1 \to (n-i+1)^*$ provided that the game is on at the start of this turn, whatever the value of $i \in \intint{n}$. It is natural to suspend currency revaluation in analysing the game and instead take the view that the stake strategies are indexed by possibly random vectors $\bar{x} = \big( x_i: i \in \intint{n} \big)$ (for Mina)  and  $\bar{y} = \big( y_i: i \in \intint{n} \big)$ (for Maxine) of non-negative entries that respectively sum to one and $\lambda$. At the start of the game, Mina places $x_i$ units and Maxine $y_i$ units against the vertex $i$ for each $i \in \intint{n}$. These deposits are the stakes to be played should the counter be at vertex~$i$ at the start of a turn. For $d \geq 0$, write $\hat{d}$ for the constant $n$-vector with sum~$d$, so that each entry is $d/n$. Let $\bar{y}$ be a non-random strategy for Maxine. Suppose that she plays against Mina's constant strategy $\hat{1}$. Note that the mean payment $M( \hat{1} , \bar{y} )$ equals $\prod_{i=1}^n \tfrac{y_i}{1/n+y_i}$. It follows readily from the strict concavity of $(0,\infty) \lora (0,\infty): z \mapsto -\log (1 + z^{-1})$ that the maximum over $\bar{y}$ of $M( \hat{1},\bar{y})$   equals $\lambda^n(1 + \lambda)^{-n}$ and is attained uniquely by $\bar{y} = \hat{\lambda}$. This is equally true if $\bar{y}$ is permitted to vary over random strategies for Maxine. Now suppose instead that Mina plays a strategy $\bar{x}$ against Maxine's constant $\hat{\lambda}$.
We have $M(  \bar{x},\hat\lambda) = \prod_{i=1}^n \tfrac{\lambda/n}{\lambda/n + x_i}$. 
Permit $\bar{x}$ to vary over non-random strategies for Mina. By the strict convexity of $z \mapsto -\log (1+z)$, this expression is minimized 
with value $\lambda^n(1+\lambda)^{-n}$ and a unique minimizer $\bar{x} = \hat{1}$. The conclusion holds equally when variation is permitted over random choices of~$\bar{x}$. These are the essentials of the proof of the next result.

\begin{proposition}
Let $n \in \nwithoutzero$. The value of the regular game $\game(1,\lambda,n)$ on the half-ladder $H_n$ equals $\lambda^n (1+\lambda)^{-n}$. The game has a unique Nash equilibrium. 
Let $i \in \intint{n}$. 
Should play under the Nash equilibrium continue to the start of the $i$\textsuperscript{th} turn, the counter will lie at $n+1-i$ at this time. The proportion of her remaining reserves that each player will then stake is equal to $(n+1-i)^{-1}$. The index of the final turn of the game has the law $G \wedge n$, where $G \geq 1$ is a geometric random variable of parameter $\PP(G =1)= (1+\lambda)^{-1}$.
\end{proposition} 
Note that the stake formula in its guise~(\ref{e.altstake}) is satisfied for $\game(1,\lambda,n)$ because~(\ref{e.altstake}) states that $\stake(1,\lambda,n) = \frac{\lambda^{n-1}(1+\lambda)^{-(n-1)}}{\sum_{i=1}^n \lambda^{i-1}(1+\lambda)^{-(i-1)} \cdot  \lambda^{n-i}(1+\lambda)^{-(n-i)}}= 1/n$. The summand 
is the product of the probability that play continues to the $i$\textsuperscript{th} turn and the value $\Delta(\lambda,n+1-i)$.

\subsection{Structure of the remainder of the article}

The $T$-graph at the point $\lambda=\lambda_c$ has illustrated that, if the set of optimal move nominations that a player may make 
changes with $\lambda$, 
complications arise for analysis, even in the idealized Poisson game. In working with root-reward trees, we are choosing to sidestep this difficulty, by using a framework in which such special values of $\lambda$ do not exist. In Section~\ref{s.tools}, we review 
the Peres-\v{S}uni\'c algorithm from~\cite{PeresSunic} for finding  biased infinity harmonic functions on graphs. We will establish that a related decomposition of the graph into paths (which has already been suggested in our discussion of the $T$ graph) has no dependence on $\lambda \in (0,\infty)$ in the case of root-reward trees. We will prove the $h$-differentiability Proposition~\ref{p.differentiable}
and several other needed facts about $h$ on root-reward trees.

Section~\ref{s.finitehorizon} concerns finite time-horizon games $\game_n(\e,\lambda,v)$ whose gameplay is in essence the truncation of that in $\game(\e,\lambda,v)$ to the first $n$ turns. Our finite-horizon games will be set up, in Section~\ref{s.basics}, so that each has value $h(\lambda,v)$, with no error needed to account for the finiteness of~$n$: see Lemma~\ref{l.tugnash} for the case of constant-bias tug-of-war; later, Theorem~\ref{t.leisurely.theorem}(1) will assert this for the stake-governed version.  As such, finite-horizon constant-bias tug-of-war becomes a tool for deriving the stake function formulas in Proposition~\ref{p.totvar}. Indeed,  we prove Proposition~\ref{p.totvar}  in Section~\ref{s.totvar}, by an argument concerning tug-of-war that is inspired by the heuristic perturbative identification of the stake function in the stake-governed version of the game seen in Section~\ref{s.stakeperturb}.

 Section~\ref{s.finitehorizonsolution} 
presents a fundamental element, the solution of finite-horizon games: it proves a counterpart Theorem~\ref{t.leisurely.theorem} 
for $\game_n(\e,\lambda,v)$  of Theorems~\ref{t.leisurely} and~\ref{t.nashform}.
This theorem shows that the global saddle hope is realized, and 
finds the location of the global saddle point, in the finite games. The result will be proved by the fundamental technique of backward induction,
for which we need to consider finite-horizon games; in this proof,
the special form of the putative stake function as a ratio of polynomials available for root-reward trees is harnessed in the inductive step to demonstrate that a global saddle point lies where it should.  
 By taking the horizon to infinity, we will derive Theorems~\ref{t.leisurely} and~\ref{t.nashform} in Section~\ref{s.finitetoinfinite}. We conclude in Section~\ref{s.directions} by reviewing some proof aspects and indicating several directions for developing theory and applications of stake-governed games.

\section{The Peres-\v{S}uni\'c decomposition and root-reward trees}\label{s.tools}

Here we recall, derive and harness properties of biased infinity harmonic functions.
In Section~\ref{s.ps} we review the Peres-\v{S}uni\'c algorithm for computing such functions and record useful consequences. In Section~\ref{s.stake.explicit},
we utilize a product formula for the biased functions to prove  the explicit formula for the stake function on line graphs given in Proposition~\ref{p.stake},
 as well as a further such result on root-reward trees.

\subsection{Peres and \v{S}uni\'c's algorithm for finding $\lambda$-biased infinity harmonic functions}\label{s.ps}

Peres and \v{S}uni\'c~\cite{PeresSunic} present an algorithm for computing biased infinity harmonic functions on a finite graph. 
Our treatment adds nothing essential to that offered by~\cite{PeresSunic}; but,  in our context of root-reward trees, it is convenient to recast their `PS' algorithm a little to a form we call mainline-sidings or `MS'. By a reduction to~\cite{PeresSunic}, Proposition~\ref{p.rhoslope} will show that the MS algorithm computes biased infinity harmonic functions. 
The mainline-sidings framework 
leads in Proposition~\ref{p.rootrewardinfinitybias} to a product formula for $\lambda$-biased infinity harmonic functions on root-reward trees;
and  a `journey-data' graph offers  a presentation of  a given root-reward tree in a simpler form, by identifying paths that are in essence equal under relevant symmetry. 
In Proposition~\ref{p.meanindep}, the journey-data graph will be exploited to  show that the stake function formula in Proposition~\ref{p.totvar} is well defined.  

\subsubsection{The MS algorithm}
A root-reward tree $T = (V,E,{\bf 1}_r)$ is akin to a network of  tracks that merge and lead to a grand central hub at~$r$.  
The mainline-sidings~[MS] partition is a
 simple decomposition of~$T$ along these lines. 

The network is composed of line graphs $L_k$. We begin by noting a formula for  $\lambda$-biased infinity harmonic functions on such graphs. 

\begin{definition}\label{d.hkell}
For $k \in \nwithoutzero$ and $\ell \in \llbracket 0 , k \rrbracket$, set  
$$
 H(\lambda,k,\ell) \, = \,
\begin{cases}
\, \,  \frac{1 - \lambda^{-(k-\ell)}}{1 - \lambda^{-k}}  & \textrm{for $\lambda \in (0,\infty) \setminus \{ 1 \}$}  \, , \\
\, \,  \frac{k-\ell}{k}  & \textrm{for $\lambda = 1$}  \,   .
\end{cases}
$$
\end{definition}

\begin{lemma}\label{l.hfact}
For such $k$ and $\ell$ as above, $H(\lambda,k,\ell)$ equals the $\lambda$-biased infinity harmonic
function $h(\lambda,\ell)$ for the line graph $L_k$ with root~$0$. 
\end{lemma}

{\bf Proof.} By solving~(\ref{e.h}) on $L_k$ recursively, it is easy to check that the system of equations has a unique solution, which is $H(\lambda,k,\ell)$.
\qed


\begin{definition}
Recalling the in- and out-distances from Definition~\ref{d.inout}, we set $\thespan:V \lora \N$ according to
 $\thespan(v) = \din(v) + \dout(v)$.

Define the close boundary  $B_C = \left\{ v \in \nonrootmac : d(r,v) = \min_{w \in \nonrootmac} d(r,w) \right\}$ to be that part of~$B^*$ closest to the root. 
The mainline $M = M(T)$ is the subgraph of $T$ induced by the union of paths $[r,u]$ for $u \in B_C$; $M$ is a root-reward tree with root~$r$. A siding is a connected component of the forest 
induced by the edge set $E(T) \setminus E(M)$.
 Each siding contains a unique element of $V(M)$ and may be viewed as a root-reward tree whose root is this element.  

The value of $\thespan(v)$ is independent of $v \in V(M)$, so we may define the mainline span by setting $\thespan(M)$ equal to the common value, which is simply $d(r,u)$ for $u \in B_C$.
More generally, whenever~$S$ is any root-reward subtree arising in the MS decomposition (as a siding or sub-siding), the same definition applied within $S$ yields a well-defined quantity $\thespan(S)$, equal to the span of its own mainline $M(S)$.
\end{definition}


Let $\lambda \in (0,\infty)$. The mainline-sidings~[MS] algorithm is an iterative procedure for determining a function $ \hms: V \lora [0,1]$.

First we set $\hms(r) = 1$ and $\hms(v) = 0$ for $v \in \nonrootmac$.
Then we set $ \hms:M \lora [0,1]$ on the mainline via 
\begin{equation}\label{e.ms.firststep}
\hms(v) = H\big(\lambda,\thespan(M),d(r,v) \big) \, \, \,  \textrm{for} \, \, \,  v \in V(M)
\end{equation}
 where $H$ is from Definition~\ref{d.hkell}. 
The mainline is removed and the graph fragments into sidings. Each siding $S$ has a root $r(S) \in M = M(T)$  on the mainline. 
Since $S$ is a root-reward tree, this graph  has its own mainline~$M(S)$.
The function~$\hms$ is now specified on the vertices of 
every siding mainline via
\begin{equation}\label{e.hmsassign}
 \hms(v) = \hms(v^\uparrow) \cdot H \big( \lambda, \thespan(S) , d(v^\uparrow,v) \big) \, \, \, \, \textrm{for} \, \, \, \, v \in V \big( M(S) \big) \, ,
\end{equation}
where $v^\uparrow$ denotes  $r(S) \in M(T)$. 

Next consider the graph formed by the removal of the union $\bigcup \big\{ M(S): S \in \mc{S} \big\}$
of the siding mainlines. Sub-sidings are the connected components of the resulting forest. 
Each contains a unique vertex in  $\bigcup \big\{ M(S): S \in \mc{S} \big\}$. Rooted at these vertices, the sub-sidings are root-reward trees. At the next step, values of $\hms$ are assigned to vertices in the mainlines of the sub-sidings~$S$ via the formula~(\ref{e.hmsassign}) whose first right-hand factor is an already ascribed value.

Sub-siding mainlines are then removed, and further offshoots in the network may result. The iteration continues with values assigned by the above formula at each step. The process stops when~$\hms$ has been defined on all of $V$. The mainline~$M(T)$ may be said to have depth zero; its sidings, depth one; the sub-sidings, depth two; and so on. The MS-depth of~$T$ set equal to
zero if $M(T)=T$ and to the 
  maximum siding depth otherwise.
  (The notion of MS-depth applies only to~$T$, 
  in contrast to the notion of depth, which qualifies  the mainlines of any siding or sub-siding of~$T$.)  
 In this way, the procedure terminates 
at the $k$\textsuperscript{th} stage after $\hms$-values have been ascribed if, for
every siding~$S$ of depth~$k$, the mainline $M(S)$  equals~$S$.
 We set the depth of any mainline set equal to that of the siding of which it forms part. Then the procedure's output is 
the {\it MS partition}: namely, 
 a partition of~$E(T)$
into edge-sets of mainlines 
with depths in~$\llbracket 0 , k \rrbracket$, with $k$ equal to the MS-depth. Every $v \in V$ except~$r$  lies in a unique mainline of which it is not the root, which mainline we denote by~$M_v$.

\subsubsection{The PS algorithm}

The Peres-\v{S}uni\'c~[PS] algorithm ascribes interpolated  $\lambda$-biased infinity harmonic values along paths in an iterative scheme in which paths are selected in descending order of  $\rho$-slope, a $\lambda$-determined real-valued function on paths which generalizes the standard notion of slope (for $\lambda=1$) to all $\lambda \in (0,\infty)$. What follows is a description of how the algorithm acts on root-reward trees.

The algorithm will output $\hps:V \lora [0,1]$ by iteratively specifying this function on increasing subsets of~$V$. 
Let $k:W\cup B^* \lora [0,1]$ for $W \subseteq V \setminus B^*$ be such a partially defined function, with $k(b)=0$ for all $b\in B^*$. For $w \in W$ and $b \in B^*$, the $\rho$-slope of $k$ on $[w,b]$ equals  
\begin{equation}\label{e.rhoslope}
 \rho \big( k, [w,b] \big) \, = \,
\begin{cases}
\, \,   \frac{k(w) (\lambda^{-1} - 1)}{\lambda^{-\ell(w,b)} - 1}  & \textrm{for $\lambda \in (0,\infty) \setminus \{ 1 \}$}  \, , \\
\, \,  \frac{k(w)}{\ell(w,b)}  & \textrm{for $\lambda = 1$}  \,   ,
\end{cases}
\end{equation}
where $\ell(w,b)$ denotes the length of the path~$[w,b]$ or the distance in $T$ between $w$ and $b$. 
Our $\rho$-slope is~\cite{PeresSunic}'s $r$-slope, but in any case the above formula is a special case of~\cite[Definition~$8$]{PeresSunic}.

Choose $u$ to be a vertex among $b \in \nonrootmac$ that  maximizes $\rho({\bf 1}_r, [r,b])$.
 Specify the function $\hps$ on the vertices in the path $[r,u]$ via $\hps(v) = H \big( \lambda, \ell(r,u), \ell(r,v) \big)$ for $v \in V[r,u]$. (This definition is consistent with $\hps(r)=1$.)  Remove $[r,u]$ from~$T$. Each connected component~$C$ is a root-reward tree rooted at a vertex $r(C)$ on~$[r,u]$ with non-root boundary $B^*(C) = B(C) \setminus \{ r(C) \}  \subseteq \nonrootmac =  B \setminus \{ r(T) \}$. Form a list of these components~$C$.
Next vary over pairs $(r',u')$ where $r' = r(C)$ for some component $C$ and $u' \in B(C)$, $u' \neq r'$. Make a choice for which  $\rho \big( k, [r',u'] \big)$ is maximum, where note that the value $\hps(r')$ entering into the above specification of the $\rho$-slope is an already assigned value of $\hps$.  Then set $\hps(v) = \hps(r') H \big( \lambda, \ell(r',u'), \ell(r',v) \big)$ for $v$ along $[u',r']$. 
Remove $C$ from the list of components. Then remove the path $[u',r']$ from $C$.
What results are certain connected components, which we may view as root-reward trees. Add these to the list of components. With this update, we continue with another variation to locate a pair of maximum~$\rho$-slope. The process continues until $\hps$ is defined on all of~$V = V(T)$. 
By~\cite[Claim~$18$]{PeresSunic}, the sequence of paths to which the PS algorithm assigns value has non-increasing $\rho$-slopes.
One may say that there is a progressively falling level for $\rho$-slope admission, in which a path whose endpoints have already been assigned value, so that the path's $\rho$-slope is well defined, is admitted (and its internal vertices assigned value) as soon as the admission value falls past the path's $\rho$-slope (with a tie-break rule applied, should other candidate entrants have the same $\rho$-slope).

(This is a description of the algorithm stated after~\cite[Definition~$17$]{PeresSunic}. Note that Step~$6$ there is a treatment of `$\rho$-slope zero' paths at the procedure's end. This step does not occur in our case.)

The resulting $\hps:V \lora [0,1]$ is a $\lambda$-biased infinity harmonic function on $(V,E,{\bf 1}_r)$ by~\cite[Claim~$19$]{PeresSunic}; and it is the only such, as noted at the start of~\cite[Section~$1.5$]{PeresSunic}.

\subsubsection{Reducing to the Peres-\v{S}uni\'c treatment}


By reducing to the analysis of~\cite{PeresSunic}, we will  find that, in the root-reward context, $\hms$ equals the $\lambda$-biased interpolation of boundary data. Proposition~\ref{p.rhoslope} to this effect builds on the relation shared by PS and MS stated in Lemma~\ref{l.hmsps}.

Let $S$ be a siding of depth one. Then $S$ is a root-reward tree with root $r(S)$. The MS and PS algorithms may act on this tree, interpolating the boundary data ${\bf 1}_{r(S)}:B(S) \lora \{ 0,1\}$.
We record the output functions in the form $\hms(\cdot,S)$ and $\hps(\cdot,S)$. These are maps from $V(S)$ to $[0,1]$.

\begin{lemma}\label{l.hmsps}
Let $S$ be a siding of depth one.  Write $v^\uparrow = r(S)$ for $v \in V(S)$. The functions $\hms$ and $\hps$ both satisfy the relation $h(v) = h(v^\uparrow) h(v,S)$ for such~$v$. 
\end{lemma}

{\bf Proof.} For $h = \hms$, the definition made by (\ref{e.hmsassign}) for $v \in M(S)$ satisfies the stated relation, because $h(v,S) = H \big( \lambda, \thespan(S) , \ell(v^\uparrow,v) \big)$ holds in view of~(\ref{e.ms.firststep}). A simple induction then serves to extend the validity of the relation to all of $V(S)$.
For $h = \hps$, consider how the PS algorithm operates from the moment that $\hps\big(r(S)\big)$ is assigned by the introduction of a connecting path in~$M(S)$. Monitor the algorithm only at times when  it adds a path between endpoints in the siding~$S$. The process of iteratively assigning value along these paths is identical to that of the algorithm were it to be applied solely to~$S$, except that the boundary condition is $\hps\big(r(S)\big){\cdot \bf 1}_{r(S)}: B(S) \lora [0,1]$ in place of ${\bf 1}_{r(S)}$. The relation for $\hps$ results. \qed

\begin{proposition}\label{p.rhoslope}
Let $T = (V,E,{\bf 1}_r)$ be a root-reward tree, and  
let $\lambda \in (0,\infty)$. 
The function $\hms: V \lora [0,1]$ computed by the above procedure is the unique $\lambda$-biased infinity harmonic function that interpolates ${\bf 1}_r$ on~$T$.
\end{proposition}

{\bf Proof.}
The task 
reduces to showing that $\hps = \hms$ in view of $\hps(\cdot) = h(\lambda,\cdot)$. We will obtain this by an inductive argument. 
Our inductive hypothesis~IH($\ell$) is that $\hps = \hms$ on the subtree induced by the mainlines of 
 depth at most $\ell$ in any root-reward tree $T$. 

First, the base case, $\ell=0$. This is equivalent to proving that,  after a certain number of steps in the PS algorithm, the union of admitted paths equals $M(T)$, and $\hps$ and $\hms$ agree there. 


Consider ${\bf 1}_r:B \lora \{ 0,1 \}$, vanishing on $B^*$, and recall the notion of $\rho$-slope from~(\ref{e.rhoslope}). Let $\rho_0$ denote the maximum value that $\rho \big( {\bf 1}_r, [r,b] \big)$ assumes for $b \in B_*$. Clearly, maximizing paths run from $r$ to the close boundary~$B_C$.

Let $U$ be any subset of $V\big(M(T)\big)\cap B^* = B_C$, and let $T_U$ be the tree induced by the union of paths $[r,u]$ as $u$ ranges over $U$. 
Consider also the function $k=k_U:V(T_U) \lora [0,1]$ given by the $\lambda$-biased interpolation of ${\bf 1}_r: \{ r \} \cup U \lora \{ 0,1 \}$ on this set.
A {\em leg} is a path from an element of $V(T_U)$ to $B^*$ that contains no edge in $T_U$. 
All out paths in $T_U$ that start at $v$ and end in $B^*$ have a common length. 
A leg $[v,b]$ is {\em short} if its length is this shared value; otherwise it is {\em long}.

We may speak of the $\rho$-slope $\rho \big( k, [v,b] \big)$ of any leg, because $k$ is already defined at the leg's endpoints. We {\em claim} that, for an arbitrary such subset~$U$, the $\rho$-slope of any leg is at most $\rho_0$, with this value assumed if and only if the leg is short. 

Indeed,  let $Q = [q,b]$ be a short leg. Let $P =[ r,b']$ be a path in $T_U$ that contains~$q$ and runs to some $b' \in B_C$. 
  The $\rho$-slope of $P$ is~$\rho_0$.
  Write $P'$ for the subpath $[q,b']$ of $P$.
  Since $k$ is specified on~$P$,  we may speak of the $\rho$-slope of $Q$ and~$P'$. These values are equal because 
$\rho$-slope is determined by endpoint values and path length, data which is shared by the two paths. And the $\rho$-slope of $P'$ equals that of $P$ because $k$ on~$P$ is $\lambda$-biased infinity harmonic and $P'$ is a subpath of~$P$: see~\cite[Proposition~$15$]{PeresSunic}. So the $\rho$-slope of $Q$ equals that of~$P$, namely~$\rho_0$. Now suppose instead that $Q = [q,b]$ is a long leg. Let $P = [r,b']$ be as above. Since $\ell(q,b) > \ell(q,b')$, the $\rho$-slope of $Q$ is less than $\rho_0$ in view of the formula~(\ref{e.rhoslope}).  
This completes the proof of the claim.

Now, the claim easily implies {\rm IH}($0$). At the start of the PS algorithm, a path $[r,u]$ with $u \in B_C$ is added, with the $\lambda$-biased values along this path in $M(T)$. Taking $U=\{u\}$, the claim shows that the next added path also lies in $M(T)$. And so on, until all of $M(T)$, adorned with $\lambda$-biased values, has been added. These values are also those for $\hms$ on $M(T)$ as given in~\eqref{e.ms.firststep}, hence {\rm IH}($0$) holds.

We now show that  {\rm IH}($i$), $i \in \llbracket 0,\ell-1 \rrbracket$, imply  {\rm IH}($\ell$).  
Let $T_\ell$ be the union of all mainlines of depth at most $\ell$, and suppose that $w \in V$ is on a mainline of depth $\ell$. So $w \in V(S)$ for some siding $S$, and $S\cap T_\ell$, viewed as a root-reward tree itself, rooted at $w^\uparrow \in M(T)$, has MS-depth $\ell-1$. Then we have
$$
 \hms(w) = \hms(w^\uparrow) \hms(w,S) =  \hps(w^\uparrow) \hps(w,S) = \hps(w) = h(\lambda,w) \, ,
$$ 
which establishes  {\rm IH}($\ell$); the successive equalities are due to Lemma~\ref{l.hmsps} for~$\hms$;  {\rm IH}($i$) 
for $i \in \{ 0,\ell-1\}$;  
 Lemma~\ref{l.hmsps} for~$\hps$; and $\hps = h(\lambda,\cdot)$ by~\cite[Claim~$19$]{PeresSunic}. 
This completes the proof of  Proposition~\ref{p.rhoslope}. \qed

 

\subsubsection{Journey data and the product formula}

A path from the root to a given vertex~$v$ may leave the mainline and enter a siding, and 
it may be further sidetracked by
 entering sub-sidings of progressively greater depth. Such moments of departure may naturally be viewed as junctions on a journey along the path. We present a recursion for $\lambda$-biased functions in terms of junction points; formalize the list of such points as the `journey data' of the destination~$v$; present the product formula Proposition~\ref{p.rootrewardinfinitybias} in terms of journey data; and record some ramifications. 

For $v \in V \setminus \{ r\}$, recall that $M_v$ denotes the unique mainline in the MS algorithm of which $v$ is a non-root element.
We now set $v^\uparrow_{\rm loc}$ equal to $r(M_v)$. Where $v^\uparrow$ has been specified to be the element on the principal or depth-zero mainline~$M(T)$ which roots the depth-one siding that contains $v$, the new definition is a local manifestation of this notion which refers to the root of~$M_v$.

\begin{lemma}\label{l.rootrewardinfinitybias.developed}
Let $(V,E,{\bf 1}_r)$ be a root-reward tree. 
The $\lambda$-biased infinity harmonic function $h(\lambda,\cdot): V \lora [0,1]$ specified by~(\ref{e.h}) with $p = {\bf 1}_r$ satisfies the following recursion. 
Let $v \in \openmac =  V \setminus B$. Then
$$
 h(\lambda,v) \, = \, h \big( \lambda,v^\uparrow_{\rm loc} \big)  \cdot H\Big( \lambda, \thespan(M_v) , d \big( v^\uparrow_{\rm loc},v \big) \Big) \, .
$$
\end{lemma}
{\bf Proof.} By Peres-\v{S}uni\'c, $h(\lambda,v) = \hps$. By Proposition~\ref{p.rhoslope}, $\hps = \hms$. So it suffices to demonstrate the recursion for~$\hms$. 
But this recursion is simply~(\ref{e.hmsassign}). \qed

For  $v \in V \setminus \{ r \}$, recall that $v_+$ is the unique neighbour of~$v$ that is closer to $r$ than $v$ is.

\begin{definition}\label{d.journeydata}
 Let $v \in V \setminus \{ r \}$. 
\begin{enumerate}
\item
A vertex $u \in V$ is a junction of $v$  if either
\begin{itemize}
\item $u$ equals $r$ or $v$; or
\item  $u$ is a vertex (not equal to~$r$) on the path $[r,v_+]$ 
such that $\thespan(u_-) > \thespan(u)$, where~$u_-$ denotes the successor of $u$ on $[r,v]$. 
 \end{itemize} 
\item  Let $\big\{ j_i: i \in \llbracket 0,k+1 \rrbracket \big\}$, with $k \in \N$, denote the junctions of~$v$ in the order that they are encountered along $[r,v]$, with $j_0 = r$ and $j_{k+1} =v$. Set $d_i = \ell(j_i,j_{i+1})$ and $s_i = \thespan\big(M_{j_{i+1}}\big)$ for $i \in \llbracket 0,k \rrbracket$.
The journey data of $v$ is the list {\rm JD}($v$) $= \big\{  (s_i,d_i): i \in \llbracket 0,k \rrbracket \big\}$.
\end{enumerate}
\end{definition}

See the left picture in Figure~\ref{f.essencetree} below for an example of journey data.

A product formula for biased infinity harmonic functions may be given in terms of journey data.

\begin{proposition}\label{p.rootrewardinfinitybias}
Let $(V,E,{\bf 1}_r)$ be a root-reward tree.
The $\lambda$-biased infinity harmonic function $h(\lambda,\cdot): V \lora [0,1]$ specified by~(\ref{e.h}) with $p = {\bf 1}_r$ takes the form
$$
 h(\lambda,v) \, = \, \prod_{i=0}^k H(\lambda,s_i,d_i) \, ,
$$
where $\big\{ (s_i,d_i): i \in \llbracket 0, k \rrbracket \big\}$ is the journey data of $v \in V \setminus \{ r  \}$. (Naturally, $h(\lambda,r)=1$.)
\end{proposition}
{\bf Proof.} Note that $j_k = v^\uparrow_{\rm loc}$ and $j_{k+1} = v$, so that $s_k = \thespan(M_v)$ and $d_k = d\big( v^\uparrow_{\rm loc},v \big)$.
The formula is obtained by induction on the depth of  $M_v$, with the recursion in Lemma~\ref{l.rootrewardinfinitybias.developed} thus
enabling the inductive step. \qed

\begin{corollary}\label{c.hinc}
Let $(V,E,{\bf 1}_r)$ be a root-reward tree. For $v \in V$, $h(\cdot,v): [0,\infty) \lora [0,1]$ is increasing. 
\end{corollary}

{\bf Proof.} This follows directly from Proposition~\ref{p.rootrewardinfinitybias}. \qed

Another consequence of the proposition is the existence of the derivative of $h(\lambda,v)$:

{\bf Proof of Proposition~\ref{p.differentiable}.} Journey data depends on the tree but not on~$\lambda$. Thus $h(\lambda,v)$ is a finite product of factors  $H(\lambda,s,d)$ with fixed $s$ and $d$ which is differentiable in $\lambda$. \qed


We now introduce a notation with similarities to $\mc{V}_-(v)$ from Definition~\ref{d.inout}.
\begin{definition}\label{d.vn}
Let  $v \in V \setminus \{ r \}$. Write 
$$
 \mc{V}_-(\lambda,v) \,= \, \left\{ \, w \in V: w \sim v \, , \, h(\lambda,v) = \min_{u \sim v} h(\lambda,u) \, \right\} \, .
$$
For $v \in V \setminus B^*$, where recall $B^* = B \setminus \{ r \}$,
let $\mc{V}_+(\lambda,v)$ be the counterpart set where maximum replaces minimum on the right-hand side.
\end{definition}

\begin{lemma}\label{l.a}
Let $v \in V(T)$ have journey data  $\big\{ (s_i,d_i): i \in \llbracket 0, k \rrbracket \big\}$, and let $\lambda \in (0,\infty)$.
\begin{enumerate}
\item The set $\mc{V}_+(\lambda,v)$ is a singleton whose element is $v_+$, the parent of $v$. The data {\rm JD}($v_+$) is formed from {\rm JD}($v$) $= \big\{ (s_i,d_i): i \in \llbracket 0, k \rrbracket \big\}$ by removing the final pair $(s_k,d_k)$ if $d_k=1$; if $d_k > 1$, by the replacement of $d_k$ by $d_k - 1$.
\item If $v \in O = V \setminus B$, then for all $\lambda \in (0,\infty)$ and $w \in \mc{V}_-(\lambda,v) $, {\rm JD}($w$) is formed from {\rm JD}($v$) by the replacement $d_k \to d_k+1$. 
\item For  $v \not\in B^*$, the sets  $\mc{V}_-(v)$  and  $\mc{V}_-(\lambda,v)$ are equal. 
\item For $v \neq r$, $h(\lambda,v_+) > h(\lambda,v)$.
\end{enumerate}
\end{lemma}
{\bf Proof: (1).} In the MS algorithm, $\hms$-values decrease as distance increases from the root, because this property holds whenever new values are recorded. But $h(\lambda,\cdot) = \hms$ by Proposition~\ref{p.rhoslope}. So $v_+$ is the unique neighbour of $v$ that is accorded a higher $h(\lambda,\cdot)$-value.

{\bf (2).} This again follows from Proposition~\ref{p.rhoslope}, but let us give a direct computational proof that does not rely on \cite{PeresSunic}. We work with the formula for $h(\lambda,\cdot)$ given in Proposition~\ref{p.rootrewardinfinitybias}, using the function $H(\lambda,\cdot,\cdot)$ from Definition~\ref{d.hkell}. Suppose that $v \in O$ has journey data with final pair $(s_k,d_k) = (q,\ell)$. The parent $v_+$ cannot minimize $h(\lambda,\cdot)$ among neighbours of~$v$ by the preceding part. Let $u$ and $u'$ be children of $v$ that respectively do and do not lie in the mainline of which $v$ is a non-root vertex. Note that JD($u$) is formed from JD($v$) by the replacement of the final pair $(q,\ell)$ with $(q,\ell+1)$, while JD($u'$) is formed from JD($v$) by adding a final pair $(q',1)$. 
Note that $q' > q-\ell$, since in the MS algorithm a new mainline branching from $v$ must have span strictly larger than the remaining outward distance $q-\ell$ along the mainline containing $v$ as a non-root vertex.
By Proposition~\ref{p.rootrewardinfinitybias}, the condition that $h(\lambda,u) < h(\lambda,u')$ is equivalent to 
$$
  H(\lambda,q,\ell) \cdot H(\lambda,q',1) > H(\lambda,q,\ell+1)
$$
which in turn is equivalent to $t(q',\lambda) > t(q-\ell,\lambda)$ where $t(a,\lambda) = \frac{1 - \lambda^{1-a}}{1-\lambda^{-a}}$. Now $q'$ strictly exceeds $q-\ell$, and the latter quantity is positive, because $v \in O$ is not a boundary vertex. Hence, the requisite comparison of $t$-quantities follows from the claim that for $\lambda \in (0,\infty) \setminus \{ 1\}$, the map $(0,\infty) \lora \R: a \mapsto t(a,\lambda)$ is increasing. The claim is validated by noting that  $t(a,\lambda) = 1 - \frac{\lambda -1}{\lambda^a - 1}$, with $\lambda^a - 1$ positive and increasing for $\lambda > 1$, and negative and decreasing for $\lambda \in (0,1)$. 

{\bf (3,4).} These assertions have been obtained in the preceding proof. \qed
  
 \subsubsection{The journey-data tree and the uniqueness of lazy biased random walk thereon}\label{ss.JDtree}

Let $T = (V,E, {\bf 1}_r)$.  
Write $\Phi$ for the map  defined on $V$ that sends any vertex $v \in V$ to its journey data. This map naturally extends to be defined on~$E$. The triple
$\big(\Phi(V),\Phi(E), {\bf 1}_{\Phi(r)}\big)$ is a root-reward tree that we call the {\em journey-data tree} of~$T$ and denote by $\tjdg  = (\vjdg,\ejdg,{\bf 1}_{r(\tjdg)})$. Thus~$\Phi$ is a graph homomorphism from $T$ to $\tjdg$. A description of $\ejdg$ has been given in Lemma~\ref{l.a}.



The MS partition of a root-reward tree is called {\em simple} if all of its mainlines are copies of a line graph. 

 Here is a useful inference of Lemma~\ref{l.a}.
\begin{corollary}\label{c.uniqueplay}
Let $T = (V,E,{\bf 1}_r)$ denote a root-reward tree.
 Then $\tjdg$ is simple and the set $\mc{V}_-(\lambda,v) \subset \vjdg$ is a singleton set for all $v \in \vjdg$ and $\lambda \in (0,\infty)$. The element is $v_-$, the unique child of $v$ such that $\dout(v_-) = \dout(v) -1$.
\end{corollary}
{\bf Proof.} By Lemma~\ref{l.a}(2) applied to $\tjdg$, all elements of $\mc{V}_-(\lambda,v) \subset \vjdg$ are described by a single piece of journey data. Thus, $\tjdg$ is simple. \qed

\begin{figure}[htbp]
\begin{center}
\SetLabels
(0.11*0.34) \color{red}{$v=j_3$}\\
(0.15*0.54) $j_2$\\
(0.33*0.48) $j_1$\\
(0.32*0.69) $r=j_0$\\
(0.01*0.63) \color{blue}{$w$}\\
(0.78*0.34) \color{red}{$\Phi(v)$}\\
(0.78*0.06) \color{blue}{$\Phi(w)$}\\
(0.96*0.7) $\Phi(r)$\\
\endSetLabels
\centerline{
\AffixLabels{
\includegraphics[width=0.75\textwidth]{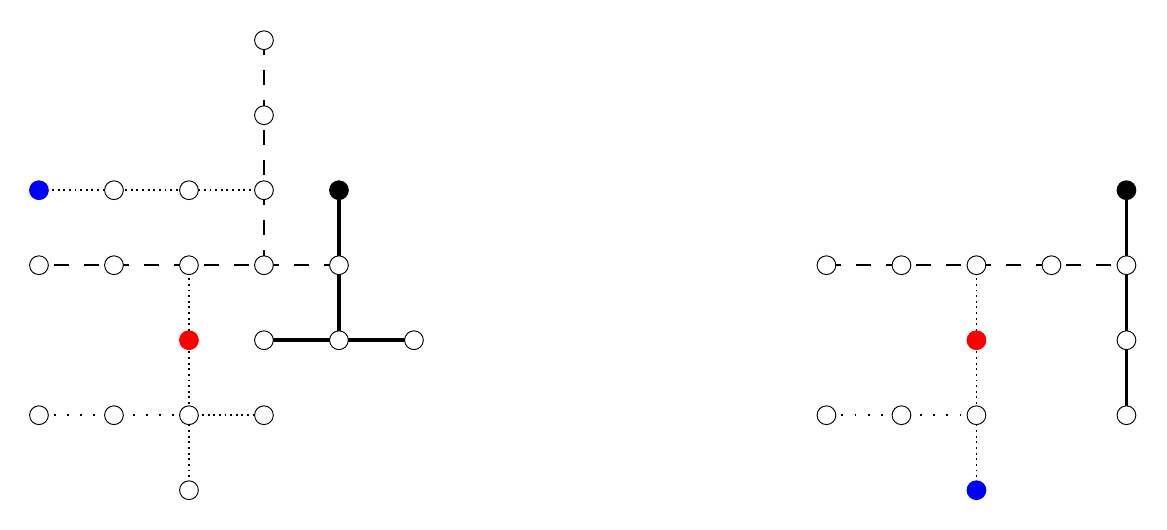}
}}
\caption{The right root-reward tree is the journey-data tree 
of the left one. 
Each tree has depth three;
the depth-zero mainlines are emboldened, the depth-one mainlines are dashed, the depth-two and three mainlines are densely and sparsely dotted, respectively. The junctions of vertex $v$ on the left are indicated. The vertices $v$ and $\Phi(v)$ have journey data $\big\{ (3,1),(4,2),(3,1) \big\}$ in their respective trees. For Maxine to win from $v$ on the left, she must triumph in three rounds: by reaching $j_2$, $j_1$ and~$j_0$. For $\lambda \in (0,\infty)$, the $\lambda^{-1}$-slope of the path in the MS
decomposition on which the counter is found rises with each round at which
Maxine succeeds:
 on root-reward trees, tug-of-war for Maxine resembles a gameshow where a wrong answer means the contestant's exit, and the pitch of whose music rises by a semitone with each passing round.}\label{f.essencetree}
 \end{center}
\end{figure}

Recall from Definition~\ref{d.theta} the $(1-\e)$-lazy $\lambda$-biased random walks $X_\theta: \llbracket 0,F \rrbracket \lora V$ on root-reward trees, indexed by $\theta \in \Theta$.
    In the journey data tree, $\mc{V}_-(v)$ is a singleton~$\{ v_-\}$ for each $v \in \openmac$
    by Lemma~\ref{l.a}(3) and Corollary~\ref{c.uniqueplay}.
     The space $\Theta$ is thus also a singleton in this case. 
 We may denote the process~$X_\theta$, where $\Theta = \{ \theta \}$, by $X_{\rm JD}$. To be explicit, 
$X_{\rm JD}: \llbracket 0,F_{\rm JD} \rrbracket \lora V_{\rm JD}$ is  the Markov process such that  $ X_{\rm JD}(0) =v \in V_{\rm JD}$  and, for $k \in \nwithzero$ such that $X_{\rm JD}(k) \in \openmac$,
$$
 X_{\rm JD}(k+1) =
\begin{cases}
\,  X_{\rm JD}(k)   & \textrm{with probability $1-\e$}  \, , \\
\,  \big( X_{\rm JD}(k)\big)_+  & \textrm{with probability $\e \tfrac{\lambda}{1 + \lambda}$}  \, , \\
\, \big(X_{\rm JD}(k)\big)_-  & \textrm{with probability $\e \tfrac{1}{1 + \lambda}$}  \,   .
\end{cases}
$$
In Figure~\ref{f.essencetree}, a root-reward tree appears alongside its journey-data tree.
 Note that although on the path in the latter tree from $\Phi(w)$ to $\Phi(r)$ there are branches without counterparts on the path from~$w$ to $r$ in the original tree, the biased random walk $X_{\rm JD}$ from $\Phi(w)$ is not influenced by these.
 
  \begin{proposition}\label{p.meanindep}
 Let $(V,E,{\bf 1}_r)$ be a root-reward tree.
 \begin{enumerate}
 \item
  The processes $X_\theta$, $X_\theta(0) = v \in \openmac$, are such that the $V_{\rm JD}$-valued processes $\Phi \circ X_\theta : \llbracket 0, F_\theta \rrbracket \lora V_{\rm JD}$ are each equal in law to  $X_{\rm JD} : \llbracket 0, F_{\rm JD} \rrbracket \lora V_{\rm JD}$, $X_{\rm JD}(0) = \Phi(v)$. 
  \item 
  The value
  $$
  \E \, \totvar \big(\e,\theta,\lambda,v \big) \, = \, \E \sum_{i=0}^{F_\theta-1} \bigg( h\Big(\lambda,\big(X_\theta(i)\big)_+\Big) - h\Big(\lambda, \theta\big(\lambda,X_\theta(i),i\big)\Big) \bigg)
  $$ 
  is independent of $\theta \in \Theta$. 
  \end{enumerate}
 \end{proposition}

 Proposition~\ref{p.meanindep}(2) resolves the second criticism of the heuristic perturbation argument, which was voiced in Subsection~\ref{s.criticismtwo}: despite the apparent dependence of the total variation on the choice of $\theta \in \Theta$, its expectation, which is  the denominator in~(\ref{e.altstake}), is well-defined and independent of this choice.

  {\bf Proof: (1).}  Lemma~\ref{l.a} implies that $\Phi \circ X_\theta$  evolves according to the rules specified there, independently of which  $h(\lambda,\cdot)$-minimizing neighbour is chosen by $\theta$
  when a move away from the root is made.
  
{\bf (2).} Let $\theta \in \Theta$.
Since $\Phi \circ X_\theta$ and $X_{\mathrm{JD}}$ have the same transition probabilities and initial states, they may be coupled to agree almost surely.
By Lemma~\ref{l.a}(1), we see then that, almost surely,
$$
h\Big(\lambda,\big(X_\theta(i) \big)_+ \Big) \, = \,
h\Big(\lambda,\big(X_{\rm JD}(i) \big)_+\Big) \, .
$$
By Lemma~\ref{l.a}(2) and Corollary~\ref{c.uniqueplay}, we also have that
$$
h\Big(\lambda,\theta\big(\lambda,X_\theta(i),i \big) \Big) \, = \, h\Big(\lambda,\big( X_{\rm JD}(i)\big)_-\Big) 
$$
almost surely. Thus,
  $$
  \E \, \totvar \big(\e,\theta,\lambda,v \big) \, = \, \E \sum_{i=0}^{F_{\rm JD}-1} \bigg( h\Big(\lambda,\big(X_{\rm JD}(i) \big)_+\Big)  -h\Big(\lambda,\big( X_{\rm JD}(i)\big)_-\Big) 
 \bigg) 
  $$ 
whatever the value of $\theta \in \Theta$. And so we obtain Proposition~\ref{p.meanindep}(2). \qed

\subsection{The stake function computed explicitly}\label{s.stake.explicit}

Here is the explicit formula for the stake function that we promised at the end of the introduction. 

\begin{theorem}\label{t.stakeformula}
For $\lambda \in (0,\infty) \setminus \{ 1\}$ and $\ell \in \nwithoutzero$, set $\Psi(\lambda,\ell) = \tfrac{\lambda^{-\ell}}{1 - \lambda^{-\ell}}$.
Let $(V,E,{\bf 1}_r)$ be a root-reward tree. 
Let $v \in V$ have journey data $\big\{ (s_i,d_i): i \in \llbracket 0, k \rrbracket \big\}$ for $k \in \nwithzero$. 
\begin{enumerate}
\item When $\lambda \in (0,\infty) \setminus \{ 1\}$, we have that
$$
\stake(1,\lambda,v) \, = \, 
 \frac{(\lambda - 1) \Psi(\lambda,s_k-d_k)}{(\lambda+1)\sum_{i=0}^k \Big(  (s_i - d_i) \Psi(\lambda,s_i - d_i) - s_i  \Psi(\lambda,s_i) \Big)} \, .
$$
\item When $\lambda =1$,
$$
 \stake(1,1,v) \, = \, \bigg( (s_k - d_k) \sum_{i=0}^k d_i \bigg)^{-1} \, .
$$
\end{enumerate}
\end{theorem}

Before attempting the proof, we derive the stake function on line graphs.

{\bf Proof of Proposition~\ref{p.stake}.} Here, $h(\lambda,i)$ is evaluated for the line graph $(L_n,\sim,{\bf 1}_n)$ at vertex $i \in \llbracket 1,n-1 \rrbracket$.

{\bf (1).}
 We have that $h(\lambda,i)$ equals $H(\lambda,n-i,n)$ from Definition~\ref{d.hkell}, where the reversed position $n-i$ is considered because, in the proposition, the root of $L_n$ is at $n$, whereas this root is at zero in Lemma~\ref{l.hfact}.
 The sought formula follows from this lemma and by a computation of~(\ref{e.stake}).
 
 {\bf (2).} This follows by computing the right-hand side of the formula in the first part. A game-theoretic interpretation is available, given Theorem~\ref{t.nashform}. If Maxine and Mina switch roles, the root-reward tree $L_n$ is reflected about $n/2$, and the initial fortune $\lambda$ is replaced by $\lambda^{-1}$, then the payment made or received by either player is unchanged in law. Thus the sought formula reduces to that in the preceding part.
  
  {\bf  (3).} The function $(0,\infty) \lora [0,1]: \lambda \mapsto h(\lambda,i)$ is seen to be continuously differentiable  in view of Lemma~\ref{l.hfact};
  when $\lambda =1$, $\tfrac{\partial}{\partial\lambda}h(\lambda,i)$ equals $\tfrac{i(n-i)}{2n}$, as we find by computing $\lim_{\lambda \searrow 1} \tfrac{\partial}{\partial\lambda} h(\lambda,i)$.
  The claimed stake formula thus results by computing the limit as $\lambda \searrow 1$ of the formula obtained in the proposition's first part.
  
  {\bf (4).} In computing~(\ref{e.stake}), a factor of $p(n) - p(0)$ appears in numerator and denominator, so that the formula coincides with the special case already treated. Alternatively, a game-theoretic interpretation, which is contingent on Theorem~\ref{t.nashform}: if we deduct the constant $p(0)$ from $\pay$, and then revalue currency so that $\pay \to \pay/(p(n) - p(0))$, we do not affect strategy. \qed
  
  {\bf Proof of Theorem~\ref{t.stakeformula}: (1).} When $\lambda \in (0,\infty)$ is not equal to one, we compute terms in the expression $\stake(1,\lambda,v)$ as specified in~(\ref{e.stake}) 
by means of the formula for $h(\lambda,v)$ in Proposition~\ref{p.rootrewardinfinitybias}. Whatever the value of $\lambda \in (0,\infty)$,
the numerator in the formula satisfies
\begin{equation}\label{e.deltaformula}
 \Delta(\lambda,v)  \, = \,   \big( H(\lambda,s_k,d_k - 1) -H(\lambda,s_k,d_k + 1) \big)  \prod_{i=0}^{k-1} H(\lambda,s_i,d_i) \, ;
 \end{equation}
 when $\lambda \neq 1$, it thus equals
$$
   \tfrac{\lambda^{-(s_k - d_k)}}{1 - \lambda^{-s_k}}\big( \lambda - \lambda^{-1}\big)
  \prod_{i=0}^{k-1} H(\lambda,s_i,d_i)  \, = \,  \tfrac{\lambda^{-(s_k - d_k)}}{1 - \lambda^{-(s_k-d_k)}}\big( \lambda - \lambda^{-1}\big) h(\lambda,v) \, .
$$
That is,  $\Delta(\lambda,v)  = \Psi(\lambda,s_k - d_k) (\lambda - \lambda^{-1}) h(\lambda,v)$.
To compute the denominator of~(\ref{e.stake}), note that 
$$
 \tfrac{\partial}{\partial \lambda} h(\lambda,v) \, = \, h(\lambda,v) \sum_{i=0}^k  \tfrac{\partial}{\partial \lambda}  (\log H)(\lambda,s_i,d_i)
 \, = \, h(\lambda,v) \lambda^{-1} \sum_{i=0}^k  \Big( (s_i - d_i) \Psi(\lambda,s_i - d_i) - s_i  \Psi(\lambda,s_i) \Big) \, .
$$
Dividing the obtained expression for $\Delta(\lambda,v)$ by the product of $(\lambda+1)^2$ and the last display, we obtain the formula for $\stake(1,\lambda,v)$ asserted by Theorem~\ref{t.stakeformula}(1).

{\bf (2).} Now $\lambda$ equals one. By~(\ref{e.deltaformula}) with $\lambda = 1$, we find that $\Delta(\lambda,v) = 2 s_k^{-1}  \prod_{i=0}^{k-1} H(1,s_i,d_i)$. Since $H(1,s_k,d_k)$ equals $d_k/s_k$,
we obtain $\Delta(\lambda,v) = 2 (s_k - d_k)^{-1} h(\lambda,v)$.  

We showed in the proof of Proposition~\ref{p.stake}(3) that $\tfrac{\partial}{\partial \lambda} H(\lambda,s,d) \big\vert_{\lambda =1} = \tfrac{d(s-d)}{2s}$ for $s \in \nwithoutzero$ and $d \in \llbracket 0,s\rrbracket$. 
Since $H(1,s,d) = (s-d)/s$, $\tfrac{\partial}{\partial \lambda}  (\log H) (\lambda,s,\ell) \big\vert_{\lambda =1} =  d/2$. Thus,
$$
\tfrac{\partial}{\partial \lambda} h(\lambda,v)  \big\vert_{\lambda =1} \, = \, h(\lambda,v) \sum_{i=0}^k  \tfrac{\partial}{\partial \lambda}  (\log H)(\lambda,s_i,d_i)\big\vert_{\lambda =1} \, = \,  2^{-1}h(1,v)   \sum_{i=0}^k d_i \, .
$$
From these inputs, we obtain Theorem~\ref{t.stakeformula}(2). \qed

We may now derive the second consequence of Theorem~\ref{t.stakeformula} stated in the introduction.

{\bf Proof of Proposition~\ref{p.threelambda}.}
We may write the quantities $\dout(v)$ and $\din(v)$ from Definition~\ref{d.inout}
in terms of the journey data  $\big\{ (s_i,d_i): i \in \llbracket 0, k \rrbracket \big\}$  of $v \in V$. Indeed, it is readily verified that $\dout(v) = s_k - d_k$
and $\din(v) = \sum_{i=0}^k d_i$. For $l \in \nwithoutzero$, $\lim_{\lambda \searrow 0} \Psi(\lambda,\ell) = -1$.
 The first two assertions in Proposition~\ref{p.threelambda} arise from Theorem~\ref{t.stakeformula}(1) in view of these facts. The third assertion is implied by the next corollary. \qed
 
The next result treats the high $\lambda$ asymptotic for the stake function. 

\begin{corollary}\label{c.highlambda}
Let $(V,E,{\bf 1}_r)$ be a root-reward tree, and  let $v \in V$  have journey data $\big\{ (s_i,d_i): i \in \llbracket 0, k \rrbracket \big\}$. 
Let $J$ denote the set of $j \in \llbracket 0, k \rrbracket$ such that $s_j - d_j$ is minimal. 
As $\lambda \to \infty$,
$$
\stake (1,\lambda,v) \, = \,  \frac{1}{\vert J \vert (s_i - d_i)} \cdot \lambda^{s_i - d_i - (s_k - d_k)}\big( 1 + O(\lambda^{-1}) \big) \, ,
$$
where $i$ is any element of $J$.
\end{corollary}

{\bf Proof.} This is due to Theorem~\ref{t.stakeformula}(1); $d_i > 0$ for each $i \in \llbracket 0,k \rrbracket$; and $\lim_{\lambda \nearrow \infty}\Psi(\lambda,\ell) \lambda^\ell = 1$ for $\ell \in \nwithoutzero$.

A further consequence of Theorem~\ref{t.stakeformula} will be needed.
\begin{corollary}\label{c.smallstakes}
For $L \in (0,\infty)$, there exists $h > 0$ such that $\stake(\e,\lambda,v) \geq h \e$ for $\e \in (0,1]$, $\lambda \in (0,L)$ and $v \in \openmac$.
\end{corollary} 
{\bf Proof.} Since $\stake(\e,x,v) = \e \, \stake(1,x,v)$, it suffices to consider $\e =1$. Note that the function $(0,\infty) \lora [0,\infty): \lambda \mapsto \stake(1,\lambda,v)$ is continuous: this follows directly from Theorem~\ref{t.stakeformula}(1) when $\lambda \neq 1$, and from a short computation involving both parts of this theorem in the remaining case. The function may be extended to a continuous function on $[0,\infty)$ by the first assertion of Proposition~\ref{p.threelambda}. Since the function is strictly positive, it is bounded away from zero uniformly on $[0,L]$ for any $L > 0$. \qed

 \section{Games of finite horizon}\label{s.finitehorizon} 

In Section~\ref{s.basics} we specify finite-horizon versions of tug-of-war, with and without stakes. In Section~\ref{s.totvar}, we use the stakeless versions  to prove Proposition~\ref{p.totvar} on the equality of our stake function formulas. Then, in Section~\ref{s.finitehorizonsolution}, we solve the finite-horizon games with stakes.

\subsection{The basics}\label{s.basics}  
  
Let $n \in \nwithoutzero$. In $\game_n(\e,\lambda,v)$, the rules of $\game(\e,\lambda,v)$ are followed, but the game is forcibly ended at the end of the $n$\textsuperscript{th} turn if the counter has then yet to arrive in $\boundarymac$ (so that $X_n \in \openmac$). If the game ends for this reason, the value of $\pay$ equals $h(\lambda,X_n)$. A given payment---say zero or one---would seem to be a more basic choice, but our choice will 
in essence 
permit us to view   $\game_n(\e,\lambda,v)$ as a truncation of  $\game(\e,\lambda,v)$ to the first $n$ turns without the truncation introducing any mean error.

To modify the definition of the mixed and pure strategy spaces $\mc{S}_\pm$ and $\mc{P}_\pm$ in treating the finite game $\gameneps$ for $n \in \nwithoutzero$,
we merely replace the turn index set~$\nwithoutzero$ by $\intint{n}$ in Definition~\ref{d.strategy}, calling the new spaces $\mc{S}_\pm(n)$ and $\mc{P}_\pm(n)$. With this change made, the value $\val_n(\e,\lambda,v)$ of $\game_n(\e,\lambda,v)$ remains specified by Definition~\ref{d.purevalue}.


 \begin{definition}\label{d.hn}
For a boundary-payment graph $(V,E,p)$, parameters $\e \in (0,1]$, $\lambda \in [0,\infty)$, and $n \in \nwithzero$, let $h_n(\e,\cdot,\cdot):(0,\infty) \times V \lora (0,\infty)$
be iteratively specified by $h_0(\e,\lambda,v) = h(\lambda,v)$ for $v \in V$
and, for $n \in \nwithoutzero$, 
$$
 h_n(\e,\lambda,v) :=
\begin{cases}
\, p(v), \textrm{ for } v \in \boundarymac, \\
\, \tfrac{\e\lambda}{\lambda +1} \displaystyle\max_{u \sim v} h_{n-1}(\e,\lambda,u) +  \tfrac{\e}{\lambda +1} \displaystyle\min_{u \sim v} h_{n-1}(\e,\lambda,u) + (1-\e) h_{n-1}(\e,\lambda,v), \textrm{ for } v \in \openmac.
\end{cases}
$$
\end{definition}

\begin{lemma}\label{l.hnconv}
Let $(V,E,p)$ be a boundary-payment graph. For $\e \in (0,1]$, $\lambda \in [0,\infty)$ and $v \in V$, 
 $h_n(\e,\lambda,v)$ equals $h(\lambda,v)$.
\end{lemma}

{\bf Proof.} An induction on $n \in \nwithzero$. The base case $n=0$ holds by definition. Let $n \in \nwithoutzero$. The inductive hypothesis at index $n - 1$ and Definition~\ref{d.hn} imply that
$$
 h_n(\e,\lambda,v) - h(\lambda,v) =
\begin{cases}
\, \, 0 \, \, \textrm{for} \, \, v \in \boundarymac  \, , \\
\, \,  \tfrac{\e\lambda}{\lambda +1} \displaystyle\max_{u \sim v} h(\lambda,u) +  \tfrac{\e}{\lambda +1} \displaystyle\min_{u \sim v} h(\lambda,u)  -\e h (\lambda,v) \, \, \textrm{for} \, \, v \in \openmac \,   .
\end{cases}
$$
Thus, the inductive hypothesis at index $n$ holds by~(\ref{e.h}). \qed

 In root-reward trees, $\mc{V}_+(\lambda,v)$ from Definition~\ref{d.vn} is a singleton, the parent $v_+$, for any $v \in V$. This property is transmitted to the specification of lazy biased walks that we saw in Definition~\ref{d.theta}. We are led to a redefinition of these walks in the more general context of boundary-payment graphs (which, briefly at least, we now consider), in which $\mc{V}_+(\lambda,v)$ may have several elements; we also work now with a finite horizon~$n$.
\begin{definition}\label{d.thetafinitehorizon}
Let $\Theta_n(\lambda)$ denote the set of mappings $(\theta_-,\theta_+): \openmac \times \llbracket {0,n-1 \rrbracket} \lora V \times V$ such that $\theta_-(v,i) \in \mc{V}_-(\lambda,v)$ and  $\theta_+(v,i) \in \mc{V}_+(\lambda,v)$ for each $v \in \openmac$. For $\theta = (\theta_-,\theta_+) \in \Theta_n(\lambda)$, 
let $X_\theta: \llbracket 0,F \rrbracket \lora V$
 denote the Markov process such that  $X_\theta(0) =v \in \openmac$  and, for $k \in \nwithoutzero$,
$$
 X_\theta(k+1) =
\begin{cases}
\,  X_\theta(k) \, \, \textrm{with probability $1-\e$}  \, , \\
\,  \theta_+\big(X_\theta(k),k\big) \, \, \textrm{with probability $\e \tfrac{\lambda}{1 + \lambda}$}  \, , \\
\, \theta_-\big(X_\theta(k),k\big) \, \, \textrm{with probability $\e \tfrac{1}{1 + \lambda}$}  \,   .
\end{cases}
$$
The process is stopped either at time $n$ or on arrival at $\boundarymac$, so that $F$ is the minimum of $n$ and  $\min \big\{ j \geq 0: X_\theta(j) \in \boundarymac \big\}$.
\end{definition}

\begin{definition}\label{d.leisurelytugofwar}
Let $(V,E,p)$ be a boundary-payment graph. Let $\e \in (0,1]$, $q \in [0,1]$ and $v \in \openmac$.
The leisurely version~$\tugofwar(\e,q,v)$ of tug-of-war may be specified by adapting Definition~\ref{d.tugofwar}. At any given turn, a move takes place with probability $\e$, the decision taken independently of other randomness; if a move does take place, the existing rules apply.
\end{definition} 

\begin{lemma}\label{l.tugnash}
Let $(V,E,p)$ be a boundary-payment graph.  A pure Nash equilibrium exists for the game $\tugofwar_n\big(\e,\tfrac{\lambda}{1 + \lambda},v\big)$, the finite horizon version of the one in Definition~\ref{d.leisurelytugofwar}. Under any pure Nash equilibrium, the gameplay process has the law of $X_\theta: \llbracket 0,F \rrbracket \lora V$ for some    
$\theta \in \Theta_n(\lambda)$. The mean payment at any such equilibrium equals  $h(\lambda,v)$.
\end{lemma}
{\bf Proof.} We prove this by induction on $n$. If Maxine and Mina propose respective random moves $v_+$ and $v_-$ at the first turn in  $\tugofwar_n\big(\e,\tfrac{\lambda}{1+\lambda},v\big)$, and then adhere to strategies in a Nash equilibrium, then, by the inductive hypothesis for index~$n-1$, the mean payment in the resulting gameplay will equal the mean of
$$
 \e \tfrac{\lambda}{1 + \lambda} h(\lambda,v_+) + \e  \tfrac{1}{1 + \lambda} h(\lambda,v_-) + (1-\e)h(\lambda,v) \, ,
$$
where the mean is taken over the randomness in their move choices.
We see then that it is necessary and sufficient for a pure strategy pair to be a Nash equilibrium that $v_+ \in \mc{V}_+(\lambda,v)$ and $v_- \in \mc{V}_-(\lambda,v)$ hold alongside the players adhering jointly to play at a pure Nash equilibrium in the subgame copy of $\tugofwar_{n-1}(\e)$ that takes place after the first turn. By taking $v_+$ and $v_-$ to be deterministic elements of the respective sets, we obtain the inductive hypothesis's claim about the form of Nash equilibria; in view of~(\ref{e.h}), the mean payment formula then results from the displayed equation for the above such choices of $v_+$ and~$v_-$. \qed

\subsection{The lambda-derivative of game value: proving Proposition~\ref{p.totvar}}\label{s.totvar}

Here we use probabilistic reasoning for constant-bias tug-of-war inspired by the perturbative argument presented heuristically for the stake-governed version in Section~\ref{s.stakeperturb}
in order to prove that  $\tfrac{\partial}{\partial  \lambda}h(\lambda,v)$ equals $(\lambda +1)^{-2} \e \, \E \,  
\totvar \big( \e, \lambda,\theta,v \big)$. 
In order to prove Proposition~\ref{p.totvar}, of which the above statement forms part, we will first argue that  $\stake(\e,\lambda,v)$ from~(\ref{e.stake}) has an alternative, $n$-dependent, form~$\stake_n(\e,\lambda,v)$. We thus have cause to introduce some notation that relates the infinite-horizon game with its finite-horizon counterpart. The usage of the notation $\cdot \rfloor n$ will be characteristic when  truncation of an aspect $\cdot$ of infinite gameplay to the first $n$ turns is concerned. 

\begin{definition}\label{d.thetan}
Let $\theta$ be an element of the set $\Theta$ specified in Definition~\ref{d.theta}. For $\lambda \in (0,\infty)$, $v \in \openmac$ and $n \in \nwithoutzero$, let $X_{\theta \rfloor n}(\lambda,\cdot): \llbracket 0, F_{\theta \rfloor n}\rrbracket \lora V$, $X_{\theta \rfloor n}(0) =v$,
be given by $X_{\theta \rfloor n}(\lambda,k) = X_\theta(\lambda,k)$ for $k \in \llbracket 0, F_{\theta \rfloor n} \rrbracket$. Here the truncated finish time $F_{\theta \rfloor n}$ equals $F_\theta \wedge n$. 
If the counter is in open play at time $n$, so that $X_\theta(n) \in \openmac$, we formally set $F_{\theta \rfloor n}$ equal to $n+1$.

Extending the shorthand in Definition~\ref{d.theta}, we will write $X_{\theta \rfloor n}(\cdot)  = X_{\theta \rfloor n}(\lambda,\cdot)$.
\end{definition}

To specify the quantity $\stake_n(\e,\lambda,v)$, let $\theta \in \Theta$.
For $\lambda \in (0,\infty)$, set
$$
\totvar \big( \e, \lambda,\theta \rfloor n,v \big) \, = \, \sum_{i=0}^{F_{\theta \rfloor n} -1} \Delta \big(\lambda,X_{\theta \rfloor n}(i)\big) \, ,
$$ 
 where recall that $X_{\theta \rfloor n}(0) = v$, and  $\Delta(\lambda,u) = h \big( \lambda, u_+ \big) - h \big( \lambda, u_- \big)$ with $u_-$ an element of the set $\mc{V}_-(u)$ from Definition~\ref{d.inout}.
 A device in fact permits simpler notation. We specify $X_\theta(i)$ to equal its terminal value $X_\theta(F_\theta)$ whenever $i>F_\theta$; thus, $X_\theta$ is now a random map from $\nwithzero$ to $V$. And we set $\Delta(\lambda,v) = 0$ for $v \in \boundarymac$, so that $\Delta(\lambda,\cdot)$ now maps $V$ to $[0,1]$. 
Under these conventions, 
\begin{equation}\label{e.totvaralt}
\totvar \big( \e,\lambda,\theta \rfloor n,v \big) \, = \, \sum_{i=0}^{n -1} \Delta \big(\lambda,X_\theta(i)\big) \, .
\end{equation}

 This quantity is a finite-horizon counterpart to 
$\totvar \big( \e,\lambda,\theta ,v \big)$, which is defined before Proposition~\ref{p.totvar}. 
We set 
\begin{equation}\label{e.stakeformula}
  \stake_n(\e,\lambda,v) \, = \, \frac{\Delta(\lambda,v)}{\totvar \big( \e, \lambda,\theta \rfloor n,v \big)  
   +  \e^{-1} (\lambda +1)^2 \sum_{w \in \openmac}
  \tfrac{\partial}{\partial  \lambda}h(\lambda,v) \cdot \mu_n(v,w)} 
   \, ,\end{equation}
where  
 $\mu_n(v,w) = \PP \big( X_\theta(n) = w 
  \big)$, so that $\big\{ \mu_n(v,w): w \in \openmac \big\}$ is the sub-probability measure of the counter's location in the event that the game is unfinished at time~$n$.

(It is perhaps worth noting that, when $n$ is high,   $\stake_n(\e,\lambda,v)$ is equal to the product of~$\e$ and a unit-order quantity. Indeed, $\totvar \big( \e,\lambda,\theta \rfloor n,v \big)$  typically reports $\Theta(\e^{-1})$ terms for each jump that $X_\theta$ makes, because this process is a $(1-\e)$-lazy walk. The implied unit-order quantity is apparent in~(\ref{e.stake}).)

The next result shows that the new stake formula~(\ref{e.stakeformula}) equals the original one~(\ref{e.stake}). Recall from Proposition~\ref{p.differentiable} that 
 the function $(0,\infty) \lora (0,1): \lambda \mapsto h(\lambda,v)$ is differentiable for $v \in \openmac$.
\begin{proposition}\label{p.nashformprep}
\leavevmode
\begin{enumerate} 
\item For $n \in \nwithoutzero$, $$
  \tfrac{\partial}{\partial  \lambda}h(\lambda,v) \, = \, \frac{\e}{(\lambda + 1)^2} 
\cdot \E \, 
\totvar \big(\e, \lambda,\theta \rfloor n,v \big)
   \, + \, \sum_{w \in \openmac}
  \tfrac{\partial}{\partial  \lambda}h(\lambda,w) \cdot \mu_n(v,w) 
   \, .
$$
\item Let $\theta \in \Theta$.
 The expression  $\E \sum_{i=0}^\infty \Delta\big(\lambda,X_\theta(\e,\lambda,i)\big)$ is finite. For each $n \in \nwithoutzero$, it equals
\begin{equation}\label{e.stability}
  \E \,
\totvar \big( \e,\lambda,\theta \rfloor n,v \big) 
   \, + \, \e^{-1}(\lambda + 1)^2 \sum_{w \in \openmac}
  \tfrac{\partial}{\partial  \lambda}h(\lambda,w) \cdot \mu_n(v,w) \, .
\end{equation}
\item This expression also equals $\e^{-1} (\lambda + 1)^2  \tfrac{\partial}{\partial  \lambda}h(\lambda,v)$.

\item The stake function~(\ref{e.stakeformula}) is independent of $n \in \nwithoutzero$ and equals $\stake(\e,\lambda,v)$ as specified in~(\ref{e.stake}).
\end{enumerate}
\end{proposition}
{\bf Proof of Proposition~\ref{p.totvar}.}
Recall that  $\E \, \totvar(\e,\lambda,v) = \E \sum_{i=0}^\infty \Delta\big(\lambda,X_\theta(i) \big)$.
By Proposition~\ref{p.nashformprep}(3), this expression equals  $\e^{-1}(\lambda + 1)^2\tfrac{\partial}{\partial \lambda} h(\lambda,v)$.  Naturally then the alternative stake  formula~(\ref{e.altstake}) is obtained from the original formula~(\ref{e.stake}). \qed

{\bf Proof of Proposition~\ref{p.nashformprep}: (1).} 
By Lemma~\ref{l.hnconv}, $h(\lambda,v)$ equals $h_n(\e,\lambda,v)$. The latter quantity is the value of $\tugofwar_n\big(\e,\tfrac{\lambda}{1 + \lambda},v\big)$ by Lemma~\ref{l.tugnash}. 
It is this interpretation that we exploit in this proof. Let $\theta \in \Theta$. By Lemma~\ref{l.tugnash}, 
there exists a Nash equilibrium in  $\tugofwar_n\big(\e,\tfrac{\lambda}{1 + \lambda},v\big)$  such that the counter evolution has the law $X_{\theta \rfloor n}: \llbracket 0 , F_{\theta \rfloor n} \rrbracket \lora V$, $X_{\theta \rfloor n}(0) = v$. (The notational convention leading to~(\ref{e.totvaralt}) will however be in force, so that the domain of this process may formally be taken to be $\llbracket 0, n \rrbracket$.)
 We provide a coupling of these counter evolutions 
 that permits them to be realized simultaneously for differing values of $\lambda$. To this end, let $\big\{ U_i: i \in \llbracket 0, n-1 \rrbracket \big\}$ denote an independent sequence of random variables, each having the uniform law on $[0,1]$.  Let $i \in \intint{n}$.
 If a move takes place at the $i$\textsuperscript{th} turn in the coupled copy of 
  $\tugofwar_n\big(\e,\tfrac{\lambda}{1 + \lambda},v\big)$,
  then the condition that ensures that it is Maxine who wins the right to move is that $U_{i-1} \leq \tfrac{\lambda}{1 + \lambda}$.

Let $\phi > 0$. In order to study the difference   $h_n(\e,\lambda+ \phi,v) -  h_n(\e,\lambda,v)$, we call the coupled copy of the $\lambda$-biased game the {\em original} game; and the $(\lambda + \phi)$-biased game, the {\em alternative game}. 
In addition, we specify the {\em intermediate} game. To do so, let the {\em disaccord} set $D \subseteq \llbracket 1, n \rrbracket$ denote the set of indices of turns at which a move occurs that Maxine wins in the alternative game, but not in the original one.
Then let 
 $\sigma \in \intint{n}$ denote the minimum of $D$.  That is, $\sigma$ is the smallest $i \in \llbracket 1 ,n \rrbracket$ such that $U_i \in \big(\tfrac{\lambda}{\lambda + 1},\tfrac{\lambda + \phi}{\lambda  + \phi +1} \big]$. (Which is to say, in essence: $\sigma$ is the first turn that Maxine wins in the alternative game but not in the original one. Note however that $\sigma$ is defined even if the game has already finished by the time that its value is reached. This device will permit a convenient independence property.)
  If $D$ is empty, set $\sigma = n + 1$. Gameplay in the intermediate game follows that of the alternative game until move~$\sigma$ when, at least if $\sigma \leq n$, this gameplay diverges from that in the original game; after this, the gameplay coincides with that in the original game. In other words, the cutoff $\chi$ in the condition $U_{i-1} \leq \chi$ 
takes the value $\chi = \tfrac{\lambda + \phi}{\lambda + \phi + 1}$ until move $i=\sigma$ and then changes to $\chi = \tfrac{\lambda}{\lambda + 1}$ for higher-indexed moves. 

With the indices $0$, $1$ and $2$ denoting the original, intermediate and altered games, write $P_i$, $0 \leq i \leq 2$, for the mean value of $\pay$ under the corresponding game. As we have noted, $h(\lambda,v) = h_n(\lambda,v)$ is the value of $\tugofwar_n \big(\e,\tfrac{\lambda}{1+\lambda},v \big)$; thus,
\begin{equation}\label{e.pzeropone}
 P_0 = h(\lambda,v)  \, \, \, \, \textrm{and} \, \, \, \, P_2 = h(\lambda+\phi,v) \, .
\end{equation}
Let $j \in \intint{n}$. The occurrence of $\sigma = j$ imparts no information regarding the trajectory $X_\theta$ for any value of $\lambda$, because $\sigma$ takes the form $G \wedge (n+1)$ where $G \geq 1$
is a geometric random variable of parameter $\tfrac{\lambda + \phi}{\lambda  + \phi +1} - \tfrac{\lambda}{\lambda +1}$ that is independent of other randomness. We see then that
$$
 \E \Big[ \big(  P_1 - P_0 \big) {\bf 1}_{\sigma =j} \Big]  \, = \, \E  \Big[ \big(  P_1 - P_0 \big) \, \Big\vert \, \sigma =j \Big]  \,  \PP \big( \sigma =j \big) 
 \, = \, \E  \Delta \big( \lambda, X_\theta(j-1) \big) \,  \PP \big( \sigma =j \big) \, ,
$$
where the just noted independence was used in the latter equality to write 
$\E  \big[ (  P_1 - P_0 ) \, \big\vert \, \sigma =j \big]$ as the mean  
$\E \Delta \big( \lambda , X_\theta(j-1) \big)$
of the difference in terminal payment of
 $X_\theta$ with parameter $\lambda$ according to whether Maxine or Mina wins at turn~$j$, the mean taken over the location $X_\theta(j-1)$ as specified in the original dynamics (and thus with parameter $\lambda$).
Recalling $F_{\theta \rfloor n}$ from Definition~\ref{d.thetan}, note then that
$$
 P_1 - P_0 \, = \, \E  \sum_{i=0}^{n-1} \Delta \big( \lambda, X_\theta(i) \big) \, \zeta_i \, \, + \, \, \PP \big( \sigma = n+1 \big) \E \, \bigg[ \Big( h\big(\lambda + \phi,X_\theta(n)\big) - h\big(\lambda,X_\theta(n) \big) \Big) {\bf 1}_{F_{\theta \rfloor n} = n+1} \bigg] \, ,
$$
where $\zeta_i = \PP (\sigma = i+1)$. 
We {\em claim} that
$$
 \zeta_i \, = \, \e \tfrac{\phi}{(\lambda + 1)^2} + i \e^2 \phi^2 \Theta(1) \, .
$$
To verify this, note that $\zeta_i = (1 - \alpha)^i \alpha$, where $\alpha > 0$ is the probability that, at a given turn, a move takes place that Maxine wins in the altered game but not in the original one.  
This quantity satisfies
$$
\alpha \, = \, \e \cdot \Big(\tfrac{\lambda + \phi}{\lambda + \phi +1} - \tfrac{\lambda}{\lambda + 1} \Big) \,  = \, \e\tfrac{\phi}{(\lambda + 1)^2} +  \e\tfrac{\lambda}{\lambda + 1}  \phi^2 \Theta(1) \, .
$$
From this, the claim readily follows. The claim yields the consequence that
$$
 \PP \big( \sigma = n+1 \big) \, = \, 1 - n  \e \tfrac{\phi}{(\lambda + 1)^2} +  n^2 \e^2 \phi^2 \Theta(1) \, .
$$
Applying the claim and its consequence leads to the formula
\begin{eqnarray*}
 P_1 - P_0 & = & \e \, \Big( \tfrac{\phi}{(\lambda + 1)^2} +   n \e  \phi^2 O(1) \Big) \cdot \E  \sum_{i=0}^{n-1} \Delta \big(\lambda, X_\theta(i) \big) \\
 & & \qquad    + \, \, \Big( 1 - n  \e \tfrac{\phi}{(\lambda + 1)^2} +  n^2 \e^2 \phi^2 \Theta(1) \Big) \sum_{w \in \openmac}
   \big( h(\lambda + \phi,w) - h(\lambda,w) \big) \cdot \mu_n(v,w) 
 \, .
\end{eqnarray*}
Note also that 
$$
 P_2 - P_1  \, \leq \, \PP \big( \vert D \vert \geq 2 \big) \, ,
$$
whose right-hand side is at most $n^2 \alpha^2 = n^2 \e^2 O(\phi^2)$. 

Since $P_2 - P_0$ equals $h(\lambda + \phi,v) - h(\lambda,v)$ by~(\ref{e.pzeropone}), 
 the two preceding displays yield Proposition~\ref{p.nashformprep}(1) when a limit $\phi \searrow 0$ is taken.

{\bf 
(2).} Proposition~\ref{p.nashformprep}(1) implies that the expressions~(\ref{e.stability}) are equal to $\e^{-1}(\lambda + 1)^2   \tfrac{\partial}{\partial  \lambda}h(\lambda,v)$, whatever the value of $n \in \nwithoutzero$. This common value is finite by Proposition~\ref{p.differentiable}. Naturally, the sequence $\E \, \sum_{i=0}^{n-1} \Delta \big(\lambda, X_\theta(i) \big)$ converges pointwise in $\lambda \in (0,\infty)$ as $n \to \infty$ to   $\E \, \sum_{i=0}^\infty \Delta \big(\lambda, X_\theta(i) \big)$. Since the $n$-indexed sequence  $\e^{-1}(\lambda + 1)^2 \sum_{w \in \openmac}
  \tfrac{\partial}{\partial  \lambda}h(\lambda,w) \cdot \mu_n(v,w)$ is non-negative, we find then that   $\E \, \sum_{i=0}^\infty \Delta \big(\lambda, X_\theta(i) \big)$ is finite. Note that this expression at least 
  $q(\lambda) \E F_\theta$, where $q(\lambda) = \min_{u \in \openmac}  \Delta (\lambda, u)$ is positive by Lemma~\ref{l.a}(4). Thus, $\E \, F_\theta$ is finite, so that  
  \begin{equation}\label{e.ftail}
  \sum_{w \in \openmac} \mu_n(v,w) = \PP (F_\theta > n) \, \lora \, 0 \, \, \, \, \textrm{as} \, \, \, \,  n \to \infty \, ,
    \end{equation} 
    and
$$
   \sum_{w \in \openmac}
  \tfrac{\partial}{\partial  \lambda}h(\lambda,w) \cdot \mu_n(v,w) \,\leq \, \sup_{w \in \openmac}
  \tfrac{\partial}{\partial  \lambda}h(\lambda,w) \cdot \sum_{w \in \openmac} \mu_n(v,w) \, \to \, 0 \, \, \, \textrm{as $n \to \infty$} \, .
  $$
  This implies that the common value of the expressions~(\ref{e.stability}) is   $\E \, \sum_{i=0}^\infty \Delta \big(\lambda, X_\theta(i) \big)$. 
  
  {\bf (3).} This follows from Proposition~\ref{p.nashformprep}(1,2).
  
 
 {\bf (4).} This is due to Lemma~\ref{l.hnconv} and the proposition's first part.  \qed
 
 An alternative characterization of the stake function will be useful. 
\begin{lemma}\label{l.bsolution}
 The unique solution $b$ of the equation
\begin{equation}\label{e.oneminusb}
 (1 -b)   \tfrac{\partial}{\partial  \lambda}h(\lambda,v) =   \tfrac{\lambda}{\lambda + 1}  \tfrac{\partial}{\partial  \lambda}h(\lambda,v_+) +  \tfrac{1}{\lambda + 1}  \tfrac{\partial}{\partial  \lambda}h(\lambda,v_-)    
\end{equation}
is $b = \stake(1,\lambda,v)$.
\end{lemma}
{\bf Proof.} For $v \in V$ and $i \in \nwithzero$, set $W(\lambda,v,i) = \E \, \sum_{j=i}^\infty \Delta(\lambda,X_\theta(1,\lambda,i))$, where $X_\theta(1,\lambda,0) = v$. (Note that $W(\lambda,v,0)$ equals $\E \, \totvar(1,\lambda,v)$.)
By Proposition~\ref{p.nashformprep}(3), $W(\lambda,v,0) = (\lambda + 1)^2  \tfrac{\partial}{\partial  \lambda}h(\lambda,v)$. By multiplying~(\ref{e.oneminusb}) by  $(\lambda + 1)^2$, we thus find that
$$
 (1-b) W(\lambda,v,0) = W(\lambda,v,1)
$$
This implies that $b = \tfrac{\Delta(\lambda,v)}{W(\lambda,v,0)}$. Thus, $b$ equals $\stake(1,\lambda,v)$ as it is specified in the alternative formula~(\ref{e.altstake}) in Proposition~\ref{p.totvar}.
This completes the proof. \qed

  
\subsection{Solving the finite-horizon games}\label{s.finitehorizonsolution}

We now derive a counterpart to our main results for games of finite horizon, for which the fundamental technique of backward induction is available. The principal result, Theorem~\ref{t.leisurely.theorem}, indicates that the global saddle hope offered in Section~\ref{s.picture} is realized for finite horizon leisurely games.
This theorem asserts the existence of a non-random stake proportion that the players share when they adhere to any Nash equilibrium. The result moreover offers a formula for this shared proportion. This formula, which is simply $\stake(\e,\lambda,v)$ from~(\ref{e.stake}), does not depend on $n$, as  the choice of horizon-time reward of $h(\lambda,v)$ in the specification of $h_n(\e,\lambda,v)$ in Definition~\ref{d.hn} serves to ensure.

Recall that $\game_n(\epsilon,\lambda,v)$ denotes the $n$-finite horizon version of $\game(\epsilon,\lambda,v)$.

  \begin{theorem}\label{t.leisurely.theorem} 
 Let $T = (V,E,{\bf 1}_r)$ be a root-reward tree. There exists $\epsilon_0 \in (0,1)$ such that, when $\epsilon \in (0,\e_0)$, $\lambda >0$ and $n \in \nwithoutzero$, the following hold for each $v \in \openmac$. 
\begin{enumerate}
\item The value of ${\rm Game}_n(\e,\lambda,v)$ exists and equals $h(\lambda,v)$.
\item
The game $\game_n(\epsilon,\lambda,v)$ has a Nash equilibrium.  Any Nash equilibrium is a pair of conforming strategies as specified by Definition~\ref{d.conforming}.  
Conversely, any such pair is a strict Nash equilibrium of $\game_n(\epsilon,\lambda,v)$. Indeed, 
let $(\Pminus,\Pplus)$ be a pair of conforming strategies. Then
\begin{equation}\label{strict Nash equilibrium for conforming}
    \begin{aligned}
        M_n(\lambda,v,\Pminus,S_+) < h(\lambda,v), \quad \text{if $S_+$ is not conforming against $\Pminus$ }, \\
        M_n(\lambda,v,S_-, \Pplus) > h(\lambda,v), \quad \text{if $S_-$ is not conforming against $\Pplus$ },
    \end{aligned}
\end{equation}
where $M_n(\lambda,v,S_-,S_+)$ denotes the mean payment of $\game_n(\epsilon,\lambda,v)$ arising from play under the strategy pair $(S_-,S_+)$.
\end{enumerate}
\end{theorem}

In Section~\ref{s.localglobal}, we considered the value     $\operatorname{Val}(\epsilon,\lambda,v,a,b)$ of a game constrained so that Maxine stakes~$a$, and Mina~$b$, at the first turn.
To prove Theorem~\ref{t.leisurely.theorem}, we define
 \begin{equation}\label{constrained value function}
    \operatorname{V}(\epsilon,\lambda,v,a,b) \, = \, (1-\epsilon)h(\lambda_1,v) \, + \, \epsilon \, \bigg(\frac{a}{a+b}h(\lambda_1,v_+)+\frac{b}{a+b} h(\lambda_1,v_-) \bigg) \, , 
\end{equation}
where $v \in \openmac$, $\lambda >0$, $\lambda_1 = \frac{\lambda-a}{1-b}$, and $v_+$ and~$v_-$ are neighbours of $v$ which maximize and minimize $h(\lambda,u)$ (so that $v_-$ is an element of $\mc{V}_-(\lambda,v)$ from Definition~\ref{d.vn}).
The object defined in (\ref{constrained value function})
 coincides with    $\operatorname{Val}(\epsilon,\lambda,v,a,b)$ if we admit that the value of the unconstrained $\game(\e,\lambda,v)$ exists and equals $h(\lambda,v)$.     
In contrast to the heuristic Section~\ref{s.picture}, we now make no conditional statements, so that $\operatorname{V}(\epsilon,\lambda,v,a,b)$ is simply defined by the formula above.

A key tool for proving Theorem~\ref{t.leisurely.theorem} says that this function has a global minimax.


\begin{proposition}\label{p.globalconstrainedminimax}
There exists $\epsilon_0 \in (0,1)$ such that for $\epsilon \in (0,\epsilon_0)$, and $(\lambda,v) \in (0,\infty)\times \openmac$, the function $\operatorname{V}(\epsilon,\lambda,v,a,b)$ has a unique global minimax at $(a,b)=\big(\lambda \stake(\e,\lambda,v),\stake(\e,\lambda,v)\big)$; which is to say, 
\begin{align*}
    \operatorname{V}\big(\epsilon,\lambda,v,a,\epsilon S(\lambda)\big) & \, \leq \, h(\lambda,v) \quad \text{for } 0\leq a \leq \lambda \, , \\
    \operatorname{V}\big(\epsilon,\lambda,v,\epsilon \lambda S(\lambda),b\big) & \,  \geq \, h(\lambda,v) \quad \text{for } 0 \leq b \leq 1 \, . 
\end{align*}
Here $S(\lambda) = \stake(1,\lambda,v)$ from~\eqref{e.stake}. The global minimax is {\em strict}: equality is achieved only when $a = \epsilon\lambda S(\lambda) $ and $b = \epsilon S(\lambda)$.
\end{proposition}
We will prove Proposition~\ref{p.globalconstrainedminimax} via the next result on decay of stake function difference quotients.

Taking $v \in V$ given, we maintain the shorthand $S(x) = \stake(1,x,v)$. Further set $T:(0,\infty) \lora (0,1]$, $T(x) = S(x^{-1})$.
It is worth noting for later use that $S$ and $T$ are at most one, by Proposition~\ref{p.totvar}.

\begin{lemma}\label{l.stakedifferencequotient}
The quantity 
$$ 
 \sup \, \bigg\{ \, (1+x)(1+y) \frac{\vert S(y)-S(x) \vert}{\vert y-x \vert}  : x,y \in (0,\infty) \, , \, x \neq y \, \bigg\}
$$ 
is finite. It remains unchanged if $S$ is replaced by $T$.
\end{lemma}

{\bf Proof.}
First note that the expression $\frac{\vert x - y \vert}{(1+x)(1+y)}$ is invariant under the transformation $(0,\infty)^2 \lora (0,\infty)^2: (x,y) \mapsto (x^{-1},y^{-1})$. Thus the claim made in regard to $T(x) = S(x^{-1})$ reduces to that made for~$S$. 

The assertion made in regard to~$S$ will be proved by establishing a pair of claims.

{\em Claim~{\rm I}.} The function $S:(0,\infty) \lora (0,\infty)$ takes the form
\begin{equation}\label{e.S}
    S(x) = \frac{\sum_{k=0}^m a_kx^k }{ \sum_{k=0}^m b_k x^k } \, , \quad \text{with } b_k>0 \text{ for } 0\leq k \leq m \,,
\end{equation}
for some $m \in \N$. The coefficients $a_k$ are real but may be zero, so that the degree of the numerator is at most that of the denominator. 

{\em Claim~{\rm II}.} Any function $s:(0,\infty) \lora (0,\infty)$ having the above form~\eqref{e.S} satisfies 
$$
\sup \, \Big\{ \, (1+t)^2 s'(t): t \in [0,\infty) \, \Big\} \, < \, \infty \, . 
$$

Indeed, admitting the claims, Lemma~\ref{l.stakedifferencequotient} for the function~$S$ is obtained simply by integrating its derivative and invoking Claim~{\rm II}.

{\em Deriving Claim~{\rm I}.} Some notation is needed.
 Set $P_n(x)   = $$  \sum_{i=0}^n x^i$ for  $k>0$ and $P_0 \equiv 1$, as well as
$$
        Q_{n,i}(x)  \, = \, \sum_{\ell=0}^{i-1}\sum_{j=0}^{n-i-1}x^jP_{n-i-1+\ell-j}(x) \, \, \, \, \textrm{for $n \geq 2$ and $i \in \intint{n-1}$}  \, .
        $$
    Here $Q_{n,i}(x)$ is a polynomial of order $n-2$ each of whose coefficients is positive. 

    
    Note that 
    \begin{align*}
        i-nx^{n-i}+(n-i)x^n &= (1-x)\big(iP_{n-i-1}(x) -(n-i)x^{n-i}P_{i-1}(x)\big) \\
        =(1-x)\sum_{\ell=0}^{i-1}\sum_{k=0}^{n-i-1}\Big(x^k- x^{n-i}x^{\ell} \Big) 
        &= (1-x)^2\sum_{\ell=0}^{i-1}\sum_{k=0}^{n-i-1}x^k P_{n-i-1+\ell-k}(x) = (1-x)^2Q_{n,i}(x) \, .
    \end{align*}
    Recalling from Theorem~\ref{t.stakeformula} the notation  $\Psi(x,\ell ) = \frac{x^{-\ell}}{1-x^{-\ell}}$, note that 
    $$
            (n-i)\Psi(x,n-i) - n\Psi(x,n) = \frac{n-i}{x^{n-i}-1 }-\frac{n}{x^n-1} = \frac{i - n x^{n-i}+(n-i)x^n }{(x^n-1)(x^{n-i}-1) } 
     $$
     equals $\frac{Q_{n,i}(x)}{P_{n-1}(x)P_{n-i-1}(x)}$. 
    Using Theorem~\ref{t.stakeformula}, we find that, for $x\neq 1$, \begin{equation}\label{stake formula as rational function}
        \begin{aligned}
            S(x) &= \frac{(x - 1) \Psi(x,s_k-d_k)}{(x+1)\sum_{i=0}^k \Big(  (s_i - d_i) \Psi(x,s_i - d_i) - s_i  \Psi(x,s_i) \Big)}  = \frac{\frac{1 }{P_{s_k-d_k-1}(x)}}{(x+1)\sum_{i=0}^{k}\frac{Q_{s_i,d_i}(x)}{P_{s_i-1}(x)P_{s_i-d_i-1}(x) } }\\
            &= \frac{ P_{s_k-1}(x) \prod_{\ell =0}^{k-1}P_{s_{\ell} -1}(x)P_{s_{\ell} -d_{\ell} -1}(x) }{(1+x)\sum_{i=0}^{k}\Big( Q_{s_i,d_i}(x) \prod_{\ell\neq i}P_{s_{\ell} -1}(x)P_{s_{\ell} -d_{\ell} -1}(x)\Big)  } \,,
        \end{aligned}
    \end{equation}
where $\big\{ (s_i,d_i): i \in \llbracket 0, k \rrbracket \big\}$ is the journey data of $v$ as specified in Definition~\ref{d.journeydata}. 
    
    With $x=1$, Theorem~\ref{t.stakeformula} tells us that $S(x) = \Big( (s_k - d_k) \sum_{i=0}^k d_i \Big)^{-1}$.  
    
    Since $P_k(1) = k+1$ and  $Q_{n,i}(1) = ni(n-i)/2$, we have $2\frac{Q_{s_i,d_i}(1)}{P_{s_i-1}(1)P_{s_i-d_i-1}(1) } = d_i$.  So  letting $x=1$ in the second inequality in \eqref{stake formula as rational function} yields $\Big( (s_k - d_k) \sum_{i=0}^k d_i \Big)^{-1}$.  Thus, \eqref{stake formula as rational function} is true for any $x\geq 0$.
    In this way, $S(x)$ is a rational function $\frac{f(x)}{g(x)}$. Note that $P_k(x)$ has order $k+1$, while $Q_{n,i}(x)$ has order $n-2$; thus, $f(x)$ is of order 
    $$
    \sum_{\ell=0}^{k}(2s_{\ell} -d_{\ell} -2) \;  - \; (s_k-d_k-1) \, = \,  -2k-1 +  d_k - s_k  \, + \,    \sum_{\ell=0}^{k}(2s_{\ell} -d_{\ell})    \, ,
    $$
     and $g(x)$ is of order \begin{eqnarray*}  
     & &   1 \, + \, \max_{0\leq i \leq k} \, \bigg( \, s_i-2+ \sum_{\ell=0}^{k}(2s_{\ell} -d_{\ell} -2) - (2s_i-d_i-2) \, \bigg) \\
       & = &  -2k-1\, + \, \sum_{\ell=0}^{k} \big(2s_{\ell} -d_{\ell}  \big)   \, + \, \max_{0\leq i \leq k}\big( d_i - s_i   \big) \, .
            \end{eqnarray*}
Hence, the order of $f(x)$ is at most that of $g(x)$; and clearly all coefficients of $g(x)$ are positive.  We can thus write $S(x) = \frac{\sum_{k=0}^m a_kx^k }{ \sum_{k=0}^m b_k x^k }$ with $b_k >0$ for $k \in \llbracket 0, m \rrbracket$ and thereby obtain  Claim~{\rm I}.

{\em Deriving Claim~{\rm II}.}  For $s(x) = f(x)/g(x)$ with $f(x) = \sum_{k=0}^m a_k x^k $ and $g(x) = \sum_{k=0}^m b_k x^k$, we have $s'(x) = \frac{f'(x)g(x)-f(x)g'(x)}{g^2(x)}$, where the function $f'(x)g(x)-f(x)g'(x)$ is a polynomial of order at most $2m-2$ because the  order $2m-1$ monomials cancel. Thus,
    \begin{align*}
        |s'(x)| = \left| \frac{\sum_{k=0}^{2m-2}{A_k}x^k}{\sum_{k=0}^{2m}{B_k}x^k} \right|\leq \frac{C}{(1+x)^2}, \quad \text{for some constant } C>0 \, ,
    \end{align*} 
    where we also used $B_k>0$ for $0\leq k\leq 2m$, a fact direct from $b_k>0$. Moreover, $(1+x)^2\sum_{k=0}^{2m-2}{A_k}x^k$ is a polynomial of order $2m$.

This completes the proof of Lemma~\ref{l.stakedifferencequotient}. \qed

{\bf Proof of Proposition~\ref{p.globalconstrainedminimax}.}   To prove that $\big(\epsilon\lambda S(\lambda), \epsilon S(\lambda)\big)$ is a strict global minimax of $\operatorname{V}(\epsilon,\lambda,v,a,b)$, we need only show that \begin{equation}\label{condition of first derivative}
        \begin{cases}
            \, \frac{\partial}{\partial a}\operatorname{V}(\epsilon,\lambda,v,a,\epsilon S(\lambda))>0, & 0<a< \epsilon\lambda S(\lambda) \, , \\
            \, \frac{\partial}{\partial a}\operatorname{V}(\epsilon,\lambda,v,a,\epsilon S(\lambda) )<0, &  \epsilon\lambda S(\lambda)<a <\lambda \, , \\
            \, \frac{\partial}{\partial b}\operatorname{V}(\epsilon,\lambda,v,\epsilon \lambda S(\lambda) ,b)<0, & 0<b< \epsilon S(\lambda) \, , \\
            \, \frac{\partial}{\partial b}\operatorname{V}(\epsilon,\lambda,v, \epsilon \lambda S(\lambda),b)>0, & \epsilon S(\lambda)<b <1 \, . 
        \end{cases}
    \end{equation}
(Neither the point~$\big(\epsilon\lambda S(\lambda), \epsilon S(\lambda)\big)$ nor the four boundary cases are problematic, because~$\operatorname{V}$ is continuous at the concerned points. For example, $[0,\e \lambda S(\lambda)]\ni a \mapsto \operatorname{V}(\epsilon,\lambda,v,a,\epsilon S(\lambda))$ is continuous on the boundary, since $\lambda_1=\tfrac{\lambda -a}{1-\e S(\lambda)}$ is continuous at $a \in \{0,\e \lambda S(\lambda) \}$, and $h(\cdot,w):[0,1] \lora [0,1]$ in~(\ref{constrained value function}) is also continuous.) 
    
Set $\omega = \frac{\lambda-a}{1-\epsilon S(\lambda)}$. This is the value $\lambda_1 = \lambda_1\big(a,\e S(\lambda)\big)$ of Maxine's fortune as the second turn arrives in the $\e$-leisurely game if, at the first turn, Maxine stakes~$a$ in response to Mina's putatively endorsed~$\epsilon S(\lambda)$. Maxine's choice of $a = \lambda \epsilon S(\lambda)$ would yield $\omega = \lambda$, with $\omega-\lambda$ falling from positive to negative as $a$ rises through this value. Algebraically, this is easy to check, since $\omega - \lambda = \frac{\lambda \e S(\lambda) - a}{1 - \e S(\lambda)}$, an expression whose denominator is safely positive, since $\e < 1$ and $S \leq 1$.  It is also worth noting that~$\omega$ is positive in the regime~(\ref{condition of first derivative}) in question, because $a <\lambda$ there.

Consequently, we see that the first two inequalities in  \eqref{condition of first derivative}  are implied by 
\begin{equation}\label{1st derivative inequality}
    (\omega-\lambda) \tfrac{\partial}{\partial a} \operatorname{V}(\epsilon,\lambda,v,a,\epsilon S(\lambda)) > 0 \, , \quad \text{for } \omega \neq \lambda \, ,
\end{equation}
an estimate having the virtue of being a lower bound for all $\omega \in (0,\infty) \setminus \{ \lambda \}$, this achieved by the  sign change in $\omega - \lambda$ engineered by introducing the variable~$\omega$.

Something similar can be done in regard to the latter two inequalities in~(\ref{condition of first derivative}). Reflecting a form of player role-reversal, we take $\rho = \lambda^{-1}$
and, in counterpart to~$\omega$, introduce 
$$
\eta  \, = \,  \Big(\frac{\lambda - \epsilon\lambda S(\lambda)}{1-b}\Big)^{-1} \, = \,  \frac{\rho-\rho b}{1-\epsilon T(\rho)} \, ,
$$
where the latter equality is due to $T(\rho) = S(\lambda)$. There is a sign change in $\rho - \eta$ at $\eta$, with the upshot that the latter two inequalities we have mentioned are implied by
\begin{equation}\label{2nd derivative inequality}
    (\rho-\eta) \tfrac{\partial}{\partial b} \operatorname{V}(\epsilon,\lambda,v,\epsilon \lambda S(\lambda),b )>0 \, ,\quad \text{for } \eta \neq \rho \, .
\end{equation}    
Our task then is to prove the last two labelled bounds. We begin with  \eqref{1st derivative inequality}.

Writing $h'$ for the partial derivative in the first coordinate, note that   $\tfrac{\partial}{\partial a}h(\lambda_1,v) = -\frac{1}{1-b}h'(\lambda_1,v)$. Using this and 
 \eqref{constrained value function}, we find that 
 \begin{eqnarray}
         \tfrac{\partial}{\partial a} \operatorname{V}(\epsilon,\lambda,v,a,b) & = & \frac{\epsilon b}{(a+b)^2}\Big( h(\lambda_1,v_+) - h(\lambda_1,v_-) \Big) \, - \, \frac{1-\epsilon}{1-b}h'(\lambda_1,v) \label{e.abderivative} \\
         & & \qquad 
         - \, \frac{\epsilon}{1-b}\Big(\frac{a}{a+b}h'(\lambda_1,v_+)+\frac{b}{a+b} h'(\lambda_1,v_-) \Big) \, , \nonumber
\end{eqnarray}
where recall that  $\lambda_1 = \frac{\lambda-a}{1-b}$.

The right-hand expression has dependence on the values of $h$ or its derivative at the vertices $v$, $v_-$, and $v_+$.
However,  Lemma~\ref{l.bsolution} relates $h'(\lambda,\cdot)$ at these three locations. 
As the next result records, it may be invoked to find a lower bound on the quantity in~\eqref{1st derivative inequality}  that is expressed purely in terms of the stake function $S(\cdot) = S(1,\cdot,v)$ at the vertex~$v$.

\begin{lemma}\label{l.pdlb}
When $\omega \neq \lambda$, the quantity $(\omega - \lambda) \tfrac{\partial}{\partial a} \operatorname{V}\big( \epsilon,\lambda,v,a,\e S(\lambda) \big)$ is bounded below by  
$$ 
h'(\omega,v) (\omega-\lambda)^2 \big( \tfrac{\lambda - \omega}{\epsilon} + S(\lambda)  (1+\omega)\big)^{-2} \big( 1 - \e S(\lambda) \big)^{-1} \cdot \alpha \, , 
$$
 where $\alpha$ equals
\begin{equation}\label{e.alpha}
  S(\lambda)(1+\omega)^2\frac{S(\omega)-S(\lambda)}{\omega-\lambda} + S(\lambda)(1+\omega)\Big(\frac{1}{\epsilon}-S(\omega)\Big) +\frac{1-\epsilon}{\epsilon} \Big(\frac{\lambda - \omega}{\epsilon} + S(\lambda)(1+\omega) \Big)  \, .
\end{equation}
\end{lemma}
{\em Remark.} The lemma introduces~$\alpha$. Once this result is proved,   \eqref{1st derivative inequality} will follow by showing that this quantity is positive, since
 $h'(x,v)>0$ for any $(x,v) \in (0,\infty)\times \openmac$.
 
{\bf Proof of Lemma~\ref{l.pdlb}.}
Note that  $a = \lambda - \omega + \epsilon S(\lambda)\omega$ and $a+\epsilon S(\lambda) = \lambda - \omega + \epsilon S(\lambda)(1+\omega)$. Moreover, $\omega = \frac{\lambda-a}{1-\epsilon S(\lambda)}$ equals $\lambda_1$ when we take $b = \e S(\lambda)$ in~(\ref{e.abderivative}).  Thus we have
\begin{equation*}
\begin{aligned}
    \tfrac{\partial}{\partial a}\operatorname{V}(\epsilon,\lambda,v,a,\epsilon S(\lambda) ) \, = \, \frac{\epsilon^2 S(\lambda)}{\Big(\lambda - \omega + \epsilon S(\lambda)(1+\omega)\Big)^2}\Big( h(\omega,v_+) - h(\omega,v_-) \Big) - \frac{(1-\epsilon)}{1-\epsilon S(\lambda)}h'(\omega,v) \\
     \qquad - \, \frac{\epsilon}{1-\epsilon S(\lambda)}\frac{\epsilon S(\lambda)}{\lambda - \omega + \epsilon S(\lambda)(1+\omega)}\Big( \omega h'(\omega,v_+) + h'(\omega,v_-) \Big)\\
    - \, \frac{\epsilon}{1-\epsilon S(\lambda)}\frac{1}{\lambda - \omega + \epsilon S(\lambda)(1+\omega)}(\lambda-\omega)h'(\omega,v_+) \, .
\end{aligned}    
\end{equation*}

By Lemma~\ref{l.bsolution}, we know that $(1 - S(\lambda) )   \tfrac{\partial}{\partial  \lambda}h(\lambda,v) =   \tfrac{\lambda}{\lambda + 1}  \tfrac{\partial}{\partial  \lambda}h(\lambda,v_+) +  \tfrac{1}{\lambda + 1}  \tfrac{\partial}{\partial  \lambda}h(\lambda,v_-)$.  Applying this to $\omega$, and recalling that, by definition and from Corollary~\ref{c.uniqueplay}, $v_+$ and $\mc{V}_-(\lambda,v)$ do not depend on the value of $\lambda$, we get that $\omega h'(\omega,v_+) +  h'(\omega,v_-) =  (1+\omega)\big(1-S(\omega)\big) h'(\omega,v)$.

Let $A$ denote the quantity of~\eqref{1st derivative inequality}, which we seek to bound.
  Since $h$ is increasing (and differentiable), $(\omega-\lambda)^2 h'(\omega,v_+) \geq 0$.  Recall that $1-\epsilon S(\lambda) >0$ and $\lambda - \omega + \epsilon S(\lambda)(1+\omega) = a +\epsilon S(\lambda) >0$. Thus,
\begin{multline*}
       A \geq  (\omega-\lambda) \, \Bigg\{ \, \frac{\epsilon^2 S(\lambda)}{\big(\lambda - \omega + \epsilon S(\lambda)(1+\omega)\big)^2}\Big( h(\omega,v_+) - h(\omega,v_-) \Big) \, - \, \frac{(1-\epsilon)}{1-\epsilon S(\lambda)}h'(\omega,v)\\
          - \, \frac{\epsilon}{1-\epsilon S(\lambda)}\frac{\epsilon (1+\omega) S(\lambda) \big(1-S(\omega)\big)}{\lambda - \omega + \epsilon S(\lambda)(1+\omega)} h'(\omega,v) \Bigg \} \, .
\end{multline*}
Furthermore, we have $S(x) = \frac{h(x,v_+) - h(x,v_-)}{(1+x)^2h'(x,v)}$.  So we can factor out $h'(\omega,v)$ to get 
\begin{multline*}
    A \geq h'(\omega,v) (\omega-\lambda) \, \Bigg\{ \, \frac{ S(\lambda)}{\big(\frac{\lambda - \omega}{\epsilon} + S(\lambda)(1+\omega)\big)^2}(1+\omega)^2 S(\omega) - \frac{(1-\epsilon)}{1-\epsilon S(\lambda)}\\
     - \,  \frac{1}{1-\epsilon S(\lambda)}\frac{\epsilon (1+\omega) S(\lambda) \big(1-S(\omega)\big)}{\frac{\lambda - \omega}{\epsilon} + S(\lambda)(1+\omega)} \, \Bigg \} \, .
\end{multline*}
Adopt the shorthand   $\zeta = \frac{\lambda - \omega}{\epsilon} + S(\lambda)(1+\omega)$.
Then 
write the preceding right-hand side in the form $h'(\omega,v) \zeta^{-2} \big( 1 - \e S(\lambda) \big)^{-1}    (\omega-\lambda) \cdot \gamma$. As such, the new variable~$\gamma$ takes the form
$$
\gamma =(1+\omega)^2 S(\lambda)S(\omega)\big(1-\epsilon S(\lambda)\big)-(1-\epsilon)\zeta^2 \, -\, \epsilon(1+\omega)S(\lambda)\big(1- S(\omega)\big) \zeta  \, .
$$
Elementary calculations show that  $\gamma = (\omega - \lambda) \alpha$ for $\omega \neq \lambda$, where $\alpha$ is specified in Lemma~\ref{l.pdlb}. 
This proves the lemma. \qed

We now argue as promised that $\alpha$ is positive.
By $a \geq 0$, we have $\frac{\lambda - \omega}{\epsilon} + S(\lambda)(1+\omega) = \frac{1}{\epsilon}\big(a+\epsilon S(\lambda) \big) \geq S(\lambda)$,
whence~$\alpha$ is at least 
\begin{multline*}
S(\lambda) \cdot \Bigg( (1+\omega)^2\frac{S(\omega)-S(\lambda)}{\omega-\lambda} + (1+\omega)\Big(\frac{1}{\epsilon}-S(\omega)\Big)  +\frac{1-\epsilon}{\epsilon}\Bigg) \\
\geq \  S(\lambda) (1+\omega) \cdot \Bigg( -(1+\omega) \frac{\big\vert S(\omega)-S(\lambda)\big\vert}{\vert \omega-\lambda \vert} \, + \, \frac{1}{2\epsilon} \, \Bigg)  \, ,
\end{multline*} 
where the displayed bound is due to $\e \leq 1/2$ and $S \leq 1$. 
We need only show then that the last factor in parentheses is positive. A weakened form of Lemma~\ref{l.stakedifferencequotient} for $S$ in which the factor of $1+y$ is omitted  asserts that the supremum of $(1+x)\tfrac{\vert S(x) - S(y)\vert}{\vert x - y \vert}$ as $x$ and $y$ range over $(0,\infty)$ is finite. By choosing $\e > 0$ to be smaller than one-half of the reciprocal of this supremum, the desired positivity is obtained. We have proved~\eqref{1st derivative inequality}.

Next we prove~\eqref{2nd derivative inequality}. 
This inference may be reduced to the preceding argument by noting suitable symmetries, as we now explain. 
First note that Lemma~\ref{l.bsolution} is invariant under the transformation of replacing $\lambda$ by $\rho=1/\lambda$, and interchanging the goals of Mina and Maxine in the following sense: interchange $v_+$ and $v_-$, and interpret $T(\rho) =$ Stake$(1,1/\lambda,v)$ 
   \ignore{ \yujienew{(Do you mean the stake ratio in the following game after interchanging the goals? If so,  do you want to say:"interpret $T(\rho) = T(1/\lambda) = \stake(1,\lambda,v)$ as the stake ratio...?" The stake ratio should be the same as the equivalent original game after interchanging the goals, which should be Stake$(1,\lambda,v)$, right? Or do you mean this stake function Stake$(1,1/\lambda,v)$ is just defined in the new game, with a similar formula using new harmonic function $\tilde{h}$, Like $T(\rho) = \tfrac{\tilde{h}(\rho,v_-) - \tilde{h}(\rho,v_+)}{(1+\rho)^2\tfrac{\partial}{\partial \rho} \tilde{h}(\rho,v)}$? But will this cause confusion for readers?)}}as the stake ratio in the game with both fortunes multiplied by $\lambda$, so $\lambda$ for original Mina and one for Maxine; but after interchanging the goals this is the same game as one for Mina and $\lambda$ for Maxine. Then note that this is also the transformation taking~\eqref{1st derivative inequality} to~\eqref{2nd derivative inequality}, including taking $\omega$ to $\eta$, but with a global sign change; moreover, every other step in the deduction of~\eqref{1st derivative inequality}  also transforms well, including a sign change from $\frac{\partial}{\partial a} h(\lambda_1,v)$ to $\frac{\partial}{\partial b} h(\lambda_1,v)$. Hence~\eqref{2nd derivative inequality} also follows.  

It is perhaps a little hard to keep track of the transformations in the above reduction, and the reader may prefer a direct argument that follows the steps that led to~\eqref{1st derivative inequality}, the varied stake now being Mina's where before it was Maxine's. That said, giving this argument here would be rather repetitive. It appears in the corresponding place in the arXiv version of this paper,
at \href{https://arxiv.org/abs/2206.08300}{arXiv:2206.08300}.

\ignore{
It is perhaps a little hard to keep track of the transformations in the above reduction, and the reader may prefer a direct argument, next presented, that follows the steps that led to~\eqref{1st derivative inequality}, the varied stake now being Mina's where before it was Maxine's.

\ignore{\yujienew{yujie: Here are some of my understanding and thoughts about this transformation argument.  Let me know if I understand it wrongly.}

\yujienew{yujie: If I am understanding correctly, when we interchange the goal of Mina and Maxine, interchange $v_-$ and $v_+$, in Lemma~\ref{l.bsolution} and in the new game described above, we also use the new harmonic function $\tilde{h}(\rho,v)$, where $\tilde{h}(x,v) = h(1/x,v)$, which is decreasing, and leads to the sign change in~\eqref{2nd derivative inequality}.  Do you think we should mention that?}

\yujienew{For Lemma~\ref{l.bsolution}: $
 (1 - \stake(1,\lambda,v) )   \tfrac{\partial}{\partial  \lambda}h(\lambda,v) =   \tfrac{\lambda}{\lambda + 1}  \tfrac{\partial}{\partial  \lambda}h(\lambda,v_+) +  \tfrac{1}{\lambda + 1}  \tfrac{\partial}{\partial  \lambda}h(\lambda,v_-)$ becomes $
 (1 - \stake(1,\lambda,v) )   \tfrac{\partial}{\partial  \rho}\tilde{h}(\rho,v) =   \tfrac{1}{\rho + 1}  \tfrac{\partial}{\partial  \rho}\tilde{h}(\rho,v_+) +  \tfrac{\rho}{\rho + 1}  \tfrac{\partial}{\partial  \rho} \tilde{h}(\rho,v_-)$, and we do not change $b = \stake(1,\lambda,v)$ (Or say $b = S(\lambda)$ goes to $T(\rho)$), (Or say the Stake function is defined in new ways?) the transformation argument is in this way, right?}

\yujienew{When I prove the second claim~\eqref{2nd derivative inequality}, besides the following way of calculation, I also calculated with another method. In that method, I make a substitution in $\operatorname{V}$ at the start of calculation, to denote it using $\tilde{h}$. Then denote~\eqref{2nd derivative inequality} using $\tilde{h}$, $\rho,\eta$, which will have a similar form as~\eqref{1st derivative inequality} after symbolic change, interchanging $v_+$ and $v_-$ and a global sign change.
 I think this method is equivalent to Gabor's transformation argument, in essence.  So I believe this transformation argument is correct.  } }

\ignore{\textcolor{blue}{Yujie: I think this transformation argument is great, do you think we still need the following second part, the detailed but repetitive computation or not?}}


\ignore{\yujie{I think I agree. Suggested from Gabor: "Replace the second set of calculation by the following sentence? In case the Reader is not convinced by this transformation argument, they may deduce~\eqref{2nd derivative inequality} by following the same steps as above, just varying the stake of Mina instead of Maxine."  }}

Using~\eqref{constrained value function}, and noting that $\tfrac{\partial}{\partial b}h(\lambda_1,v) = \frac{\lambda-a}{(1-b)^2}h'(\lambda_1,v)$, where $\lambda_1 = \frac{\lambda-a}{1-b}$, we find that \begin{equation*}
    \begin{aligned}
                 \tfrac{\partial}{\partial b} \operatorname{V}(\epsilon,\lambda,v,a,b) \, = \, -\frac{\epsilon a}{(a+b)^2}\Big( h(\lambda_1,v_+) - h(\lambda_1,v_-) \Big) + (1-\epsilon)\frac{\lambda-a}{(1-b)^2}h'(\lambda_1,v) \\+ \,\frac{\epsilon(\lambda-a)}{(1-b)^2}\bigg(\frac{a}{a+b}h'(\lambda_1,v_+)+\frac{b}{a+b} h'(\lambda_1,v_-)  \bigg) \, .
    \end{aligned}
\end{equation*}
For $a = \epsilon\lambda S(\lambda) = \epsilon \frac{T(\rho)}{\rho}$, and $\eta = \frac{\rho-\rho b}{1-\epsilon T(\rho)}$, where $T$ and $\rho$ are defined above, we have $\rho b = \rho -\eta +\epsilon \eta T(\rho)$, $\lambda_1 = \frac{\lambda-a}{1-b} = \eta^{-1} $, and \begin{align*}
    &\frac{a}{a+b} = \frac{\epsilon T(\rho)}{\rho - \eta + \epsilon T(\rho)(1+\eta) }, \quad  \frac{b}{a+b} = \frac{\rho - \eta + \epsilon T(\rho)\eta}{\rho - \eta + \epsilon T(\rho)(1+\eta)},\\
    &\frac{a}{(a+b)^2} = \rho \frac{\epsilon T(\rho)}{\Big(\rho - \eta + \epsilon T(\rho)(1+\eta)\Big)^2 }, \quad \frac{\lambda-a}{(1-b)^2} = \frac{\rho}{\eta^2}\frac{1}{1-\epsilon T(\rho)}.
\end{align*}
Thus, we have \begin{multline*}
        \tfrac{\partial}{\partial b} \operatorname{V}(\epsilon,\lambda,v,\epsilon \lambda S(\lambda),b) \, = \, \frac{\rho}{\eta^2}\Bigg\{-\frac{\epsilon^2 T(\rho)}{\Big(\rho - \eta + \epsilon T(\rho)(1+\eta)\Big)^2}\eta^2\Big( h(\eta^{-1},v_+) - h(\eta^{-1},v_-) \Big) \\
        + \, \frac{1-\epsilon}{1-\epsilon T(\rho)}h'(\eta^{-1},v) \, + \, \frac{\epsilon}{1-\epsilon T(\rho)}\frac{\epsilon T(\rho)\eta}{\rho - \eta + \epsilon T(\rho)(1+\eta)}\Big(\eta^{-1} h'(\eta^{-1},v_+) + h'(\eta^{-1},v_-) \Big) \\
         + \, \frac{\epsilon}{1-\epsilon T(\rho)}\frac{1}{\rho - \eta + \epsilon T(\rho)(1+\eta)}\Big( (\rho-\eta)h'(\eta^{-1},v_-) \Big) \, \Bigg\} \, .
\end{multline*}

By Lemma~\ref{l.bsolution}, 
we have $\eta^{-1} h'(\eta^{-1},v_+) +  h'(\eta^{-1},v_-) =  (1+\eta^{-1})\big(1-T(\eta)\big) h'(\eta^{-1},v)$.

\begin{lemma}\label{l.pdlb.b}
The quantity $ (\rho-\eta) \tfrac{\partial}{\partial b} \operatorname{V}(\epsilon,\lambda,v,\epsilon \lambda S(\lambda),b)$
 is bounded below by  
 $$
h'(\eta^{-1},v) \frac{\rho}{\eta^2}(\eta-\rho)^2 \Big(\frac{\rho -\eta}{\epsilon}+T(\rho)(1+\eta) \Big)^{-2}\big(1-\epsilon T(\rho)\big)^{-1} \cdot \beta
$$
where $\beta$ equals
\begin{equation}\label{e.beta}
T(\rho)(1+\eta)^2\frac{T(\eta)-T(\rho)}{\eta-\rho} + T(\rho)(1+\eta)\Big(\frac{1}{\epsilon}-T(\eta)\Big) +\frac{1-\epsilon}{\epsilon} \Big(\frac{\rho - \eta}{\epsilon} + T(\rho)(1+\eta) \Big) \, .
\end{equation}
\end{lemma}
{\bf Proof.} Write $B$ for the quantity to be bounded below. 
Since $h$ is increasing, $(\rho-\eta)^2 h'(\omega,v_-) \geq 0$. And we also note that $1-\epsilon T(\rho) >0$ and $\rho - \eta + \epsilon T(\rho)(1+\eta)  = \rho b + \epsilon T(\eta) > 0$.  Hence,
\begin{multline}\label{Inequality of B-1}
         B \, \geq \, \frac{\rho}{\eta^2}(\rho-\eta)\, \Bigg\{ \, -\frac{\epsilon^2 T(\rho)}{\big(\rho - \eta + \epsilon T(\rho)(1+\eta)\big)^2}\eta^2\Big( h(\eta^{-1},v_+) - h(\eta^{-1},v_-) \Big) \, + \, \frac{1-\epsilon}{1-\epsilon T(\rho)}h'(\eta^{-1},v) \\
         + \, \frac{\epsilon}{1-\epsilon T(\rho)}\frac{\epsilon T(\rho)\eta}{\rho - \eta + \epsilon T(\rho)(1+\eta)}(1+\eta^{-1})\big(1-T(\eta)\big) h'(\eta^{-1},v) \, \Bigg\} \, .
\end{multline}

Since $S(x) = \frac{h(x,v_+) - h(x,v_-)}{(1+x)^2h'(x,v)}$, we know that $T(x) = S(x^{-1}) = x^2\frac{h(x^{-1},v_+) - h(x^{-1},v_-)}{(1+x^2)h'(x^{-1},v)}$.  Then we can factor out $h'(\eta^{-1},v)$ in~\eqref{Inequality of B-1} to get
\begin{multline*}
        B \geq \frac{\rho}{\eta^2}h'(\eta^{-1},v)(\eta-\rho) \, \Bigg\{ \, \frac{ T(\rho)}{\big(\frac{\rho - \eta}{\epsilon} + T(\rho)(1+\eta)\big)^2}(1+\eta^2)T(\eta) \, - \, \frac{1-\epsilon}{1-\epsilon T(\rho)} \\
         - \, \frac{\epsilon}{1-\epsilon T(\rho)}\frac{T(\rho)}{\frac{\rho - \eta}{\epsilon} + T(\rho)(1+\eta)}(1+\eta)\big(1-T(\eta)\big) \,  \Bigg\} \, .
\end{multline*}
Write this right-hand side in the form $\rho \eta^{-1} h'(\eta^{-1},v) \psi^{-2} \big( 1 - \e T(\rho) \big)^{-1}  (\eta-\rho) \cdot \kappa$ with  $\psi$ set equal to $\frac{\rho - \eta}{\epsilon} + T(\rho)(1+\eta)$. As such,
$$
\kappa =(1+\eta)^2 T(\rho)T(\eta)\big(1-\epsilon T(\rho)\big)-(1-\epsilon)\psi^2 \, -\, \epsilon(1+\eta)T(\rho)\big(1- T(\eta)\big) \psi  \, .
$$
Elementary calculations show that  $\kappa = (\eta - \rho) \beta$ for $\eta \neq \rho$, where $\beta$ is specified in Lemma~\ref{l.pdlb}. 
This proves the lemma. \qed 

Since $h'(x,v) > 0$, $\beta > 0$ implies \eqref{2nd derivative inequality}. But $\beta$ is at least
\begin{multline*}
(\eta-\rho)^2 T(\rho) \cdot \Bigg( (1+\eta)^2\frac{T(\eta)-T(\rho)}{\eta-\rho} + (1+\eta)\Big(\frac{1}{\epsilon}-T(\eta)\Big)  +\frac{1-\epsilon}{\epsilon}\Bigg) \\
     \geq \   (\eta-\rho)^2 T(\rho) (1+\eta) \cdot \Bigg( -(1+\eta) \frac{\big\vert T(\eta)-T(\rho)\big\vert}{\vert \eta-\rho \vert} \, + \, \frac{1}{2\epsilon} \, \Bigg)  \, ,
\end{multline*} 
the bound due to $\e \leq 1/2$ and $T \leq 1$. Lemma~\ref{l.stakedifferencequotient} for $T$  implies that  
$$
\sup \big\{ (1+x)\tfrac{\vert T(x) - T(y)\vert}{\vert x - y \vert}: x,y \in [0,\infty) \big\} < \infty \, .
$$
 When $\e > 0$ is less than  one-half of the reciprocal of this supremum, the last-displayed bracketed quantity is positive. 
 Thus $\beta > 0$ and~\eqref{2nd derivative inequality} is obtained.
}

In obtaining \eqref{condition of first derivative} via 
\eqref{1st derivative inequality} and \eqref{2nd derivative inequality},
we have proved that $(a^*,b^*) = \big(\epsilon\lambda S(\lambda),\epsilon S(\lambda) \big)$ is a strict global minimax of $\operatorname{V}(\epsilon,\lambda,v,a,b)$. 
It remains to show that this global minimax is unique. Let $(a_1,b_1)$ be another. Without loss of  generality, suppose that $a_1 \neq a^*$. We invoke that $(a^*,b^*)$ is a strict global minimax in the second bound as we write 
\begin{align*}
\operatorname{V}(\epsilon,\lambda,v,a_1,b_1)\leq \operatorname{V}(\epsilon,\lambda,v,a_1,b^*)<\operatorname{V}(\epsilon,\lambda,v,a^*,b^*)\leq \operatorname{V}(\epsilon,\lambda,v,a^*,b_1)\leq \operatorname{V}(\epsilon,\lambda,v,a_1,b_1).
\end{align*}
A contradiction! The point $(a^*,b^*) = \big(\epsilon\lambda S(\lambda),\epsilon S(\lambda) \big)$ is a unique global minimax and the proof of Proposition~\ref{p.globalconstrainedminimax} is complete. \qed

Before making the dependence of $\e_0$ on the tree $T$ more explicit, we derive the main result of the section. Given Proposition~\ref{p.globalconstrainedminimax}, this will be a rather standard game theory argument.

{\bf Proof of Theorem~\ref{t.leisurely.theorem}.}
We use $M_n(\lambda,v,S_-,S_+)$ to denote the mean payment of $\game_n(\epsilon,\lambda,v)$ under the strategy pair $(S_-,S_+)$.  And in Definition~\ref{d.conforming}, we expand the notion of ``conforming strategy" so that it is also defined for the finite horizon game $\game_n(\epsilon,\lambda,v)$ for any $n$.

The proposition will be proved by induction on $n$.  We write $\mathsf{Hyp}(n,i)$, $n \in \nwithoutzero$ and $i \in \{1,2\}$, to refer to Theorem~\ref{t.leisurely.theorem}(1) with index value $n$. The base case may be viewed as $n=0$: $\game_0(\e,\lambda,v)$ is moveless, with $\pay=h(\lambda,v)$, so that  $\mathsf{Hyp}(0,1)$ holds trivially; we may formally specify  $\mathsf{Hyp}(0,2)$ to be the vacuous statement. 

Now we assume that $\mathsf{Hyp}(n-1,1)$ and $\mathsf{Hyp}(n-1,2)$ hold.  Take $\epsilon_0$ to be the minimum of this value as specified in Proposition~\ref{p.globalconstrainedminimax} and as it is specified in $\mathsf{Hyp}(n-1,1)$ and $\mathsf{Hyp}(n-1,2)$. (Actually by the inductive argument, we know our choice of $\epsilon_0$ is just the one specified in Proposition~\ref{p.globalconstrainedminimax}.)  

Take $\e \in (0, \epsilon_0]$
and
let $(\Pminus,\Pplus)$ be a pair of conforming strategies.  We note that $\Pminus$ and $\Pplus$ are well defined in $\game_n(\epsilon,\lambda,v)$ for any $n$.  So in the ensuing argument, we may use the notation $\Pminus$ or $\Pplus$ to denote a strategy in either $\game_n(\epsilon,\lambda,v)$ or $\game_{n-1} (\epsilon,\lambda,v)$, without risk of confusion.

                                                      By $\mathsf{Hyp}(n-1,2)$, we know that for 
                                                      $(\lambda,v) \in (0,\infty) \times \openmac$, $(\Pminus,\Pplus)$ is a strict Nash equilibrium of $\game_{n-1} (\epsilon,\lambda,v)$ in the sense of \eqref{strict Nash equilibrium for conforming}.  We will show that this pair is also a strict Nash equilibrium of $\game_{n}(\epsilon,\lambda,v) $ in the same sense.

       We can divide the $\game_n(\epsilon,\lambda,v)$ into two parts: what happens at the first turn, and the remainder.  Let $S_+$ be any strategy of Maxine in $\game_n(\epsilon,\lambda,v)$ and denote by $S_+\rangle_2^n$ the resulting strategy of Maxine in the game from the second turn. We may write $S_+$ in the form $\big((A,u), S_+\rangle_2^n \big)$; this means that, at the first turn, Maxine will stake $A$ and nominate a move to $u$, and, in the remaining game, she will stick to the strategy $S_+\rangle_2^n$.  Here the quantity $A$ is permitted to be random since $S_+$ may be a mixed strategy, though note that  $u$ remains deterministic, in accordance with our specification of random strategies in  Definition~\ref{d.strategy}.  Let $S_-$ be any strategy of Mina in $\game_n(\epsilon,\lambda,v)$. We write $S_-=\big( (B,u),S_-\rangle_2^n\big)$, where $S_-\rangle_2^n$ is defined similarly as was $S_+\rangle_2^n$. 

Conditionally on the fortune ratio and counter location at the start of the second turn of $\game_{n} (\epsilon,\lambda,v)$ under the pair of strategies $(\Pminus,S_+)$ assuming a given value~$(\mu,w)$, the game remaining from this turn is a copy of the $(n-1)$-horizon $\game_{n-1}(\epsilon,\mu,w)$. And  $\Pminus\rangle_2^n$ and $S_+\rangle_2^n$ may be regarded as strategies of Mina and Maxine in $\game_{n-1}(\epsilon,\mu,w)$.  The conditional mean payment of $\game_n(\epsilon,\lambda,v)$ under $(\Pminus,S_+)$ is equal to $M_{n-1}\big(\mu,w,\Pminus\rangle_2^n,S_+\rangle_2^n\big)$, and $\Pminus\rangle_2^n$ is equal to $\Pminus$ in~$\game_{n-1} (\epsilon,\mu,w)$.  

Write $S(\lambda) = \stake(1,\lambda,v)$ and $\lambda_1 = \frac{\lambda-A}{1-\epsilon S(\lambda)}$. Note that
\begin{equation*}
    \begin{aligned}
        M_n\big(\lambda,v,\Pminus,S_+\big) & \, = \, \E \, \Bigg[ \,  (1-\epsilon)M_{n-1}\big (\lambda_1, v, \Pminus,S_+\rangle_2^n\big) + \epsilon \bigg(\frac{A}{A +    \epsilon S(\lambda)}M_{n-1}\big(\lambda_1, u,\Pminus,S_+\rangle_2^n\big) \\
        &\qquad \qquad \qquad \qquad \qquad \qquad \qquad \qquad  + \, \frac{\epsilon S(\lambda)}{A+    \epsilon S(\lambda)}M_{n-1}\big(\lambda_1, v_-,\Pminus,S_+\rangle_2^n\big) \bigg) \Bigg] 
    \end{aligned}
\end{equation*}
where  the expectation arises from the randomness of $A$ in strategy $S_+$. 
We find then that
\begin{equation*}
    \begin{aligned}
        M_n\big(\lambda,v,\Pminus,S_+\big) 
        &\, \leq \, \E \, \Bigg[ \, (1-\epsilon)M_{n-1} \big(\lambda_1, v, \Pminus,\Pplus\big) + \epsilon \bigg(\frac{A}{A+    \epsilon S(\lambda)}M_{n-1}\big(\lambda_1, u,\Pminus,\Pplus\big) \\
        &\qquad \qquad \qquad \qquad \qquad \qquad \qquad \qquad + \, \frac{\epsilon S(\lambda)}{A+    \epsilon S(\lambda)}M_{n-1}\big(\lambda_1, v_-,\Pminus,\Pplus\big) \bigg)\Bigg] \\
        & \, = \,  \E \, \Bigg[ \, (1-\epsilon)h(\lambda_1,v) + \epsilon \bigg(\frac{A}{A+    \epsilon S(\lambda)}h(\lambda_1,u) + \frac{\epsilon S(\lambda)}{A+    \epsilon S(\lambda)}h(\lambda_1,v_-)\bigg) \Bigg] \\
        & \, \leq \, \E \, \Bigg[ \, (1-\epsilon)h(\lambda_1,v) + \epsilon \bigg(\frac{A}{A+    \epsilon S(\lambda)}h(\lambda_1,v_+) + \frac{\epsilon S(\lambda)}{A+    \epsilon S(\lambda)}h(\lambda_1,v_-) \bigg) \Bigg]\\
        & \, = \, \E \operatorname{V}(\epsilon,\lambda,v,A,\epsilon S(\lambda)) \, \leq \, h(\lambda,v) \, .
    \end{aligned}
\end{equation*}
Here  the first bound is due to $(\Pminus,\Pplus)$ being a pure Nash equilibrium in $\game_{n-1}(\epsilon,\lambda,v)$ for any $(\lambda,v) \in (0,\infty) \times \openmac$, and the equality obtains only if $S_+\rangle_2^n$ is conforming against $\Pminus$ by $\mathsf{Hyp}(n-1,2)$.   The next inequality follows from $h(w,u) \leq h(w,v_+) $, the equality being achieved only if $u=v_+$.  The final bound, which appears in the last line, holds by Proposition~\ref{p.globalconstrainedminimax}, with equality only when $a = \lambda \stake(\epsilon,\lambda,v)$ almost surely.  Thus we find that $M_n(\lambda,v,\Pminus,S_+) = h(\lambda,v)$ only if $S_+$ is conforming against $\Pminus$ .  In other words, \begin{align*}
    M_n(\lambda,v,\Pminus,S_+) < h(\lambda,v), \quad \text{if $S_+$ is not conforming against $\Pminus$ }.
\end{align*}

Now let $S_-$ be any strategy of Mina. With an almost verbatim argument, we learn that 
\begin{equation}\label{conforming strategy inequality}
    \begin{aligned}
    &M_n(\lambda,v,S_-,\Pplus)> h(\lambda,v), \quad \text{if $S_-$ is not conforming against $\Pplus$ }.
    \end{aligned}
\end{equation}

In this way, we have proved that $\game_n(\epsilon,\lambda,v)$ has value $h(\lambda,v)$; and moreover,  any pair of conforming strategies forms a strict pure Nash equilibrium of $\game_n(\epsilon,\lambda,v)$ in the sense of \eqref{strict Nash equilibrium for conforming}.

Now let $(P_-^0,P_+^0)$ be another Nash equilibrium.  Suppose that at least one of $P_-^0$ and $P_+^0$ is not conforming. Without losing generality, we may assume that $P_-^0$ is not conforming. Then,
\begin{equation*}
    h(\lambda,v) = M_n(\lambda,v,P_-^0,P_+^0) \geq M_n(\lambda,v,P_-^0, \Pplus) > M_n(\lambda,v,\Pminus,\Pplus) = h(\lambda,v).
\end{equation*}
The first inequality holds because $(P_-^0,P_+^0)$ is a pure Nash equilibrium;  the second is due to  \eqref{conforming strategy inequality}.  Thus we reach a contradiction, and find that 
both of $P_-^0$ and $P_+^0$ must be conforming.  We have proved $\mathsf{Hyp}(n,1)$ and $\mathsf{Hyp}(n,2)$
and have thus completed the inductive proof of Theorem~\ref{t.leisurely.theorem}. \qed

We end this section by carrying out the task needed to prove the explicit value for $\e_0$ claimed in Theorem~\ref{t.leisurely}: finding a lower bound on this quantity in Theorem~\ref{t.leisurely.theorem}.  


\begin{lemma}\label{l.estimate}
    With $S(x) = \stake(1,x,v)$ for $v\in V_0$, 
    \begin{equation*}
        \vert S'(x) \vert \, \leq \, \frac{  \vert V \vert^{2\msdepth+3}}{2(1+x)^2}  \quad \text{for } x \geq 0 \, .
    \end{equation*} 
    where $\msdepth \in \N$ denotes the MS-depth of~$T$.
\end{lemma}

{\bf Proof.} Note that $\msdepth$
 is equal to the maximum value of the index  $k$ appearing in the journey data
$\big\{ (s_i,d_i): i \in \llbracket 0, k \rrbracket \big\}$ as $v$ varies over $\openmac$.

Let $\mc{C}$ denote the set of polynomials, including constants, all of whose coefficients are positive.
  Put a partial order $\preccurlyeq$ on $\mc{C}$ via
  \begin{eqnarray*}
     \sum_{i=0}^{m}a_ix^i \, \preccurlyeq \, \sum_{i=0}^{n}b_ix^i  \quad \text{if } m\leq n , \, \text{and } a_i \leq b_i \, \, \text{for } i \in \llbracket 0,m \rrbracket \,.
\end{eqnarray*}
To given $v \in V$ is associated the journey data $\big\{ (s_i,d_i): i \in \llbracket 0, k \rrbracket \big\}$, with $k \leq \Kdepth$ by the definition of~$K$. Let  $Q_{n,i}$ and $P_{\ell} $ be the polynomials defined in the
derivation of Claim~I in the proof of Lemma~\ref{l.stakedifferencequotient}. Then 
\begin{align*}
    Q_{s_i,d_i} \, &= \, \sum_{\ell=0}^{d_i-1} \sum_{j=0}^{s_i-d_i-1}x^jP_{s_i-d_i-1+\ell-j } \, = \, \sum_{\ell=0}^{d_i-1} \sum_{j=0}^{s_i-d_i-1} \big( P_{s_i-d_i-1+\ell} - P_{j-1} \big) \\
    &\preccurlyeq \, \sum_{\ell=0}^{d_i-1} \sum_{j=0}^{s_i-d_i-1} P_{s_i-2} \, \preccurlyeq \, (s_i-1)^2P_{s_i-2}\, .
     \end{align*}
Here we used $1\leq d_i \leq s_i-1$ and $s_i \geq 2$, these due to the definition of $Q_{n,i}$. Next we note that $P_m P_{m'} = \sum_{i=0}^m\sum_{j=0}^{m'} x^{i+j} \preccurlyeq  \min\{m+1,m'+1\} P_{m+m'} $. So for a sequence $P_{a_1},\dots, P_{a_n}$ ordered so that $a_1 \geq \cdots \geq a_n$, we find that, for $\ell \in \intint{n}$, 
$$ 
(a_\ell+1)\prod_{i=1}^n P_{a_i} \, \preccurlyeq \, (a_\ell+1)\prod_{i=2}^n(a_i+1) \cdot P_{\sum_i a_i} \, \preccurlyeq \, \prod_{i=1}^n(a_i+1) \cdot P_{\sum_i a_i} \, . 
$$ 
Letting $N_i = s_i-2 + \sum_{\ell \neq i}(2s_\ell-d_\ell-2) $ for $i \in \llbracket 0, k \rrbracket$, we see that
\begin{align*}
        & (s_i-1)^2 P_{s_i-2} \prod_{\ell\neq i}P_{s_\ell-1} P_{s_\ell-d_\ell-1} \, \preccurlyeq \, (s_i-1)^2\Big(\prod_{\ell \neq i}(s_\ell-d_\ell)s_\ell \Big) P_{N_i} \, \preccurlyeq \, \Big(\prod_{\ell=0}^ks_\ell^2 \Big) P_{N_i} \, .
\end{align*}

By~\eqref{stake formula as rational function}, $S(x)$ takes the rational form $f(x)/g(x)$.  Note that $$g(x) \, =  \, (1+x)\sum_{i=0}^{k}Q_{s_i,d_i}(x) \prod_{\ell\neq i}P_{s_\ell-1}(x)P_{s_\ell-d_\ell-1}(x) \, , $$
and $\vert V \vert \geq \sum_{i=0}^k s_i $,  so that any coefficient of $g(x)$ is no more than \begin{align*}
      2(k+1)\prod_{\ell=0}^ks_\ell^2 \, \leq \, 2 (k+1) \bigg(\frac{\sum_{\ell=0}^k s_\ell }{k+1}\bigg)^{2k+2} \, \leq \, \frac{2}{(k+1)^{2k+1}} \vert V \vert^{2k+2} \, \leq \, \tfrac{1}{4}\vert V \vert^{2\Kdepth+2} \, .
\end{align*}

The last bound is due to $1\leq k \leq K$ and offers a bound on the coefficients of $g(x)$ that is uniform in $v \in \openmac$.

Writing $f(x) = \sum_{i=0}^n a_i x^i$ and $g(x) = \sum_{i=0}^n b_i x^i$, we find that $n \leq \vert V\vert-1$ and $1\leq a_i \leq b_i \leq \frac{1}{4}\vert V \vert^{2\Kdepth+2}$  for $i \in \llbracket 0, n\rrbracket $, where $a_i \leq b_i$ is direct from \eqref{stake formula as rational function} and $P_{s_k-1}\preccurlyeq (1+x)Q_{s_k,d_k} $.  We can further write \begin{align*}
    S(x) = \frac{a_n}{b_n}+\sum_{i=0}^{n-1}S_i(x), \quad \text{where }S_i(x) = \frac{f_i(x)}{g(x)}, \, \text{and } f_i(x) =\Tilde{a_i}x^i, \, \Tilde{a_i} = a_i-\frac{a_n}{b_n}b_i \, .
\end{align*}

For $S_i(x) = f_i(x)/g(x)$, we have $\vert\Tilde{a_i}\vert \leq b_i \leq \frac{1}{4}\vert V \vert^{2\Kdepth+2}$.  Then for $x\geq 0$, we calculate\begin{align*}
    \Big\vert f_i'(x)g(x) - f_i(x)g'(x) \Big\vert & \, = \, \bigg\vert i\Tilde{ a_i}x^{i-1}\sum_{j=0}^n b_j x^j  - \Tilde{a_i} x^{i}\sum_{j=1}^n  j b_j x^{j-1} \bigg\vert \, = \,  \bigg\vert \Tilde{a_i} x^{i-1}\sum_{j=0}^n (i-j)b_j x^j \bigg\vert \\
    & \, \leq \, n \vert \Tilde{ a_i} \vert x^{i-1} g(x) \, \leq \, \vert V \vert \vert \Tilde{a_i} \vert x^{i-1} g(x) \, .
\end{align*}
And we note that $(1+x)^2 \vert \Tilde{a_i} \vert x^{i-1} \preccurlyeq \frac{1}{2}\vert V \vert^{2\Kdepth+2} g(x)$, for $i \in \intint{n-1}$, so that 
$$
(1+x)^2\big\vert f_i'(x)g(x) - f_i(x)g'(x) \big\vert \, \leq \, \tfrac{1}{2} \vert V \vert^{2K+2} g(x)^2, \quad \text{for } x\geq 0 \, ,
$$
whence $\vert S_i'(x) \vert \leq \frac{ \vert V \vert^{2K+2}}{2(1+x)^2}$, for $x\geq 0$.
Thus,  we find that
\begin{equation*}
    \vert S'(x)\vert \leq \sum_{i=0}^{n-1}\vert S_i'(x)\vert \leq n \frac{  \vert V \vert^{2K+2}}{2(1+x)^2} \leq \frac{ \vert V \vert^{2\Kdepth+3}}{2(1+x)^2} 
\end{equation*}
and complete the proof of Lemma~\ref{l.estimate}. \qed

In this way, 
the hypothesis $(1+x)^2 \vert S'(x) \vert \leq C$ in Claim~{\rm II} in the proof of Lemma~\ref{l.stakedifferencequotient} 
 is valid with $C = \frac{1}{2} \vert V \vert^{2K+3}$. This result thus gives that $(1+y)\big\vert\frac{S(y)-S(x)}{y-x} \big\vert \leq C$ for $x,y \in (0,\infty)$.  And this bounds holds equally with $T$ in place of~$S$.  The choice $\e_0 = (2C)^{-1}$ then validates
 Proposition~\ref{p.globalconstrainedminimax} and Theorem~\ref{t.leisurely.theorem} because this value enables the derivation of~(\ref{1st derivative inequality}) and~(\ref{2nd derivative inequality}) in the  proof of Proposition~\ref{p.globalconstrainedminimax}.


We 
obtain the following consequence.

\begin{corollary}\label{c.estimateofe0}
   The value of $\e_0$ in Theorem~\ref{t.leisurely.theorem} may be taken equal to $\vert V \vert^{-2\Kdepth-3}/4$, with $\Kdepth$ the MS-depth of~$T$. 
\end{corollary}
Recall that $\boundarymac$ is the set of boundary vertices of $(V,E,{\bf 1}_r)$.
Since  $\vert \boundarymac \vert \geq K+1$, we may equally take $\e_0 = \vert V \vert^{-2\vert \boundarymac \vert-1}/4$, which is the simple form seen in Theorem~\ref{t.leisurely}. 

\section{Proving the main results from their finite horizon counterparts}\label{s.finitetoinfinite}

Here we prove Theorems~\ref{t.leisurely} and~\ref{t.nashform} and Corollary~\ref{c.gameplay}.

In speaking of $\game(\e,\lambda,v)$, we assume when needed that $\e \in (0,\e_0)$, where $\e_0$ is specified in Theorem~\ref{t.leisurely.theorem}.


Definition~\ref{d.conforming} concerns \jointstrategies that result in gameplay where both players conform at every turn with probability one.  It is useful to consider explicit strategies that satisfy the stronger condition that they demand that the player conform at every turn regardless of the strategy adopted by her opponent. These are the strongly conforming strategies that we now specify. 
\begin{definition}\label{d.stronglyconforming}
The strongly conforming strategy for Maxine in $\game(\e,\lambda,v)$ demands that, when $\stateofplay$ equals an arbitrary $(\mu,w) \in (0,\infty) \times \openmac$, Maxine stakes $\mu \cdot \stake(\e,\mu,w)$ and nominates the move $w_+$, regardless of the index of the turn in question. If $(\mu,w) \in \{ \infty \} \times \openmac$ (so that Mina is bankrupt), Maxine stakes a given positive but finite sum and nominates~$w_+$. Strongly conforming strategies for Mina are indexed by $\theta$ in the index set $\Theta$ specified in Definition~\ref{d.thetafirst}.
The strategy so indexed by $\theta$ demands that, when $\stateofplay = (\mu,w) \in (0,\infty) \times \openmac$ and the turn index is $i \in \nwithoutzero$, Mina stakes $\stake(\e,\mu,w)$ and nominates the move $\theta (\mu,w,i-1)$. If $(\mu,w) \in \{ 0 \} \times \openmac$, Mina stakes a given quantity on $(0,1)$---one-half, say---and nominates a definite move away from the root---let's say $\theta(1,w,i-1)$. 
\end{definition}
The next result shows that it is a mistake to deviate from conforming play if one's opponent conforms.
\begin{proposition}\label{p.mistakecost}
Let $\minaconf$ and $\maxineconf$ denote strongly conforming strategies for Mina and Maxine in $\game(\e,\lambda,v)$. Let  $\sminusinfinity \in \mc{S}_-$ be a strategy for Mina that is not conforming against  $\maxineconf$. Then
$$
 M \big( \sminusinfinity , \maxineconf \big) > h(\lambda,v) \, .
$$
And likewise let  $\splusinfinity \in \mc{S}_+$ be a strategy for Maxine that is not conforming against  $\minaconf$. Then
$$
 M \big( \minaconf , \splusinfinity \big) < h(\lambda,v) \, .
$$
\end{proposition}
This result is the principal element in the passage from finite to infinite horizon that we are undertaking.  Recall from Section~\ref{s.specifying} that the finish time of $\game(\e,\lambda,v)$
is denoted by~$F$. A critical element in making this passage is gaining an understanding that under suitable strategy pairs in the infinite-horizon game, play is certain to end at a finite time: in other words, that $F < \infty$, or $X_i \in \boundarymac$ for some $i \in \nwithoutzero$, almost surely.
 If both players conform, this is easy enough: the counter evolution is in the sense of Definition~\ref{d.theta} a lazy walk $X_\theta$ for some $\theta \in \Theta$, so that~(\ref{e.ftail}) in the proof of Proposition~\ref{p.nashformprep}(2) bounds the tail of the finish time. Indeed, it would not be hard to obtain an exponentially decaying bound on the tail in this case. When one player conforms but the other may not, proving the finiteness of the finish time is more delicate. Next is a result that does so, with Mina the conforming party. There is an additional hypothesis that regularises Maxine's play when her fortunes are good.

\begin{definition}\label{d.stronglyconformingtail}
Let $D > 0$.  A strategy $S_+ \in \mc{S}$ for Maxine in $\game(\e,\lambda,v)$ is said to be {\em strongly conforming above fortune $D$} if it meets the demand set in the above specification of  the strategy that is  strongly conforming for Maxine whenever $(\mu,w) \in (0,\infty) \times \openmac$ satisfies $\mu \geq D$. 
\end{definition}
\begin{lemma}\label{l.finitefinish}
Let $D > 0$. 
Let $\splusinfinity \in  \mc{S}_+$ be  strongly conforming above fortune $D$. 
In the game $\game(\e,\lambda,v)$, $F$ is finite almost surely under the strategy pair 
$\big(  \minaconf , \splusinfinity  \big)$. 
\end{lemma}
{\bf Proof.} Consider the gameplay governed by the strategy pair  $\big(  \minaconf , \splusinfinity  \big)$. 

An index $i \in \nwithoutzero$ is called {\em Mina} if the game is unfinished at the start of the $i$\textsuperscript{th} turn; at the start of the turn, $\stateofplay = (\mu,w) \in [0,\infty] \times \openmac$, where the fortune $\mu$ is less than $D$;
a move is selected to take place at the $i$\textsuperscript{th} turn; and Mina wins the right to make the resulting move.

An index  $i \in \nwithoutzero$ is called {\em high} if the game is unfinished at the start of the $i$\textsuperscript{th} turn; and, at the start of this turn, $\stateofplay = (\mu,w) \in [0,\infty] \times \openmac$, where $\mu \geq D$.

Set
 $$
 c \, = \,  \inf \, \Big\{ \, \mu^{-1} \stake(\e,\mu,w) : \mu \in [0,D) , w \in \openmac \, \Big\} \, .
$$
By Corollary~\ref{c.smallstakes}, $c > 0$. We {\em claim} that $c/(1+c)$ is a lower bound on the probability that Mina wins any given turn at whose start the fortune is less than $D$ and at which a move is selected to take place. Indeed, if $\stateofplay = (\mu,w)$
for $\mu < D$, then Mina will stake $\stake(\e,\mu,w)$ and Maxine will stake at most $\mu$, so that Mina's win probability given a move is at least $\tfrac{\stake(\e,\mu,w)}{\mu + \stake(\e,\mu,w)} \geq \tfrac{1}{c^{-1}+1}$; whence the claim. 

 Let $d$ denote the number of vertices on the longest path in $\openmac$. Since $\minaconf$ is strongly conforming, a string of consecutive Mina indices may have at most $d$ elements,  
 because the game will finish at or before the move corresponding to the final index in the string.

Suppose that $F \geq j$ for $j \in \nwithoutzero$. It follows from the above claim that, whatever the status of turn indices on $\llbracket 0, j-1 \rrbracket$, there is conditional probability at least $\e c/(c+1)$ that $j$ is Mina or high or $F=j$. In any block $\llbracket (i-1)d+1,id  \rrbracket$,  and given that $F \geq (i-1)d +1$, there is thus conditional probability at least $\big(\e c/(c+1)\big)^d$ that one of the following occurs: every index in the block is Mina; at least one index in the block is high; or $F \leq id$. By the preceding paragraph, the first alternative implies the third. And if any index is high, $F < \infty$, because both players conform from a turn with a high index and this ensures a finite finish as we noted in the paragraph preceding Definition~\ref{d.stronglyconformingtail}. 
(Note that Definition~\ref{d.stronglyconforming} treats opponent bankruptcy in a way that ensures a finite finish time in the cases where the fortune reaches zero or infinity.)
Thus, $F < \infty$ almost surely. \qed

We extend the truncation notation $\cdot \rfloor_n$ in Definition~\ref{d.thetan} to strategies.
\begin{definition}
Let $S \in \mc{S}_-$. The horizon-$n$ truncation of $S$, denoted by $S \rfloor_n$, is the element in $\mc{S}_-(n)$ formed by restricting $S$ to the first $n$ turns of $\game(\e,\lambda,v)$. We may equally use this notation when $S \in \mc{S}_+$.
\end{definition}

\begin{lemma}\label{l.costone}
Let $D > 0$.
\begin{enumerate}
\item
Let $\sminusinfinity \in  \mc{S}_-$.
 Then
$$
 \ M \big(  \sminusinfinity ,  \maxineconf \big)    \, \geq \,  \limsup_n M \big(  \sminusinfinity \rfloor_n  ,  \maxineconfn \big)  \, .
$$
\item Now let $\splusinfinity \in  \mc{S}_+$ be  strongly conforming above fortune $D$.  We have that
$$
  M \big(  \minaconf , \splusinfinity  \big)  \, \leq \, \liminf_n M \big(  \minaconf \rfloor_n  , \splusinfinity \rfloor_n \big)  \, .
$$
\end{enumerate}
\end{lemma}
{\bf Proof: (1).}  Note that
\begin{equation}\label{e.etwom}
M \big(     \sminusinfinity ,\maxineconf \big) \, = \,  \E_{(\sminusinfinity ,\maxineconf)} \big[ \pay \cdot {\bf 1}_{F = \infty}  \big]  +   \E_{(\sminusinfinity ,\maxineconf)} \big[ \pay \cdot {\bf 1}_{F < \infty}  \big] \, .
\end{equation}
Recall that the rules of the infinite-horizon game declare that, should the game be unfinished in any finite time---that is, should $F$ be equal to $\infty$---then the terminal payment $\pay$ equals one. We invoke this rule now, to find that the first right-hand term in~(\ref{e.etwom}) equals $\PP(F = \infty)$. We may also speak of the finish time $F$ in the finite-horizon game $\game_n(\e,\lambda,v)$; when we do so, $F$ will be at most $n$ for gameplays that finish by the counter leaving open play, and we will adopt the convention of Definition~\ref{d.thetan}
by recording the event of an unfinished game in the form $F = n+1$.
Note further that $M \big(   \sminusinfinity\rfloor_n ,   \maxineconf \rfloor_n    \big) = \zeta_1(n) + \zeta_2(n)$, where $\zeta_1(n) = \E_{(\sminusinfinity\rfloor_n ,   \maxineconf \rfloor_n )} \big[ \pay \cdot {\bf 1}_{F = n+1}  \big]$ and   $\zeta_2(n) = \E_{(\sminusinfinity\rfloor_n ,   \maxineconf \rfloor_n )} \big[ \pay \cdot {\bf 1}_{F \leq n}  \big]$. Since the infinite-horizon gameplay never finishes precisely when all of its truncations are also unfinished, we find that $\limsup_n \zeta_1(n) \leq \PP(F= \infty)$. We also have that $\zeta_2(n)$ increases in $n$ to the limiting value $ \E_{(\sminusinfinity ,   \maxineconf  )} \big[ \pay \cdot {\bf 1}_{F < \infty}  \big]$, which is the second right-hand term in~(\ref{e.etwom}). Thus do we prove Lemma~\ref{l.costone}(1).

{\bf (2).}
Since $F < \infty$ under $(\minaconf ,  \splusinfinity)$ by Lemma~\ref{l.finitefinish}, we have that
\begin{equation}\label{e.payfinite}
M \big(    \minaconf ,  \splusinfinity \big) 
=
 \E_{(\minaconf ,  \splusinfinity)} \big[ \pay \cdot {\bf 1}_{F < \infty}  \big]  \, .
\end{equation}
  Note that
  $$
   M \big(   \minaconf \rfloor_n ,  \splusinfinity\rfloor_n   \big) =  \E_{(\minaconf \rfloor_n ,  \splusinfinity\rfloor_n )} \big[ \pay \cdot {\bf 1}_{F = n+1}  \big] + \E_{(\minaconf \rfloor_n ,  \splusinfinity\rfloor_n )} \big[ \pay \cdot {\bf 1}_{F \leq n}  \big] \, . 
  $$
  The former right-hand term is non-negative, and the latter increases to 
 $\E_{(\minaconf ,  \splusinfinity)} \big[ \pay \cdot {\bf 1}_{F < \infty}  \big]$. From~(\ref{e.payfinite}), we find then that
 $$
 \lim\inf_n M \big(   \minaconf \rfloor_n ,  \splusinfinity\rfloor_n   \big) \, \geq \, 
M \big(    \minaconf ,  \splusinfinity \big) \, ,
 $$
 as we sought to show. \qed
 
An asymmetry under good fortune was noted between Mina and Maxine after Proposition~\ref{p.threelambda}. 
It is in fact this asymmetry that has led us to treat the two players differently in the preceding proof. While Maxine's strategy was circumscribed at high fortunes by a definition in Lemma~\ref{l.costone}(2), it was the `Mina pays one' rule for unfinished games that enabled the proof of Lemma~\ref{l.costone}(1). Indeed, with a rule of the form `Mina pays less than one', we cannot hope to obtain Lemma~\ref{l.costone}(1) by circumscribing Mina's strategy at low fortunes similarly as we did Maxine's: see Section~\ref{s.unfinishedgames}.

\begin{definition}\label{d.stronglyconformingvariant}
Let $S_- \in \mc{S}_-$, $S_+ \in \mc{S}_+$ and $D \in (0,\infty)$. Let $S_+(D)$ denote the element of $\mc{S}_+$ that is strongly conforming above fortune $D$ for which stakes and move nominations
offered when $\stateofplay = (\mu,w) \in (0,\infty) \times \openmac$ for $\mu \in [0,D)$ are governed by $S_+$. 
\end{definition}

 Recall that $\lambda = \lambda_0$ is the initial fortune in $\game(\e,\lambda,v)$.
 
\begin{lemma}\label{l.costtwo}
Let $S_- \in  \mc{S}_-$ be a strategy for Mina that is non-conforming against $\maxineconf$. 
\begin{enumerate}
\item Then $\liminf_n   M \big( S_- \rfloor_n , \maxineconfn \big) \, > \, h(\lambda,v)$.
\end{enumerate}
Let  $S_+ \in \mc{S}_+$ be a strategy for Maxine that is non-conforming against $\minaconf$.  
\begin{enumerate}
\setcounter{enumi}{1}
\item   There exists $\delta > 0$ such that $D > \lambda$ implies that
$$
\limsup_n   M \big( \minaconfn , S_+(D) \rfloor_n  \big) \, < \, h(\lambda,v) - \delta \, ,
$$
where $S_+(D)$ is specified in Definition~\ref{d.stronglyconformingvariant}.\end{enumerate}
\end{lemma}

{\bf Proof: (1).} This proof is omitted because it is in essence a slightly simplified version of the proof of the second part.

{\bf (2).} 
We extend the notation of Definition~\ref{d.stronglyconformingtail} so that $S_+(\infty)$ denotes $S_+$. 
Consider gameplay under the \jointstrategy $\big(  \minaconf, S_+(\infty) \big)$.
 Let the {\em deviation} set $\mathsf{D}$ denote the set of indices $i \in \nwithoutzero$
such that there is positive probability that Maxine does not conform at the $i$\textsuperscript{th} turn in the sense of Definition~\ref{d.conforming}. 
Since $S_+(\infty)$ is by assumption non-conforming,  $\mathsf{D}$ is non-empty. Let $j \in \nwithoutzero$ denote the least element in $\mathsf{D}$.

To derive Lemma~\ref{l.costtwo}(2), we will now argue that it suffices to prove the bound
\begin{equation}\label{e.itsuffices}
\limsup_n   M \big( \minaconfn , S^j_+(\infty) \rfloor_n  \big) \, < \, h(\lambda,v)  \, ,
\end{equation}
where, for $D \in (0,\infty]$, we let $S_+^j(D)$ denote the element of $\mc{S}_+$
obtained by modifying $S_+(D)$ so that play by Maxine at every turn with index at least $j+1$ is conforming. We now make two claims, which will prove that deriving the bound~(\ref{e.itsuffices}) is indeed sufficient for our purpose. First,
\begin{equation}\label{e.claimone}
M \big( \minaconfn , S^j_+(D) \rfloor_n  \big) \, \geq \,  M \big( \minaconfn , S_+(D) \rfloor_n  \big) \, \, \, \textrm{for $D \in (0,\infty]$} \, .
\end{equation}
 Second,
 \begin{equation}\label{e.claimtwo}
   M \big( \minaconfn , S^j_+(D) \rfloor_n  \big) \, = \,  M \big( \minaconfn , S^j_+(\infty) \rfloor_n  \big) 
 \, \, \, \textrm{for $D > \lambda$} \, .
\end{equation}
Let $D > \lambda$. Applying~(\ref{e.claimone}),~(\ref{e.claimtwo}) and~(\ref{e.itsuffices}), we obtain 
 Lemma~\ref{l.costtwo}(2). 
 
 We now prove~(\ref{e.claimone}). Theorem~\ref{t.leisurely.theorem}(2) shows that Maxine's play under $S_+^j(D)$ from the  $(j+1)$\textsuperscript{st}  turn is governed by a strategy that, when paired with~$\minaconfn$, forms a Nash equilibrium in $\game_{n-j}(\e,\lambda^*,v^*)$, where $\lambda^*$ and $v^*$ are the values obtained after the possibly non-conforming play in the $j$\textsuperscript{th}   step. We see then that, in replacing $S_+(D)$ with $S_+^j(D)$, Maxine will not see a fall in mean terminal payment: that is,~(\ref{e.claimone}) holds.
 
  And now the proof of~(\ref{e.claimtwo}).  
 During the first \textcolor{blue}{$j-1$} turns, both players conform under the strategy pair $\big(  \minaconf, S^j_+(\infty) \big)$, so that $\lambda_i = \lambda$ for $i \in \llbracket 0,j-1 \rrbracket$. 
 Since $D > \lambda$ (in the bound we seek), the dictates of~$S^j_+(D)$ coincide with those of $S^j_+(\infty)$ for as long as the present fortune equals $\lambda$.
 Hence, gameplay coincides until the end of the $j$\textsuperscript{th} turn under the strategy pairs  $\big(  \minaconf, S^j_+(\infty) \big)$ and $\big(\minaconf, S^j_+(D) \big)$.
 But beyond the end of the $j$\textsuperscript{th} turn, Maxine conforms under either $S^j_+(D)$ or $S^j_+(\infty)$; as does Mina under $\minaconf$. So in fact the two gameplays coincide throughout the lifetime of the game. This proves~(\ref{e.claimtwo}).

 The two claims justified, it remains, in order to obtain Lemma~\ref{l.costtwo}(2), to prove~(\ref{e.itsuffices}). Set $S$ equal to~$S_+^j(\infty)$. Our task is to show that 
$\limsup_n   M \big( \minaconfn , S \rfloor_n  \big)  <  h(\lambda,v)$. We thus consider gameplay under  $\big( \minaconfn , S \rfloor_n  \big)$.
The deviation set $\mathsf{D}$ associated to this gameplay has a unique element, $j$. Let~$E_j$ denote the event that, at the $j$\textsuperscript{th} turn, Maxine's stake does not conform. 
Let $F_j$ denote the counterpart event where `move nomination' replaces `stake'. Note that $\PP(E_j)+\PP(F_j)>0$. Definition~\ref{d.conforming} 
indicates the form of these two events:
on $E_j$,
 Maxine's stake at the $j$\textsuperscript{th} turn is not equal to $\lambda_{j-1} \cdot \stake(\e,\lambda_{j-1},X_{j-1})$; 
while $F_j$ occurs when Maxine nominates a move at this turn that lies outside of the singleton set~$\mc{V}_+(X_{j-1})$.

If the strategy $S$ is modified at the  $j$\textsuperscript{th} turn so that Maxine conforms against $\minaconfmacro$ at this turn, a strategy in $\game(\e,\lambda,v)$ that is conforming against $\minaconfmacro$ results, in the sense of Definition~\ref{d.conforming}. Denote this modified strategy
by $\maxineconf$ and its index by $\theta \in \Theta$.
Note then that 
\begin{equation}\label{e.twom}
 M \big(  \minaconfn ,  S \rfloor_n \big) = \sum_{i=1}^4 \gamma_i(A,n)  \, \, \, \, \textrm{and} \, \, \, \, M \big( \minaconfn , \maxineconfn   \big) = \sum_{i=1}^4 \gamma_i(B,n) \, ,
\end{equation}
where 
$\gamma_i(A,n)$, $i \in \intint{4}$, are the expected values of 
$\pay \cdot {\bf 1}_{E_j \cap F_j^c}$, $\pay \cdot {\bf 1}_{F_j \cap E_j^c}$, $\pay \cdot {\bf 1}_{E_j \cap F_j}$  and $\pay \cdot {\bf 1}_{E_j^c \cap F_j^c}$ 
under the \jointstrategy  $\big( \minaconfn , S \rfloor_n  \big)$, and 
$\gamma_i(B,n)$, $i \in \intint{4}$,
are the counterpart quantities under $\big(  \minaconfn  , \maxineconfn \big)$.

We claim that
\begin{equation}\label{e.threeclaims}
 \gamma_i(A,n) - \gamma_i(B,n) \, \, \, \textrm{is negative and independent of $n \geq j$ for $i \in \intint{3}$} \, ,
\end{equation}
and also that $\gamma_4(A,n) = \gamma_4(B,n)$. The latter statement has the simplest proof: the discrepancy between the two concerned strategies for Maxine is manifest only at the  $j$\textsuperscript{th} turn and this discrepancy does not affect gameplay when $E_j^c \cap F_j^c$ occurs. We now prove in turn the three assertions made in~(\ref{e.threeclaims}). 
First, note that
$$
 \gamma_1(A,n) - \gamma_1(B,n) 
 =    \sum_{w \in \openmac}  \PP \big( X_\theta(j-1) = w , E_j \cap F_j^c \big) \int \alpha_n(x) \, {\rm d}\mu_{j,w}(x) \, ,
$$
where recall that $\theta \in \Theta$ is the index of Maxine's strongly conforming strategy $\maxineconf$; the quantity $\alpha_n(x)$ is specified by
\begin{equation}\label{e.alphax}
\alpha_n(x) \, = \,  \val_{n+1-j}(\e,\lambda_{j-1},w,x,\stakemacro) - \val_{n+1-j}(\e,\lambda_{j-1},w,\lambda_{j-1} \stakemacro,\stakemacro) \, ,
\end{equation} 

where $\val_n(\epsilon,\lambda,v,a,b)$ is the value of the constrained version of the finite-horizon $\game_n(\epsilon,\lambda,v)$, in which Maxine must stake $a$ and Mina must stake $b$ at the first turn (a value which, as we will note momentarily, is well defined); and the law 
 $\mu_{j,w}$ is the conditional distribution of the stake offered by Maxine at the $j$\textsuperscript{th} turn given that  $E_j$, $F_j^c$ and $X_\theta(j) = w$ occur. 
 
By Theorem~\ref{t.leisurely.theorem}, $\game_n(\epsilon,\lambda,v)$ has game value $h(\lambda,v)$ for $(n,\lambda,v) \in \N \times (0,\infty)\times \openmac$.  Thus $\val_n(\epsilon,\lambda,v,a,b)$ is well defined for $n \in \N$.  Moreover, we have $\val_n(\epsilon,\lambda,v,a,b) = \operatorname{V}(\epsilon,\lambda,v,a,b)$, where $\operatorname{V}$ is defined in~\eqref{constrained value function}. Hence, we can rewrite~\eqref{e.alphax} as \begin{equation*}
    \alpha_n(x) \, = \,  \operatorname{V}(\e,\lambda_{j-1},w,x,\stakemacro) - \operatorname{V}(\e,\lambda_{j-1},w,\lambda_{j-1} \stakemacro,\stakemacro) \, .
\end{equation*}
 
The quantity $\stakemacro$ equals $\stake(\e,\lambda_{j-1},w)$ from~(\ref{e.stake}), so that $(\lambda_{j-1} \stakemacro,\stakemacro)$
is the strict global minimax for $(0,\infty)^2 \lora [0,1]: (a,b) \mapsto \operatorname{V}(\e,\lambda_{j-1},w,a,b)$ by Proposition~\ref{p.globalconstrainedminimax}. Note that the law $\mu_{j,w}$ assigns zero mass to the point~$\stakemacro$. Note then that, for $x \in (0,\infty)$, $x \neq \lambda_{j-1} \stakemacro$, 
\begin{equation}\label{e.alphaxnindep}
\textrm{the quantity $\alpha(x) = \alpha_n(x)$ is negative, and independent of $n \geq j$} \, .
\end{equation}
Thus we obtain~(\ref{e.threeclaims}) for $i=1$.

Note next that, by Lemma~\ref{l.hnconv},  $\gamma_2(A,n) - \gamma_2(B,n)$  equals
$$
  \sum_{w \in \openmac}  \PP \big( X_\theta(j-1) = w , F_j \cap E_j^c \big)  \e \lambda_{j-1}(1+\lambda_{j-1})^{-1} \int  \big( h(\lambda_{j-1},z) - h(\lambda_{j-1},w_+)
 \big){\rm d}\zeta_{j,\lambda_j,w}(z) \, ,
$$ 
where  $\zeta_{j,\lambda_j,w}(\cdot)$ is the discrete law charging neighbours of $w$ that is given by the conditional distribution of  Maxine's move nomination at the $j$\textsuperscript{th} turn given that  $F_j$, $E_j^c$ and $X_\theta(j-1) = w$ occur and the fortune\footnote{This fortune is $\lambda_j$, not $\lambda_{j-1}$, because players learn the revised fortune before nominating a move: see the update to $\stateofplay$ in the {\em first step} in Section~\ref{s.specifying}.} is $\lambda_j$. 
 Indeed, if a {\em move takes place} at the turn  that {\em Maxine wins} at which {\em she nominates $z$}, she suffers a conditional mean change of $h(\lambda,z) - h(\lambda_{j-1},w_+)$  in terminal payment; the respective probabilities attached to the italicized terms are $\e$, $\lambda_{j-1}(1 + \lambda_{j-1})^{-1}$ and $\zeta_{j,\lambda_j,w}(z)$.
The law $\zeta_{j,\lambda_j,w}(\cdot)$ does not charge the vertex $w_+$ because Maxine's $j$\textsuperscript{th}  move nomination is non-conforming when $F_j$ occurs. 
 Thus
  $h(\lambda_{j-1},z) - h(\lambda_{j-1},w_+)$
 is seen to be negative (and independent of $n$) when $z$ is in the support of $\zeta_{j,\lambda_j,w}(\cdot)$, 
 by  Corollary~\ref{c.uniqueplay}.
Thus,~(\ref{e.threeclaims}) holds for $i=2$.

Next we consider $i=3$. In this case, Maxine makes a mistake at the $j$\textsuperscript{th} turn both in her stake and in her move nomination. Let 
$\psi_{j,w}$ denote the conditional distribution of the stake and move nomination pair offered by Maxine at the $j$\textsuperscript{th} turn given that  $E_j \cap F_j$ and $X_\theta(j) = w$ occur. 
(We could write $\psi_{j,\lambda_{j-1},\lambda_j,w}$: the stake is decided knowing $\lambda_{j-1}$; the move is nominated knowing $\lambda_j$.)
The respective marginals of $\psi_{j,w}$ charge neither the point $\lambda_{j-1} \stakemacro$ nor the vertex $w_+$. 
We have that 
$$
 \gamma_3(A,n) - \gamma_3(B,n) 
  \, = \,   \sum_{w \in \openmac}  \PP \big( X_\theta(j-1) = w , E_j \cap F_j \big)  \int  \big(  \alpha(x) + \beta_n(x,z)  \big) {\rm d}\psi_{j,w}(x,z)  \, ,
 $$
where
$\alpha(x)$ is  specified in~(\ref{e.alphax}),
and 
$$
\beta_n(x,z) \, = \, \e \tfrac{x}{x+ \stakemacro} \cdot  \Big( \val_{n-j} \big(\tfrac{\lambda_{j-1}-x}{1 - \stakemacro},z \big) - \val_{n-j} \big(\tfrac{\lambda_{j-1}-x}{1 - \stakemacro},w_+ \big) \Big) \, .
$$
As above, $\alpha(x)$ is the conditional mean change in payment resulting from Maxine's faulty decision to stake $x$ at the  $j$\textsuperscript{th} turn.
The corresponding change caused by her move nomination of~$z$, given the stake pair $(\lambda_{j-1} \stakemacro,x)$, is 
$\beta_n(x,z)$: indeed, with probability $\e$, a move takes place at the   $j$\textsuperscript{th} turn; with probability  $\tfrac{x}{x+\stakemacro}$, Maxine wins the right to make the resulting move; and, if she does so, she incurs a change in conditional mean payment of  $\val_{n-j-1} \big(\tfrac{\lambda_{j-1}-x}{1 - \stakemacro},z \big) - \val_{n-1-j}\big(\tfrac{\lambda_{j-1}-x}{1 - \stakemacro},w_+ \big)$ by nominating~$z$ instead of the optimal choice $w_+$. The quantity $\beta(z,a) = \beta_n(x,z)$ is negative, and it is independent of $n \geq j$ in view of Theorem~\ref{t.leisurely.theorem}(1). 
Recalling~(\ref{e.alphaxnindep}), we obtain (\ref{e.threeclaims}) for $i=3$.


With ~(\ref{e.threeclaims}),  and $\gamma_4(A,n) = \gamma_4(B,n)$, we return to~(\ref{e.twom}) to learn that   
$$
M \Big( S \rfloor_n , \maxineconfn \Big) - M \Big(  \minaconfn  , \maxineconfn \Big) \, = \, \sum_{i=1}^3 \big( \gamma_i(A) - \gamma_i(B) \big)  
$$ 
is negative and independent of $n \geq j$. 
Since  the \jointstrategy $\big(  \minaconfn  , \maxineconfn \big)$ is a Nash equilibrium in $\game_n(\e,\lambda,v)$ by Theorem~\ref{t.leisurely.theorem}(2),
  $M \big(  \minaconfn  , \maxineconfn \big)$ equals $h(\lambda,v)$ by Theorem~\ref{t.leisurely.theorem}(1) and Lemma~\ref{l.hnconv}. Recalling the shorthand $S = S_+^j(\infty)$, we obtain~(\ref{e.itsuffices}) and thus complete the proof of Lemma~\ref{l.costtwo}(2). \qed
 
 Recall Definition~\ref{d.stronglyconformingtail}. 
\begin{lemma}\label{l.deta}
For $\eta > 0$, any sufficiently high $D > 0$ is such that
$$
M \big(      \minaconf , \splusinfinity \big) \leq 
M \big(       \minaconf , \splusinfinity(D)  \big)  + \eta \, .
$$
\end{lemma}
{\bf Proof.} 
Proposition~\ref{p.rootrewardinfinitybias} implies that   $\lim_{\lambda \nearrow \infty} h(\lambda,v) = 1$ for $v \in \openmac$. Thus any high enough $D > 0$ satisfies 
\begin{equation}\label{e.deta}
h(D,v) \geq 1 - \eta \, \, \,  \textrm{for}  \, \, \, v \in \openmac \, .
\end{equation}
Let $\maxinehighfortuneinfinity$ denote the event that Maxine's fortune is at least $D$ before some turn in $\game(\e,\lambda,v)$.
Note that
$$
M \big(      \minaconf , \splusinfinity \big)  = \E_{(\minaconf , \splusinfinity)} \big[ \pay \cdot {\bf 1}_{\maxinehighfortuneinfinity} \big]+ \E_{(\minaconf , \splusinfinity)} \big[  \pay \cdot {\bf 1}_{\neg \, \maxinehighfortuneinfinity}\big] \, .
$$
The first right-hand term is at most $\rho := \PP_{(\minaconf , \splusinfinity)} \big( \maxinehighfortuneinfinity \big)$, and the second coincides with its counterpart in the next equality:
\begin{multline*}
 M \big(  \minaconf , \splusinfinity(D) \big) \\
  = \E_{( \minaconf , \splusinfinity(D))} \big[ \pay \cdot {\bf 1}_{\maxinehighfortuneinfinity} \big]+ \E_{( \minaconf , \splusinfinity(D))} \big[  \pay \cdot {\bf 1}_{\neg \, \maxinehighfortuneinfinity}\big] \, .
\end{multline*}
The first right-hand term in the last display satisfies
$$
\rho  \cdot \E_{(\minaconf , \splusinfinity(D))} \big[ \pay \big\vert \maxinehighfortuneinfinity \big] \geq  \rho(1 - \eta) \, ,
$$
where the inequality is due to Corollary~\ref{c.hinc}   and~(\ref{e.deta}).  Since $\rho \leq 1$, the lemma has been proved. \qed

{\bf Proof of Proposition~\ref{p.mistakecost}.} The first assertion follows from Lemma~\ref{l.costone}(1) and Lemma~\ref{l.costtwo}(1).  
To prove the second, let $D > \lambda$.
By  Lemma~\ref{l.costone}(2) with $S_+$ there taken equal to the present $S_+(D)$, and 
Lemma~\ref{l.costtwo}(2),
we find that $M \big( \minaconf , S_+(D)   \big)  <  h(\lambda,v) - \delta$. By increasing $D > 0$ if need be, we may apply Lemma~\ref{l.deta} with $\eta = \delta/2$ to obtain $M \big( \minaconf , S_+   \big)  <  h(\lambda,v) - \delta/2$. Whence the proposition's second assertion. \qed

{\bf Proofs of Theorems~\ref{t.leisurely} and~\ref{t.nashform}.} 
Now we assume that $\e_0 = \vert V \vert^{-2\vert \boundarymac \vert -1} $, and $\e \in (0,\e_0)$. By Corollary~\ref{c.estimateofe0}, we know that for $\e \in (0,\e_0)$, Theorem~\ref{t.leisurely.theorem}(1) and~(2) are both true. 

Continuing to use our notation for strongly conforming strategies, we claim that
\begin{equation}\label{e.hlambdav}
M \big( \minaconfmacro, \maxineconfmacro \big) \, = \, h(\lambda,v) \, . 
\end{equation}
Indeed, the counter evolution under gameplay governed by the \jointstrategy  $\big( \minaconfmacro, \maxineconfmacro \big)$
takes the form $X_\theta: \llbracket 0, F_\theta \rrbracket \lora V$, $X_\theta(0) = v$. The process $\llbracket 0, F_\theta \rrbracket \lora [0,1]: i \mapsto h\big(\lambda,X_\theta(i) \big)$
is a martingale. Moreover, and as we noted after Proposition~\ref{p.mistakecost},~(\ref{e.ftail})  implies that $F_\theta$ is finite almost surely.
Thus the mean terminal payment as a function of the starting location satisfies the system~(\ref{e.h}). Hence, we obtain~(\ref{e.hlambdav}). 

Proposition~\ref{p.mistakecost} thus demonstrates that  any \jointstrategy whose components are strongly conforming 
is a Nash equilibrium; we thus obtain the assertions in Theorem~\ref{t.leisurely}(1) that  the value of $\game(\e,\lambda,v)$ equals $h(\lambda,v)$ and in Theorem~\ref{t.leisurely}(2) that a Nash equilibrium exists in $\game(\e,\lambda,v)$.

Suppose now that $(S_-,S_+)$ is a Nash equilibrium 
at least one of whose components is non-conforming against the opposing element, in the sense that, under gameplay governed by $(S_-,S_+)$, there is a positive probability that Maxine or Mina does not conform at some move.
Let $\ell_+ \in \nwithoutzero$ be the minimum turn index~$i$ at which there is a positive probability that Maxine does not conform at the $i$\textsuperscript{th} turn; we take $\ell_+ = \infty$
if all such probabilities are zero.  Let $\ell_-$ denote the counterpart quantity for Mina. Suppose that $\ell_- \leq \ell_+$; we have then that $\ell_-$ is finite. 
Let $\maxineconfmacro$ denote the strategy for Maxine that coincides with $S_+$ at turns with index at most $\ell_- - 1$, and that adheres to Maxine's strongly conforming   strategy at turns with index at least $\ell_-$.  (Our notation would be inconsistent if $\maxineconfmacro$ could fail to be strongly conforming. And in fact this difficulty may occur. But this is due merely to non-conforming choices that Maxine may make in states of play that are almost surely inaccessible under the gameplay governed by $(S_-,\maxineconfmacro)$. We may harmlessly correct such deviations of Maxine, so that $\maxineconfmacro$ may indeed be supposed to be strongly conforming.) 
Mina's play under  the \jointstrategy $\big(S_-,\maxineconfmacro\big)$ is non-conforming,
 because it is with positive probability that she does not conform at the turn with index~$\ell_-$. 
We have then that 
$$
h(\lambda,v) = M(S_-,S_+) \geq M\big(S_-,\maxineconfmacro\big) > M\big( \minaconfmacro,\maxineconfmacro\big) = h(\lambda,v) \, , 
$$
where the first equality, and the first inequality, are due to $(S_-,S_+)$ being a Nash equilibrium; the strict inequality is by appeal to Proposition~\ref{p.mistakecost} in view of the strong conformity of $\maxineconfmacro$ and the just noted non-conformity of the component $S_-$ in the \jointstrategy~$\big(S_-,\maxineconfmacro\big)$; and the latter equality is due to $\big( \minaconfmacro,\maxineconfmacro\big)$  being a Nash equilibrium (a fact shown earlier in this proof).
The displayed contradiction shows that $(S_-,S_+)$ is not a Nash equilibrium, at least under the assumption that $\ell_- \leq \ell_+$; but there was no loss of generality in making this assumption, because the opposing case has a similar proof. 
 In this way, we complete the proof of Theorem~\ref{t.leisurely}(2), and show that $\e_0$ may be taken to be $\vert V \vert^{-2 \vert \boundarymac \vert -1} $. And we obtain Theorem~\ref{t.nashform}.\qed

{\bf Proof of Corollary~\ref{c.gameplay}.} Elements of the pair $\big( \minaconfmacro,\maxineconfmacro \big)$ considered in the preceding proofs are strongly conforming, and the index $\theta \in \Theta$ of $\minaconfmacro$ is arbitrary. The gameplay process~$X_\theta$ results from the choice of $\theta \in \Theta$, and thus all of the claimed processes occur as gameplay processes governed by Nash equilibria in $\game(\e,\lambda,v)$. 

In order to show that gameplay can be no other process, we develop Definition~\ref{d.conforming}(2). Let  $(S_-,S_+) \in \mc{S}_- \times \mc{S}_+$, $\theta \in \Theta$ and $i \in \nwithoutzero$.   Mina conforms with index $\theta$ at the $i$\textsuperscript{th} turn under the gameplay governed by $(S_-,S_+)$ if she stakes   $\stake \big( \e,\lambda_{i-1},X(i-1) \big)$ 
and nominates the move $\theta \big( \lambda_{i-1}, X(i-1),i-1 \big)$ at this turn when this given \jointstrategy is adopted. The strategy~$S_-$ is said to be conforming with index $\theta$ against $S_+$ if Mina almost surely conforms with index $\theta$ at every turn. (To establish the relation to Definition~\ref{d.conforming}, note that $S_-$ conforms against~$S_+$ if and only if~$S_-$ conforms with index $\theta$ against $S_+$ for some $\theta \in \Theta$.)

Theorem~\ref{t.nashform} asserts that any Nash equilibrium in  $\game(\e,\lambda,v)$ is a conforming \jointstrategyperiod Maxine's hand is thus forced at every turn. Mina's choice at a given turn is restricted to the selection of an element of $\mc{V}_-(w)$ for her move nomination, where $w \in \openmac$ is the present counter location. By Definition~\ref{d.strategy}, she may only take account of the values of $\stateofplay$ and the turn index in making her choice. This forces her play to be conforming with index $\theta$ against Maxine's chosen strategy, for some $\theta \in \Theta$. This completes the proof of  Corollary~\ref{c.gameplay}. \qed

\section{Directions and open problems}\label{s.directions}

We offer an overview of some prospects for developments of the concepts and proofs in the article in several subsections that begin with more  specific and technical aspects and end with broader themes.

\subsection{Strategy spaces and reset rules}\label{s.resetrules}

We have not hesitated to posit limitations on game design and strategy spaces in the interests of proof simplicity provided that such assumptions change nothing essential about the anticipated game values and equilibria. 
These assumptions can be reviewed. The pair of reset rules recorded in Section~\ref{s.specifying} provides easy means of settling questions in proofs that arise from joint-zero or joint-total stakes at a given turn. But in fact players are anyway tempted away from zero or total stakes, as comments under the pair of slogans in Subsection~\ref{s.sketchesandlessons} indicate and further rigorous argument may be expected to imply. 

Legal strategies in Definition~\ref{d.strategy} permit only the present values of Maxine's fortune and counter location, alongside the turn index, to inform players' choices at a given turn. These definitions could be broadened so that the game history is available, though we do not anticipate that the basic structure of value and equilibria to shift in response to this broadening. Players have been permitted to randomize stakes, but not moves, and a similar comment may be made in this regard. Choices that are mixed for both move and stake would render the mixed strategy spaces convex, which may be an attractive feature should fixed point theorems be brought to bear in efforts to prove the existence of equilibria in variants of the game: we turn to this topic in the next paragraph but one.

\subsection{The payment when the game is unfinished}\label{s.unfinishedgames}

In Section~\ref{s.specifying}, we specified that Mina will pay one unit when $\game(\e,\lambda,v)$ fails to finish. This choice is needed to permit our proof of Lemma~\ref{l.costone}(1). 
Indeed, if $\pay$ is a given value less than one when the finish time is infinite, Mina can sometimes oppose conforming play from Maxine and achieve a lower mean payment than she would obtain were she to conform. To see this, consider the graph formed by attaching a leaf $z$ along an edge to the vertex $n-1$ in the root-reward line graph $\big( \llbracket 0,n \rrbracket, \sim,{\bf 1}_n \big)$.
This root-reward tree has root $n$, and two further leafs, $0$ and $z$. The vertex $n-2$ has journey data $\big\{ (2,1),(n-1,1) \big\}$. The game played from~$n-2$ is a contest of two rounds: if Maxine is to win, the counter must first reach $n-1$; and then $n$. Suppose that $\lambda \gg 1$. In the first round, little is at stake; but the second is more contested. Indeed, if both players conform, then Corollary~\ref{c.highlambda} shows that the stake offered in $\game(1,\lambda,n-2)$, namely $\stake(1,\lambda,n-2)$, takes the form $\lambda^{3-n}\big(1+o(1)\big)$ as $\lambda \to \infty$. Suppose that Mina modifies conforming play only at vertex $n-2$, where she stakes one-half (and plays left). From $n-2$, the counter typically moves left; then it returns to $n-2$, usually after one further turn, with the value of $\lambda$ almost doubling since the last visit to $n-2$. Except with probability of order $\lambda^{3-n}$, the game will never end. Conditionally on this non-finishing event, the gameplay empirical process $h(\lambda_i,X_i)$ will converge to one almost surely. But Mina's terminal payment will be a given constant less than one, by fiat. Thus we see that the rule for unfinished games, even when it makes the seemingly mundane demand that the terminal payment be a given constant on $[0,1]$, is not merely a technical detail to ensure well-specified play: if the constant is less than one, Maxine must do something to deviate from conforming play to overcome this unending filing of extensions from Mina. From this example, we see that the `Mina pays one' rule should be viewed, alongside the use of root-reward trees and a small move probability~$\e$, as a hypothesis that substantially enables our proofs. It would be of interest to inquire whether, admitting these other assumptions, the rule `Mina pays $p$', for $p \in [0,1)$ given, leads to the biased-infinity game value found under the $p=1$~rule.

\subsection{Abstract existence results for Nash equilibria}

We have constructed Nash equilibria by concrete means rather than by seeking to use abstract results. Von Neumann~\cite{VonNeumann} proved the existence of value in two-person zero-sum games  with finite strategy spaces via his minimax theorem, and continuum strategy spaces have been treated by Glicksberg~\cite{Glicksberg}, provided that the mean payoff is continuous.
 When it exists, the value $\val(1,\lambda,v,a,b)$ of the first-turn-constrained regular game is not in most cases continuous at $(a,b) = (\lambda,1)$: the final paragraph of Section~\ref{s.localglobal} addresses this point (in a special case). 
Existence results for equilibria where the continuity hypothesis is weakened have been obtained by Simon~\cite{Simon87}, and it would be of interest to examine the applicability of such theory to stake-governed games.

\subsection{Prospects for specifying and analysing the Poisson game}

A formal analysis of the Poisson game in Section~\ref{s.poisson} provided (what we hoped to be!) a simple and attractive point of departure for our proofs concerning the leisurely game. One may wish to treat the Poisson game rigorously. Significant conceptual questions must be answered to make rigorous sense of this continuous-time game, however. One could say that each player must adhere to a strategy that is adapted to the strict past history of gameplay. With such instantaneous communication, strategies such as `I'll stake twice what she just staked' can however lead to a folie-\`a-deux. Presumably, they should be banned by suitable constraints on strategy spaces. Simon and Stinchcombe~\cite{SimonStinchcombe} have developed a framework for addressing such problems which may be applicable for stake-governed games.  

For many boundary-payment graphs, it may be that the global saddle hope is realized for the Poisson game at $\lambda$-values outside a finite set at which the Peres-\v{S}uni\'c decomposition  changes. It is conceivable, however, that a player may force the running $\lambda$-value onto the special finite set for a positive measure of times, so that the more complex behaviour apparent in Figure~\ref{f.poisson}({\rm middle}) becomes germane even if the initial value of $\lambda$ is generic. This is highly speculative, but such possibilities should be borne in mind when the Poisson game is analysed.

\subsection{Stake-governed tug-of-war beyond root-reward trees}

The absence of dependence of journey data in Definition~\ref{d.journeydata} on the value of  $\lambda \in (0,\infty)$ attests that the Peres-\v{S}uni\'c decomposition is independent of~$\lambda$ for root-reward trees. This condition is fundamental to our analysis, because it disables the capacity (seen in Subsection~\ref{s.tgraph}) of a player  to bamboozle an opponent by uncertainties over whether she will play a short or long game. Indeed, the global saddle hope appears to be false in many cases for the regular game $\game(1,\lambda,v)$. A challenge is to make sense of the complexities of behaviour apparent in such contour plots as those in Figures~\ref{f.lthree} and~\ref{f.tgraph} and to resolve questions about the existence and structure of Nash equilibria for stake-governed tug-of-war when the global saddle hope fails. 

\subsection{Predictions for the regular game}

In Section~\ref{s.picture}, we saw that the `big picture' prediction of optimal play governed by the stake formula~(\ref{e.stake}) is inaccurate in the regular game (with $\e = 1$) for at least some values of $\lambda$ in graphs as simple as the $T$ graph. But we also saw that the prediction is correct in one case at least, that of the half-ladder $H_n$ for $n \in \nwithoutzero$.
It would be interesting to find further examples where the prediction is valid and to seek to characterize the graphs for which it is.

\subsection{A one-parameter family of games}

In a `Poorman' variant of the Richman games that is analysed in \cite{LLPSU}, the higher staking player at a given turn wins the right to move, with the stake of this player surrendered to a third party. We may vary this game to bring it closer to stake-governed tug-of-war if we instead insist that the stakes of both players are thus surrendered. Call this the Poormen game. 
In $\lambda = 1$ constant-bias tug-of-war, the right to move is allocated according to a fair coin flip, independently of the stakes (which are thus unnecessary). The Poormen game and tug-of-war may be viewed as $\gamma=\infty$ and $\gamma=0$ endpoints of a one-parameter family of games. Indeed, let $\gamma \in [0,\infty]$. If Maxine stakes $a$ and Mina $b$ at a given turn, then, in the Tullock contest with exponent $\gamma$,
Maxine's win probability is $a^\gamma/(a^\gamma + b^\gamma)$; the rules of stake-governed tug-of-war are otherwise in force. 
As we mention in Section~\ref{s.econ}, Klumpp~\cite{Klumpp} has indicated that the condition $\gamma \leq 1/2$
characterises when a premise corresponding to the stake function formula~(\ref{e.stake}) holds. It would be interesting to study further 
this continuum of games.

\subsection{Continuous-game versions and PDE}

Tug-of-war attracted great interest in PDE from its inception because of its capacity to prove properties and develop intuitions about such famously subtle problems as the uniqueness and regularity theory of the infinity Laplacian on domains in Euclidean space. It is very natural to seek to take the same path for the stake-governed version. Mina and Maxine will play on a domain in $\R^d$, moving a counter a given small distance when one or other wins the right to do so. Does the stake formula~(\ref{e.stake}) take a counterpart form in the limit of small step size? Does the alternative stake formula have a counterpart whose denominator is expressed in terms of a limiting stochastic process for continuum gameplay---an $\infty$-Brownian motion? Which PDE does the continuum stake function satisfy?
These questions may also be posed for the case of the $p$-Laplacian, which arises in tug-of-war with interjections of noise at each turn~\cite{PeresSheffield,Lewicka}, because we may naturally specify a noisy version of stake-governed tug-of-war in the Euclidean setting. This setting may be more tractable than the $p=\infty$ case because the concerned objects---functions and stochastic processes---are at least somewhat more regular and, presumably, more straightforward to define.

\subsection{Stake-governed selection games}

Study the same questions---existence of pure Nash equilibria, the value of the stake function, and so on---for the stake-governed versions of random-turn selection games based on monotone Boolean functions, such as iterated majority, AND/OR trees, or critical planar percolation \cite{PSSW07}. Of interest is to investigate the possible interplay between changes in the bias $\lambda$ and phase transition phenomena \cite{OSSS,DCRT}, especially in light of our formula~(\ref{e.altstake}) connecting optimal stake and expected length of the game, and  of the role of low revealment algorithms in the sharpness of phase transitions.

\subsection{Dynamic strategies for political advertising: stake games played in parallel}

Mina wears a red hat and Maxine a blue one. On each of fifty boundary-payment graphs, a counter is placed at some vertex. In each graph, the payment function $f:\boundarymac \lora [0,\infty)$ 
takes the form $f = \kappa {\bf 1}_D$ for some $D \subseteq \boundarymac$ and constant $\kappa \in \nwithoutzero$ determined by the graph. 
 For some graphs, such as CA and~NY, $\kappa$ is around forty, and $\boundarymac \setminus D$ comprises a few isolated sites;
 for others, such as~SD and~WY, $\kappa$ is about five, and it is $D$ that is a small and isolated set; while on such graphs as NC and PA, neither $D$ nor $\boundarymac \setminus D$ is evidently more exposed in a suitable harmonic sense. The sum of $\kappa$-values over the graphs is $538$.
Maxine and Mina have respective budgets of $\lambda$ and one at the outset. At each turn, each must allocate a stake in each graph from her remaining fortune; a move takes place in each graph at each turn, and the right to move  is allocated to a given player in each graph with probability equal to the ratio of that player's stake, and the combined stakes, offered for that graph at that turn. Maxine's objective is to maximize the probability of receiving at least~$270$ units; or, if we exercise the liberty to choose our assumptions, simply to maximize the mean terminal payment.  

\addcontentsline{toc}{section}{References}

\bibliography{stake}

\end{document}